\documentclass[letterpaper]{article}

\usepackage[letterpaper,left=1.5in,right=1.5in,top=1.5in,bottom=1.5in]{geometry}

\newcommand{\myauthor}{Benjamin Antieau and David Gepner}
\newcommand{\mytitle}{Brauer groups and etale cohomology in derived algebraic geometry}
\newcommand{\pdftitle}{\mytitle}

\title{Brauer groups and \'etale cohomology\\ in derived algebraic geometry}
\author{Benjamin Antieau\footnote{This material is based upon work supported in part by the NSF under Grant No. DMS-0901373.}~ and David Gepner}

\usepackage[pdfstartview=FitH,
            pdfauthor={\myauthor},
            pdftitle={\pdftitle},
            colorlinks,
            linkcolor=reference,
            citecolor=citation,
            urlcolor=e-mail,
            backref]{hyperref}

\usepackage{mathrsfs}
\usepackage[mathscr]{euscript}
\usepackage{amsmath}
\usepackage{amscd}
\usepackage{amsbsy}
\usepackage{amssymb}
\usepackage{microtype}
\usepackage{enumerate}
\usepackage[all,cmtip]{xy}
\usepackage{tikz}
\usetikzlibrary{matrix,arrows}
\usepackage{dsfont}
\usepackage[abbrev,lite,alphabetic]{amsrefs}
\usepackage{amsthm}
\usepackage{myarrows}

\usepackage{color}
\definecolor{todo}{rgb}{1,0,0}
\definecolor{conditional}{rgb}{0,1,0}
\definecolor{e-mail}{rgb}{0,.40,.80}
\definecolor{reference}{rgb}{.20,.60,.22}
\definecolor{mrnumber}{rgb}{.80,.40,0}
\definecolor{citation}{rgb}{0,.40,.80}

\newcommand{\icat}[0]{$\infty$-category~}
\newcommand{\icats}[0]{$\infty$-categories~}

\newcommand{\itopoi}[0]{$\infty$-topoi~}
\newcommand{\li}[0]{\Omega^{\infty}}
\DeclareMathOperator{\Ho}{Ho}
\DeclareMathOperator{\K}{K}
\newcommand{\proj}{\mathrm{proj}}
\newcommand{\perf}{\mathrm{perf}}

\DeclareMathOperator{\Dual}{D}
\DeclareMathOperator{\Rex}{Rex}

\DeclareMathOperator{\id}{id}

\DeclareMathOperator{\op}{op}
\DeclareMathOperator{\gp}{gp}
\DeclareMathOperator{\im}{im}

\DeclareMathOperator{\eq}{eq}
\DeclareMathOperator{\Eq}{Eq}
\DeclareMathOperator{\sm}{small}

\DeclareMathOperator*{\colim}{colim}
\newcommand{\desc}{\mathrm{desc}}

\newcommand{\rwe}{\tilde{\rightarrow}}
\newcommand{\riso}{\overset\simeq\rightarrow}
\newcommand{\we}{\simeq}
\newcommand{\iso}{\cong}
\newcommand{\ev}{\mathrm{ev}}
\newcommand{\coev}{\mathrm{coev}}
\newcommand{\st}{\mathrm{st}}
\newcommand{\tors}{\mathrm{tors}}

\newcommand{\Gm}{\mathds{G}_{m}}

\newcommand{\Ek}{\EE_k}
\newcommand{\Einfinity}{\EE_{\infty}}

\newcommand{\LCI}{\widehat{\Cat}_{\infty}}
\newcommand{\Prl}{\mathrm{Pr}^{\mathrm{L}}}

\newcommand{\Prlst}{\mathrm{Pr}^{\mathrm{L}}_{\mathrm{st}}}

\newcommand{\Gpd}{\widehat{\mathrm{Gpd}}_{\infty}}
\DeclareMathOperator{\Ind}{Ind}

\newcommand{\Sp}{\mathrm{Sp}}
\newcommand{\Spaces}{\mathrm{Spaces}}
\newcommand{\Mod}{\mathrm{Mod}}

\newcommand{\Alg}{\mathrm{Alg}}
\newcommand{\CAlg}{\mathrm{CAlg}}

\DeclareMathOperator{\cn}{cn}
\newcommand{\Az}{\mathrm{Az}}

\newcommand{\Aff}{\mathrm{Aff}}

\newcommand{\Shv}{\mathrm{Shv}}

\newcommand{\Test}{\mathrm{Test}}

\newcommand{\Cat}{\mathrm{Cat}}

\newcommand{\Idem}{\mathrm{Idem}}

\DeclareMathOperator{\Spec}{Spec}

\DeclareMathOperator{\GL}{GL}

\DeclareMathOperator{\Sym}{Sym}
\DeclareMathOperator{\Pic}{Pic}

\DeclareMathOperator{\Hoh}{H}
\DeclareMathOperator{\Eoh}{E}

\DeclareMathOperator{\Tor}{Tor}

\newcommand{\Map}{\mathrm{Map}} 
\newcommand{\Fun}{\mathrm{Fun}}
\newcommand{\map}{\mathrm{map}} 
\newcommand{\der}{\mathrm{der}}
\DeclareMathOperator{\End}{End}

\DeclareMathOperator{\aut}{aut}

\newcommand{\ShPic}{\mathbf{Pic}}
\newcommand{\ShGL}{\mathbf{GL}}
\newcommand{\Shpic}{\mathbf{pic}}

\newcommand{\ShM}{\mathbf{M}}

\newcommand{\ShP}{\mathbf{Pr}}

\newcommand{\ShAlg}{\mathbf{Alg}}
\newcommand{\ShAz}{\mathbf{Az}}
\newcommand{\ShBr}{\mathbf{Br}}

\newcommand{\ShMr}{\mathbf{Mr}}
\newcommand{\Shbr}{\mathbf{br}}

\newcommand{\ShProj}{\mathbf{Proj}}

\newcommand{\ShMor}{\mathbf{Mor}}
\newcommand{\Shmap}{\mathbf{map}}

\newcommand{\ShB}{\mathbf{B}} 

\newcommand{\StAlg}{\mathscr{A}\mathrm{lg}}

\newcommand{\StCat}{\mathscr{C}\mathrm{at}}
\newcommand{\StAz}{\mathscr{A}\mathrm{z}}
\newcommand{\StMod}{\mathscr{M}\mathrm{od}}
\newcommand{\StProj}{\mathscr{P}\mathrm{roj}}
\newcommand{\StM}{\mathscr{M}}

\newcommand{\et}{\mathrm{\acute{e}t}}

\newcommand{\Zar}{\mathrm{Zar}}

\DeclareMathOperator{\Br}{Br}
\DeclareMathOperator{\br}{br}

\newcommand{\alg}{\mathrm{alg}}

\newcommand{\B}{\mathrm{B}\!}
\newcommand{\Lrm}{\mathrm{L}}

\newcommand{\Mrm}{\mathrm{M}}
\newcommand{\Hrm}{\mathrm{H}\,}
\newcommand{\Nrm}{\mathrm{N}}
\newcommand{\Brm}{\mathrm{B}}
\newcommand{\Prm}{\mathrm{P}}

\newcommand{\Kscr}{\mathscr{K}}

\newcommand{\Oscr}{\mathscr{O}}

\newcommand{\Sscr}{\mathscr{S}}
\newcommand{\Mscr}{\mathscr{M}}

\newcommand{\Ascr}{\mathscr{A}}
\newcommand{\Bscr}{\mathscr{B}}
\newcommand{\Cscr}{\mathscr{C}} 
\newcommand{\Dscr}{\mathscr{D}}

\newcommand{\Xscr}{\mathscr{X}}

\renewcommand{\AA}{\mathds{A}}

\newcommand{\CC}{\mathds{C}}
\newcommand{\RR}{\mathds{R}}
\newcommand{\QQ}{\mathds{Q}}
\newcommand{\ZZ}{\mathds{Z}}
\newcommand{\EE}{\mathds{E}}
\newcommand{\FF}{\mathds{F}}

\newcommand{\GG}{\mathds{G}}
\renewcommand{\SS}{\mathds{S}}
\newcommand{\WW}{\mathds{W}}

\theoremstyle{plain}
\newtheorem{theorem}{Theorem}[section]
\newtheorem{lemma}[theorem]{Lemma}
\newtheorem{proposition}[theorem]{Proposition}

\newtheorem{corollary}[theorem]{Corollary}

\theoremstyle{definition}
\newtheorem{definition}[theorem]{Definition}

\newtheorem{example}[theorem]{Example}

\newtheorem{remark}[theorem]{Remark}

\let\oldmarginpar\marginpar
\renewcommand\marginpar[1]{\-\oldmarginpar[\raggedleft\footnotesize #1]%
{\raggedright\footnotesize #1}}

\usepackage{txfonts}

\begin{document}
\maketitle

\begin{abstract}
\noindent
In this paper, we study Azumaya algebras and Brauer groups in derived
algebraic geometry. We establish various fundamental facts about
Brauer groups in this setting, and we provide a computational tool, which we use to
compute the Brauer group in several examples. In particular, we show that the Brauer group
of the sphere spectrum vanishes, and we use this to prove two uniqueness theorems for the
stable homotopy category. Our key technical results include the local geometricity, in the sense
of Artin $n$-stacks, of the moduli space of perfect modules over a smooth and proper algebra, 
the \'etale local triviality of
Azumaya algebras over connective derived schemes, and a local to global principle for the
algebraicity of stacks of stable categories.


\paragraph{Mathematics Subject Classification 2010}
Primary: \href{http://www.ams.org/mathscinet/msc/msc2010.html?t=14Fxx&btn=Current}{14F22}, \href{http://www.ams.org/mathscinet/msc/msc2010.html?t=18Gxx&btn=Current}{18G55}.
Secondary: \href{http://www.ams.org/mathscinet/msc/msc2010.html?t=14Dxx&btn=Current}{14D20},
\href{http://www.ams.org/mathscinet/msc/msc2010.html?t=18Exx&btn=Current}{18E30}.
\end{abstract}

\tableofcontents

\section{Introduction}

\subsection{Setting}

Derived algebraic geometry is a generalization of classical Grothendieck-style
algebraic geometry aimed at bringing techniques from geometry to bear on problems in
homotopy theory, and used to unify and explain many disparate results about categories of sheaves
on schemes. It has been used by Arinkin-Gaitsgory~\cite{arinkin-gaitsgory} to formulate a precise version of the geometric
Langlands conjecture, by Ben-Zvi-Francis-Nadler~\cite{bzfn} to study integral transforms and
Hochschild homology of coherent sheaves, by Lurie and others to study topological modular forms, by
To\"en-Vaqui\'e~\cite{toen-vaquie} to study moduli spaces of complexes of vector bundles, and
by To\"en~\cite{toen-derived} to study derived Azumaya algebras. Moreover, the philosophy of
derived algebraic geometry is closely related to non-commutative geometry
and to the idea of hidden smoothness of Kontsevich.

The basic objects in derived algebraic geometry are ``derived'' versions of commutative rings.
There are various things this might mean.
For instance, it could mean simply a graded commutative ring, or a
commutative differential-graded ring, such as the de Rham complex $\Omega^*(M)$ of a
manifold $M$.
Or, it could mean a commutative ring spectrum, which is to say a spectrum
equipped with a coherently homotopy commutative and associative multiplication. The
basic example of such a commutative ring spectrum is the sphere spectrum $\SS$, which is the initial commutative ring spectrum, and hence plays the role of the integers $\ZZ$ in derived algebraic geometry. Commutative ring spectra are in a precise sense the universal class of derived commutative rings. We work throughout this paper with connective commutative ring spectra, their module categories, and their associated schemes.

While a substantial portion of the theory we develop in this paper has been studied
previously for simplicial commutative rings, it is important for applications to homotopy theory
and differential geometry to have results applicable to the much broader class of commutative (or $\EE_\infty$) ring spectra,
as the vast majority of the rings which arise in these contexts are only of this more general form. Simplicial commutative rings are special cases of commutative differential graded rings,
and an $\EE_\infty$-ring spectrum admits an $\EE_\infty$-dg model if and only if it is a commutative algebra over the Eilenberg-MacLane spectrum $\Hrm\ZZ$. To give an idea of how specialized a class this is,
note that an arbitrary spectrum $M$ is an $\Hrm\ZZ$-module precisely when all of its $k$-invariants are trivial, meaning that it decomposes as a product of spectra $\Sigma^n\Hrm\pi_n M$, or that it has no nontrivial extensions in its ``composition series'' (Postnikov tower). Rather, the basic $\EE_\infty$-ring is the sphere spectrum $\SS$, which is the group completion of the symmetric monoidal category of finite sets and {\em automorphisms} (as opposed to only {\em identities}, which yields $\ZZ$) and captures substantial information from differential and $\FF_1$-geometry and contains the homological complexity of the symmetric groups. Similarly, the algebraic $K$-theory spectra, as well as other important spectra such those arising from bordism theories of manifolds, in the study of the mapping class group and the Mumford conjecture, or in topological Hochschild or cyclic homology, tend not to exist in the differential graded world.

Nevertheless, an $\EE_\infty$-ring spectrum $R$ should be regarded as a nilpotent thickening of its underlying commutative ring $\pi_0 R$, in much the same way as the Grothendieck school successfully incorporated nilpotent elements of ordinary rings into algebraic geometry via scheme theory.
Of course, this relies upon the ``local'' theory of homotopical commutative algebra, which, thanks to
the efforts of many mathematicians, is now well established.
In particular, there is a good notion of \'etale map of commutative ring spectra, and so the basic geometric objects in our paper will be glued together, in this topology, from commutative ring spectra.
We adopt Grothendieck's ``functor of points'' perspective; specifically, we fix a base $\EE_\infty$-ring $R$ and consider the category of connective commutative $R$-algebras, $\CAlg_R^{\cn}$. A sheaf is then a space-valued functor on $\CAlg_R^{\cn}$ which satisfies descent for the \'etale topology in the appropriate homotopical sense.
For instance, if $S$ is a commutative $R$-algebra, there is the representable sheaf $\Spec S$ whose space of $T$-points is the mapping space $\map(S,T)$ in the $\infty$-category of connective commutative $R$-algebras.

Just as in ordinary algebraic geometry, one is really only interested in a subclass of sheaves which are geometric in some sense.
An important feature of derived algebraic geometry is the presence of higher versions of
Artin stacks, an idea due to Simpson~\cite{simpson}; roughly, this is the smallest class of sheaves which contain the representables $\Spec S$ and is closed under formation of quotients by smooth groupoid actions. By restricting attention to these sheaves,
it is possible to prove many EGA-style statements.
The situation is entirely analogous to the classes of schemes or algebraic spaces in ordinary algebraic geometry, which can be similarly expressed as the closure of the affines under formation of Zariski or \'etale quotients, respectively.
The difference is that we allow our sheaves to take values in spaces, a model for the theory of higher groupoids, and that we require the larger class which is closed under smooth actions, so that it contains objects such as the deloopings $\B^{\, n} A$ of a smooth abelian group scheme $A$. These are familiar objects: the Artin $1$-stack $\B\, A$ is the moduli space of $A$-torsors, and the Artin $2$-stack $\B^{\, 2} A$ is the moduli space of gerbes with band $A$.

One of the main goals of this paper is to study Azumaya algebras over these derived
geometric objects. Historically, the notion of Azumaya algebra, due to
Auslander-Goldman~\cite{auslander-goldman},
arose from an attempt to generalize the Brauer group of a field. It was then globalized by
Grothendieck~\cite{grothendieck-brauer-1},
who defined an Azumaya algebra $\Ascr$ over a scheme $X$ as a sheaf of coherent $\Oscr_X$-algebras that is \'etale locally a matrix algebra.
In other words, there is a surjective \'etale map $p:U\rightarrow X$ such that $p^*\Ascr\iso\Mrm_n(\Oscr_U)$. The Brauer group of a scheme classifies Azumaya algebras up to Morita equivalence, that is, up to equivalence
of their module categories. The original examples of Azumaya algebras are central simple
algebras over a field $k$; by Wedderburn's theorem,
these are precisely the algebras $\Mrm_n(D)$, where $D$ is a division algebra of finite
dimension over its center $k$. The algebra of quaternion numbers over $\RR$ is thus an example of an Azumaya $\RR$-algebra, and represents the generator of $\Br(\RR)\cong\ZZ/2$.

In more geometric settings, the first example of an Azumaya algebra is
the endomorphism algebra of a vector bundle, though these have trivial Brauer class.
Locally, any Azumaya algebra is the endomorphism algebra of a vector bundle, but the vector
bundles do not generally glue to a vector bundle on the total space.
However, every Azumaya algebra is the endomorphism algebra of a twisted vector bundle, a perspective that has recently gained a great deal of importance. For instance, in the theory
of moduli spaces of vector bundles, there is always a twisted universal vector bundle, and the class of its endomorphism algebra in the Brauer group is precisely the obstruction to the existence of a universal (non-twisted) vector bundle on the moduli space.
Brauer groups and Azumaya algebras play an
important role in many areas of mathematics, but especially in arithmetic geometry,
algebraic geometry, and applications to mathematical physics. In arithmetic geometry, they are closely related to Tate's conjecture on $l$-adic cohomology of schemes over finite fields,
and they play a critical role in studying rational points of
varieties through, for example, the Brauer-Manin obstructions to the Hasse principle.
In algebraic geometry, Azumaya algebras arise naturally when studying moduli
spaces of vector bundles, and Brauer classes appear when considering certain constructions
motivated from physics in homological mirror symmetry. The Brauer group was also used by
Artin-Mumford~\cite{artin-mumford} to construct one of the first examples of a non-rational unirational complex variety.

As an abstract group, defined via the above equivalence relation, the Brauer group
is difficult to compute directly. Instead, one introduces the cohomological Brauer group
of a scheme, $\Br'(X)=\Hoh^2_{\et}(X,\Gm)_{\tors}$. There is an inclusion
$\Br(X)\subseteq\Br'(X)$. A first critical problem, posed by Grothendieck, is whether this
inclusion is an equality. Unfortunately, the answer is ``no'' in general, although de Jong
has written a proof~\cite{dejong} of a theorem of O. Gabber that equality holds if $X$ is quasi-projective, or more generally has an ample line bundle. However, by expanding the
notion of Azumaya algebra to derived Azumaya algebra, as done
in~\cite{lieblich-thesis}*{Chapter 3} and To\"en~\cite{toen-derived}, the answer to the
corresponding question is ``yes,'' at least for quasi-compact and quasi-separated schemes.
This was shown by To\"en, who also shows that the result holds for quasi-compact
and quasi-separated derived schemes built from simplicial commutative rings. One of the purposes of the present paper is to generalize this theorem to quasi-compact and quasi-separated derived schemes based on connective commutative ring
spectra, which is necessary for our applications to homotopy theory.
To any class $\alpha\in\Br'(X)$ there is an associated category
$\Mod_X^{\alpha}$ of complexes of quasi-coherent $\alpha$-twisted sheaves. When this
derived category is equivalent to $\Mod_{\Ascr}$ for an ordinary Azumaya algebra $\Ascr$,
then $\alpha\in\Br(X)$. However, even when this fails, as long
as $X$ is quasi-compact and quasi-separated, there is a derived Azumaya algebra $\Ascr$ such that $\Mod_X^{\alpha}\we\Mod_\Ascr$. Hence derived Azumaya algebras are locally endomorphisms algebras of complexes of vector bundles, and not just vector bundles, and therefore the appropriate notion of Morita equivalence is based on tilting complexes instead of bimodules.

One of the main features of this category $\Mod_X^{\alpha}\simeq\Mod_\Ascr$ of quasi-coherent $\alpha$-twisted sheaves is that it
allows us to define the $\alpha$-twisted $K$-theory spectrum $K^\alpha(X)$ of $X$ as the $K$-theory of the subcategory of perfect objects (see Definition \ref{def:perfectobject}).
The reason this is sensible is that, given an Azumaya $\Oscr_X$-algebra $\Ascr$, there's an Azumaya $\Oscr_X$-algebra $\Bscr$ such that $\Mod_\Ascr\otimes\Mod_\Bscr\simeq\Mod_X$; moreover,
$\Bscr$ can be taken to be the opposite $\Oscr_X$-algebra $\Ascr^{\op}$ and $\Mod_{\Ascr^{\op}}\simeq\Mod_X^{-\alpha}$.
Note that because $\Br'(X)=\Hoh_{\et}^2(X;\GG_m)$, this is entirely analogous to what happens topologically,
where the twists are typically given by elements of the cohomology group
$\Hoh^2(X;\CC^\times)\cong\Hoh^3(X;\ZZ)$
in the complex case and elements of $\Hoh^2(X;\RR^\times)$ in the real case \cite{agg}.
While we do not study the twisted $K$-theory of derived schemes in this paper,
the basic structural features (such as additivity and localization) follow from the untwisted case as in \cite{toen-derived},
using the fact that our categories of $\alpha$-twisted sheaves $\Mod_X^{\alpha}\simeq\Mod_\Ascr$ admit global generators with endomorphism algebra $\Ascr$.

\subsection{Summary}

We now give a detailed summary of the paper.
By definition, an $R$-algebra $A$ is Azumaya if it is a
compact generator of the $\infty$-category of $R$-modules and if the multiplication action 
\begin{equation*}
    A\otimes_R A^{\op}\longrightarrow\End_R(A).
\end{equation*}
of $A\otimes_R A^{\op}$ on $A$ is an equivalence.
This definition is due to Auslander-Goldman~\cite{auslander-goldman} in the case of discrete
commutative rings, and it has been studied in the settings of schemes by
Grothendieck~\cite{grothendieck-brauer-1}, $\EE_\infty$-ring spectra by
Baker-Richter-Szymik~\cite{baker-richter-szymik},
and derived algebraic geometry over simplicial commutative rings by
To\"en~\cite{toen-derived}. In a slightly different direction, it has also been studied in
the setting of higher categories by Borceux-Vitale~\cite{borceux-vitale} and Johnson~\cite{johnson-azumaya}.
All of these variations ultimately rely on the idea of an Azumaya algebra as an algebra whose module category is invertible with respect to a certain ``Morita'' symmetric monoidal structure.

Although we restrict to Azumaya algebras over commutative ring spectra, we note that the notion of Azumaya algebra makes sense over any $\EE_3$-ring spectrum.
The reason for this is that if $R$ is an $\EE_3$-ring, then $\Mod_R$ is naturally
a $\EE_2$-monoidal $\infty$-category, and so its $\infty$-category of modules is naturally
$\EE_1$-monoidal.
The theory of Azumaya algebras is closely related to the notions of smoothness and
properness in non-commutative geometry, which have been studied extensively starting from
Kontsevich~\cite{kapranov}. These and related ideas have been used to great success to prove
theorems in algebraic geometry. For instance, van den Bergh~\cite{vandenbergh} uses non-commutative algebras
to give a proof of the Bondal-Orlov conjecture, showing that birational smooth projective $3$-folds are derived equivalent.

One of main points of the paper is to establish the following theorem, which says that all
Azumaya algebras over the sphere spectrum are Morita equivalent.

\begin{theorem}
    The Brauer group of the sphere spectrum is zero.
\end{theorem}

This theorem follows from several other important results, which we now outline.
In Theorem~\ref{thm:azumayaiffinvertible}, we show that a compactly
generated $R$-linear category $\Cscr$ (a stable presentable $\infty$-category enriched in $R$-modules) is dualizable if and only if it is of the form
$\Mod_A$ for a smooth and proper $R$-algebra $A$, and $\Cscr$ is invertible as a $\Mod_R$-module if and only if
$\Cscr\we\Mod_A$ for an Azumaya $R$-algebra $A$.
The analogous results were proved for simplicial commutative rings in~\cite{toen-derived}. A final algebraic ingredient is the fact that smooth and proper
$R$-algebras are compact. In particular, Azumaya algebras are compact algebras. This is a key point
later in our analysis of the geometricity of the sheaf of perfect modules for an Azumaya algebra.
To establish it requires showing that the $\infty$-category of spectra $\Sp$ is compact in
the $\infty$-category of all compactly generated $\SS$-linear categories, which does not
follow immediately from the fact that it is the unit object in this symmetric monoidal
$\infty$-category. The theory of smooth and proper algebras is also fundamental in the theory of
non-commutative motives, and has been studied in that setting by
Cisinski-Tabuada~\cite{cisinski-tabuada} and Blumberg-Gepner-Tabuada~\cite{bgt1}.

Suppose that $R$ is an $\EE_k$-ring spectrum for $3\leq k\leq\infty$. Then, the characterization
of Azumaya algebras above lets us define the Brauer space of $R$ as the Picard space
\begin{equation*}
    \Br_{\alg}(R)=\Pic(\Cat_{R,\omega})
\end{equation*}
of the $\EE_{k-2}$-monoidal $\infty$-category of compactly generated $R$-linear categories.
This space is a grouplike $\EE_{k-2}$-space, and so is, in particular, a $(k-2)$-fold
loop-space. The Brauer group is the abelian group
\begin{equation*}
    \pi_0\Br_{\alg}(R).
\end{equation*}
When $k=\infty$, it follows that there is a Brauer spectrum $\br_{\alg}(R)$.
One strength of this definition is that it generalizes well to other settings,
such as arbitrary compactly generated $\EE_{k-2}$-monoidal stable $\infty$-categories. We do
not develop this theory in our paper, instead working only with $\EE_\infty$-ring spectra,
but it is closely related to ideas about Brauer groups of $2$-categories.

Let us take a moment to place this idea in context. We can describe the space
$\Br_{\alg}(R)$ as follows. The $0$-simplices are $\infty$-categories $\Mod_A$ where $A$ is
an Azumaya $R$-algebra. A $1$-cell from $A$ to $B$ is an equivalence $\Mod_A\we\Mod_B$;
these may be identified with certain right $A^{\op}\otimes_R B$-modules. A $2$-cell is the
data of an equivalence between bimodules, and so forth. When $R$ is an ordinary ring, there
is no interesting data in degree higher than $2$. However, when $R$ is a derived ring, the
higher homotopy groups appear in the homotopy of $\Br_{\alg}(R)$. See \eqref{eq:134} below.
Thus, our Brauer space can be viewed as a generalization of the Brauer $3$-group of
Gordon-Power-Street~\cite{gordon-power-street} and Duskin~\cite{duskin}, and as a
generalization of the approach to Brauer groups in~\cite{vitale-brauer} and~\cite{borceux-vitale}.

The subject of derived algebraic geometry is increasingly important due to its
utility in proving theorems in homotopy theory and algebraic geometry.
As we will see in this paper, even to derive
purely homotopy-theoretic results about modules over the sphere spectrum, we will need to
employ derived algebraic geometry in an essentially non-trivial way. Such methods are
essential even in ordinary algebra. For instance, the classical
proof~\cite{grothendieck-brauer-3} that the Brauer
group of the integers vanishes employs geometric methods and cohomology.

In order to utilize cohomological methods to compute Brauer groups of derived schemes, it is
necessary to show that Azumaya algebras are locally Morita equivalent to the base.
This local triviality holds in the \'etale topology, but not in the Zariski topology (as is
shown by the quaternions over $\RR$). This is not easy to prove and uses the geometry of
smooth higher Artin sheaves. Higher Artin sheaves are built inductively out of affine
schemes by taking iterated quotients by smooth equivalence relations.
We study these sheaves in Section~\ref{sec:sheaves}, and we prove the
following theorem: if $f:X\rightarrow \Spec R$ is a smooth surjection, where $R$ is a
connective commutative ring spectrum and $X$ is filtered by higher Artin stacks, then $f$
has \'etale local sections. This extends the classical result about smooth morphisms of
schemes to derived algebraic geometry, and has been established in other contexts by
To\"en-Vezzosi~\cite{hag2}. To use this result on sections of smooth morphisms, we first
need to establish the following theorem, showing that a certain moduli sheaf is sufficiently
geometric; it is due to~\cite{toen-vaquie} in the simplicial commutative setting.

\begin{theorem}
    If $\Cscr$ is a stable $\infty$-category of finite type, then the moduli space
    $\ShM_\Cscr$ of compact objects of $\Cscr$ is locally geometric.
\end{theorem}

This is the case in particular for $\Cscr=\Mod_A$ when $A$ is an
Azumaya algebra, in which case the subsheaf $\ShMor_A\subseteq\ShM_A$ that classifies
Morita equivalences from $A$ to $R$ is smooth and surjective over $\Spec R$. This is used to
prove the following theorem.

\begin{theorem}
    If $A$ is an Azumaya algebra over a connective commutative ring spectrum $R$, then there
    is an \'etale cover $R\rightarrow S$ such that $A\otimes_R S$ is Morita
    equivalent to $S$. Thus, Azumaya algebras over any derived scheme are \'etale locally trivial.
\end{theorem}

For nonconnective commutative ring spectra the question of \'etale-local triviality is more
subtle. One possibility is to use Galois descent instead of \'etale descent. This is the
subject of current work by the second author and Tyler Lawson~\cite{gepner-lawson}.

In Section~\ref{sec:gluing}, we study families of linear categories over sheaves in order to establish the following key result regarding the existence of compact generators.

\begin{theorem}
    If $\Cscr$ is an $R$-linear category with descent such that $\Cscr\otimes_R S$ has a
    compact generator for some faithfully flat \'etale $R$-algebra $S$, then $\Cscr$ has a
    compact generator.
\end{theorem}

This local-global principle is proved by establishing analogous
statements for Zariski covers, for finite flat covers, and for Nisnevich covers. The method
for showing the Zariski local-global result follows work of Thomason~\cite{thomason-trobaugh},
B\"okstedt-Neeman~\cite{bokstedt-neeman}, Neeman~\cite{neeman-1992}, and Bondal-van den
Bergh~\cite{bondal-vandenbergh} on derived categories of schemes. The local-global principle
for finite flat covers is straightforward. The real work is in establishing the principle
for \'etale covers. Lurie proves in~\cite{dag11}*{Theorem~2.9} that the Morel-Voevodsky
theorem, which reduces Nisnevich descent to affine Nisnevich excision, holds for
connective $\EE_\infty$-ring spectra. Thus, we show a local-global principle for affine
Nisnevich squares. This idea parallels work of Lurie on a local-global principle for the
compact generation of linear categories (as opposed, in our work, for compact generation by a
single object). To\"en~\cite{toen-derived}
proves a similar local-global principle for fppf covers in the setting of simplicial commutative rings,
but his proofs both of the \'etale and the fppf local-global principles do not obviously generalize to $\EE_\infty$-ring spectra
because it is typically not the case that there are algebra structures on module-theoretic quotients of ring
spectra.

The local-global principle shows that if $\Cscr$ is a linear category with descent over a
quasi-compact and quasi-separated derived scheme such that $\Cscr$ is \'etale locally
equivalent to modules over an Azumaya algebra, then $\Cscr$ is globally equivalent to
modules over an Azumaya algebra. This solves the $\Br=\Br'$ problem of Grothendieck for
derived schemes.

\begin{theorem}
    If $X$ is a quasi-compact and quasi-separated scheme over the sphere spectrum,
    then $\Br(X)=\Br'(X)$.
\end{theorem}

The local-global principle has another interesting application: if $X$ is a quasi-compact
and quasi-separated derived scheme over the $p$-local sphere, then the $\infty$-category $L_{K(n)}\Mod_X$
of $K(n)$-local objects is compactly generated.

In Section~\ref{sec:brauer}, we define a Brauer sheaf $\ShBr$. If $X$ is an \'etale sheaf,
the Brauer space $\ShBr(X)$ of $X$ is the space of maps from $X$ to $\ShBr$ in the
$\infty$-topos $\Shv_R^{\et}$.
In the case of an affine scheme $\Spec R$, combining the \'etale-triviality of Azumaya
algebras and the \'etale local-global principle, we find that $\Br_{\alg}(R)\we\ShBr(\Spec R)$.
One advantage of using the Brauer sheaf $\ShBr$ is that it is a delooping of the Picard sheaf:
$\Omega\ShBr\we\ShPic$. This allows us to compute the homotopy sheaves of $\ShBr$:
\begin{equation*}
    \pi_k\ShBr\we\begin{cases}
        0   &   \text{$k=0$,}\\
        \ZZ &   \text{$k=1$,}\\
        \pi_0\Oscr^{\times} &   \text{$k=2$,}\\
        \pi_{k-2}\Oscr      &   \text{$k\geq 3$,}
    \end{cases}
\end{equation*}
where $\Oscr$ denotes the structure sheaf of $\Shv_R^{\et}$.
We introduce a computation tool, a descent spectral sequence
\begin{equation*}
    \Eoh_2^{p,q}=\Hoh^p(X,\pi_q\ShBr)\Rightarrow\pi_{q-p}\ShBr(X),
\end{equation*}
which converges if $X$ is affine or has finite \'etale cohomological dimension. When
$X=\Spec R$, the spectral sequence collapses, and we find that
\begin{equation}\label{eq:134}
    \pi_k\ShBr(R)\iso\begin{cases}
        \Hoh^1_{\et}(\Spec\pi_0 R,\ZZ)\times\Hoh^2_{\et}(\Spec\pi_0 R,\Gm) & k=0,\\
        \Hoh^0_{\et}(\Spec\pi_0R,\ZZ)\times\Hoh^1_{\et}(\Spec\pi_0 R,\Gm) & k=1,\\
        \pi_0 R^{\times}    & k=2,\\
        \pi_{k-2}R & k\geq 3.
    \end{cases}
\end{equation}

It follows that the Brauer group vanishes in many interesting cases; for example:
\begin{equation*}
    \pi_0\ShBr(ko)=0,\hspace{1cm}\pi_0\ShBr(ku)=0,\hspace{1cm}\pi_0\ShBr(MU)=0,\hspace{1cm}\pi_0\ShBr(tmf)=0.
\end{equation*}
For examples of non-zero Brauer groups, we find that
\begin{equation*}
    \pi_0\ShBr\left(\SS\left[\frac{1}{p}\right]\right)\we\ZZ/2,
\end{equation*}
and for the $p$-local sphere spectrum, the Brauer group fits into an exact sequence
\begin{equation*}
    0\rightarrow\pi_0\ShBr(\SS_{(p)})\rightarrow\ZZ/2\oplus\bigoplus_p\QQ/\ZZ\rightarrow\QQ/\ZZ\rightarrow 0.
\end{equation*}
Note that the $p$-inverted sphere and the $p$-local sphere give examples of
non-Eilenberg-MacLane $\EE_\infty$-ring spectra with non-zero Brauer groups.
By~\eqref{eq:134}, if $R$ is a connective $\EE_\infty$-ring spectrum, we can compute the
homotopy groups of $\ShBr(R)$ whenever we can compute the relevant \'etale cohomology groups
of $\Spec\pi_0R$. For example, $\pi_0\ShBr(R)=0$ if $R$ is any connective $\EE_\infty$-ring spectrum such that $\pi_0
R\iso\ZZ$ or $\WW_k$, the ring of Witt vectors over $\FF_{p^k}$.

We state as theorems two consequences of the vanishing of the Brauer group of the sphere
spectrum.

\begin{theorem}
    Let $\Cscr$ be a compactly generated stable presentable $\infty$-category, and suppose that
    there exists a stable presentable $\infty$-category $\Dscr$ such that
    $\Cscr\otimes\Dscr\we\Mod_{\SS}$, the $\infty$-category of spectra. Then,
    $\Cscr\we\Mod_{\SS}$.
\end{theorem}

\begin{theorem}
    Let $\Cscr$ be a stable presentable $\infty$-category such that there exists a
    faithfully flat \'etale $\SS$-algebra $T$ such that $\Cscr\otimes\Mod_T\we\Mod_T$. Then,
    $\Cscr\we\Mod_{\SS}$.
\end{theorem}

The first theorem says that if $\Cscr$ is compactly generated and invertible as a
$\Mod_\SS$-module, then $\Cscr$ is already equivalent to $\Mod_{\SS}$. The second theorem
says that if $\Cscr$ is \'etale locally equivalent to spectra, then $\Cscr$ is already
equivalent to the $\infty$-category of spectra. These give strong uniqueness, or rigidity,
results for $\SS$-modules. Such statements have a long history, and are related to the
conjecture of Margolis, which gives conditions for a triangulated category to be equivalent
to the stable homotopy category. The conjecture was proven for triangulated categories with
models by Schwede and Shipley~\cite{ss-uniqueness}. Our results extend theirs and
also those of~\cite{schwede}.

We briefly outline the contents of our paper.
We start in Section~\ref{sec:ringsmodules} by giving some background on rings, modules, and
the \'etale topology in the context of derived algebraic geometry. In
Section~\ref{sec:modulecategories}, we consider the module categories of $R$-algebras $A$
under various conditions, including compactness, properness, and smoothness. We prove there
the characterization that $A$ is Azumaya (resp. smooth and proper) if and only if $\Mod_A$
is invertible (resp. dualizable) in a certain symmetric monoidal $\infty$-category of $R$-linear categories.
We develop the theory of higher Artin sheaves in derived algebraic geometry in
Section~\ref{sec:sheaves}.
In Section~\ref{sec:moduli}, we harness the notion of geometric sheaves to study the moduli
space of $A$-modules for nice $R$-algebras $A$. Specializing to the case of Azumaya
algebras, we prove that the sheaf of Morita equivalences from $A$ to $R$ is smooth and
surjective over $\Spec R$, and hence has \'etale-local sections. It follows that Azumaya
$R$-algebras are \'etale locally trivial. We consider the problem of when a stack of linear
categories over a stack admits a perfect generator in Section~\ref{sec:gluing}. In the final
section, Section~\ref{sec:brauer}, we study the Brauer group, define the Brauer spectral
sequence, and give the computations, including the important theorem that the Brauer group
of the sphere spectrum vanishes.

\subsection{Acknowledgments}

Much    of    this    paper    is    based    on    ideas    of  Lieblich and  To\"en,
and the technical framework is supported on the work of Jacob Lurie. We
thank Tyler Lawson and Jacob Lurie for various comments and suggestions.
Finally, we want to even more specifically cite the work of To\"en as influencing ours.
For simplicial commutative rings, the key insights on moduli spaces of objects, the
\'etale local triviality of Azumaya algebras, and gluing of generators are due to
To\"en~\cite{toen-derived} and To\"en-Vaquie~\cite{toen-vaquie}.
Most of the work of Sections~\ref{sec:moduli} and~\ref{sec:gluing} is a recapitulation of
these ideas in the setting of ring spectra.

\section{Ring and module theory}\label{sec:ringsmodules}

In this section of the paper, we give some background on ring spectra and their module
categories, compactness, Grothendieck topologies on commutative ring spectra, and
Tor-amplitude.

\subsection{Rings and modules}\label{sub:ringsmodules}

The work~\cite{ha} of Lurie gives good notions of module categories for ring objects in
symmetric monoidal $\infty$-categories. We refer that book for details on the construction of the
objects introduced in the rest of this section. If $\Cscr$ is a symmetric monoidal
$\infty$-category, and if $A$ is an algebra object, by which we mean an $\EE_1$-algebra
in $\Cscr$, then there is an $\infty$-category $\Mod_A(\Cscr)$ of right $A$-modules in $\Cscr$;
similarly there is an $\infty$-category of left $A$-modules
$_{A}\Mod(\Cscr)\we\Mod_{A^{\op}}(\Cscr)$. Given two
algebras $A$ and $B$, there is an $\infty$-category $_A\Mod_B(\Cscr)$ of $(A,B)$-bimodules
in $\Cscr$, which is equivalent to $\Mod_{A^{\op}\otimes B}(\Cscr)$.
The $\EE_k$-algebras in $\Cscr$ form an $\infty$-category $\Alg_{\EE_k}(\Cscr)$. When $k=1$,
we write $\Alg(\Cscr)$ for this $\infty$-category, and when $k=\infty$, we write
$\CAlg(\Cscr)$.  When $\Cscr=\Sp$, the $\infty$-category of spectra with the smash product tensor
structure, we will  write  more  simply  $\Alg_{\EE_k}$,  $\Alg$,  $\CAlg$,  $\Mod_A$,
$_A\Mod$,  $_A\Mod_B$  and  so  forth   for   the   \icats   of   $\Ek$-ring   spectra,
associative ring spectra, commutative ring spectra, etc.
When $A$ is a discrete associative ring, then the $\infty$-category of right modules $\Mod_{\Hrm A}$ over the
Eilenberg-MacLane spectrum of $\Hrm A$ is equivalent to the $\infty$-category of chain
complexes on $A$.

\subsection{Compact objects and generators}\label{sub:compact}

We introduce the notion of compactness, which will play a crucial role in everything that
follows.

\begin{definition}
    Let  $\Cscr$  denote  an  \icat  which  is  closed  under  $\kappa$-filtered
    colimits~\cite{htt}*{Section~5.3.1}.
    A functor $f:\Cscr\rightarrow\Dscr$ is said to be {\em $\kappa$-continuous} if $f$ preserves
    $\kappa$-filtered colimits.  In the special case that  $f$  preserves  $\omega$-filtered
    colimits, we simply say that $f$ is continuous.
\end{definition}

\begin{definition}
    Let  $\Cscr$  denote  an  \icat  which  is  closed  under  $\kappa$-filtered   colimits.
    An object $x$ of $\Cscr$ is said to be $\kappa$-compact  if  the  mapping  space  functor
    \begin{equation*}
        \map_{\Cscr}(x,-):\Cscr\longrightarrow\Sscr
    \end{equation*}
    is $\kappa$-continuous, where $\Sscr$ is the \icat of spaces. We say that $x$ is compact if it is
    $\omega$-compact.
\end{definition}

\begin{definition}
    Let $\Cscr$ be an \icat which is closed under geometric realizations (in other words, colimits of simplicial diagrams).
    An object $x$ of $\Cscr$ is said to be {\em projective} if  the  mapping  space  functor
    \begin{equation*}
        \map_{\Cscr}(x,-):\Cscr\longrightarrow\Sscr
    \end{equation*}
    preserves geometric realizations.
\end{definition}

A compact projective object $x$ of an \icat $\Cscr$ corepresents a functor which preserves both filtered colimits and geometric realizations.
Both filtered colimits and the simplicial indexing category $\Delta^{\op}$ are examples of {\em sifted colimits}
(that is, colimits indexed by simplicial sets $K$ such that the diagram $K\to K\times K$ if cofinal),
and $x$ is compact projective if and only if $\map_\Cscr(x,-)$ preserves sifted colimits (see~\cite{htt}*{Corollary~5.5.8.17}).

An $\infty$-category with all small colimits is said to be $\kappa$-compactly generated if
the natural map $\Ind_\kappa(\Cscr^{\kappa})\rightarrow\Cscr$ is an equivalence, where
$\Ind_\kappa$ denotes the $\kappa$-filtered cocompletion~\cite{htt}*{Section~5.3.5}.
By definition, a presentable $\infty$-category is an $\infty$-category that is
$\kappa$-compactly generated for some infinite regular cardinal $\kappa$. When $\Cscr$ is
$\omega$-compactly generated, we say simply that it is compactly generated. If $\Cscr$ is
compactly generated, and if $\Dscr$ is a full subcategory of $\Cscr^{\omega}$ such that the
closure of $\Dscr$ in $\Cscr$ under finite colimits and retracts is equivalent to $\Cscr^{\omega}$, then
we say that $\Cscr$ is compactly generated by $\Dscr$. 

\begin{lemma}\label{lem:generating}
    A stable presentable $\infty$-category $\Cscr$ is compactly generated by a set
    $X$ of compact objects if for any object $y\in\Cscr$,
    $\Map_{\Cscr}(x,y)\we 0$ for all $x\in X$ if and only if $y\we\ast$.
    \begin{proof}
        See the proof of~\cite{schwede-shipley}*{Lemma~2.2.1}.
    \end{proof}
\end{lemma}

Returning to the algebraic situation of the previous section,
if $A$ is an $\EE_1$-ring spectrum, then $\Mod_A$ is a stable presentable $\infty$-category.
Presentability follows from~\cite{ha}*{Corollary~4.2.3.7.(1)} and stability is straightforward. In particular, $\Mod_A$ admits all small colimits.
Moreover, $\Mod_A$ is compactly generated by the single object $A$.
We will often refer to compact objects of $\Mod_A$ as perfect modules, in keeping with the usual terminology of
algebraic geometry.

\begin{definition}
    A connective  $A$-module  $P$  is   projective   if   it   is   projective   as   an   object   of
    $\Mod_A^{\cn}$, the $\infty$-category of connective $A$-modules.
\end{definition}

The following argument shows why there are no non-zero projective objects of $\Mod_A$ in general.
Suppose that $M$ is a projective object of $\Mod_R$. For any $R$-module $N$, we can write
the suspension of $N$ as the geometric realization
\begin{equation*}
    \Sigma N\we\left|0\llarrow N\lllarrow N\oplus N\llllarrow\ldots\right|
\end{equation*}
Then, by stability,
\begin{equation*}
    \map_R(\Sigma^{-1}M,N)\we\map_R(M,\Sigma N)\we|\map(M,N^{\oplus n})|\we\Brm\map_R(M,N).
\end{equation*}
In particular, $\pi_0\map_R(\Sigma^{-1}M,N)=0$ for all $R$-modules $N$, so that $\id_M\we
0$. Thus, $M\we 0$.

We record here a few facts about projective and compact modules.

\begin{proposition}\label{prop:compactness}
    Let $A$ be a connective $\EE_1$-ring spectrum.
    \begin{enumerate}
        \item[\rm{(1)}] A connective right $A$-module $P$ is projective if  and  only  if  it
            is a retract of a free right $A$-module.
        \item[\rm{(2)}] A connective right $A$-module $P$ is projective if and  only  if  for
            every  surjective  map  $M\rightarrow   N$   of   right   $A$-modules   the   map
            \begin{equation*}
                \map(P,M)\rightarrow\map(P,N)
            \end{equation*}
            is surjective.
        \item[\rm{(3)}] A right $A$-module $P$ is compact if and only if  it  is  dualizable:
            there  exists  a  left  $A$-module  $P^{\vee}$  such   that   the   composition
            \begin{equation*}
                \Mod_A\xrightarrow{\otimes_A P^{\vee}}\Sp\xrightarrow{\Omega^{\infty}}\Sscr
            \end{equation*}
            is equivalent to the functor corepresented by $P$. In this case, $P^{\vee}$ is a
            perfect left $A$-module.
        \item[\rm{(4)}] If $P$  is  a non-zero compact right  $A$-module, then $P$ has a
            bottom non-zero homotopy group; that is,  there  exists  some
            integer $N$ such that
            \begin{equation*}
                \pi_n P=0
            \end{equation*}
            for $n\leq N$ and $\pi_{N+1}P\neq 0$. Moreover, $\pi_{N+1}P$ is finitely presented as a $\pi_0
            A$-module.
    \end{enumerate}
    \begin{proof}
        Part (1) is~\cite{ha}*{Proposition~8.2.2.7}. The proof of part (2) is the same as in
        the discrete case. Part (3)
        is~\cite{ha}*{Proposition~8.2.5.4}. Part (4) is~\cite{ha}*{Corollary~8.2.5.5}.
    \end{proof}
\end{proposition}

The following lemma will be used later in the paper.

\begin{lemma}\label{lem:functoriality}
    Let $R$ be a commutative ring spectrum, $A$ an $R$-algebra, and $P$ and $Q$ compact right
    $A$-modules.   Then,  for   any   commutative   $R$-algebra   $S$,   the   natural   map
    \begin{equation}\label{eq:10000}
        \Map_A(P,Q)\otimes_R  S\rightarrow\Map_{A\otimes_R  S}(P\otimes_R  S,Q\otimes_R   S)
    \end{equation}
    is an equivalence of $S$-modules.
    \begin{proof}
        The statement is clear when $P$ is a suspension of the free $A$-module $A$. We prove
        the   lemma   by    induction    on    the    cells    of    $P$.     So,    suppose
        that we have a cofiber sequence of compact modules
        \begin{equation*}
            \Sigma^n A\rightarrow N\rightarrow P
        \end{equation*}
        such that
        \begin{equation}\label{eq:10001}
            \Map_A(N,Q)\otimes_R S\rightarrow\Map_{A\otimes_R S}(N\otimes_R S,Q\otimes_R  S)
        \end{equation}
        is an equivalence. We show that~\eqref{eq:10000} holds.
        We obtain a morphism of cofiber sequences
        \begin{equation*}
            \begin{CD}
                \Map_A(\Sigma^n   A,Q)\otimes_R   S   @>>>   \Map_A(N,Q)\otimes_R   S   @>>>
                \Map_A(P,Q)\otimes_R S\\
                @VVV                            @VVV                @VVV\\
                \Map_{A\otimes_R S}(\Sigma^n A\otimes_R S,Q\otimes_R S) @>>>
                \Map_{A\otimes_R  S}(N\otimes_R  S,Q\otimes_R   S)   @>>>   \Map_{A\otimes_R
                S}(P\otimes_R S,Q\otimes_R S).
            \end{CD}
        \end{equation*}
        Since the left  two  vertical  arrows  are  equivalences,  the  right  arrow  is  an
        equivalence. To finish the proof, we show that if~\eqref{eq:10001} holds for $N$ and
        if $P$ is a retract of $N$, then~\eqref{eq:10000} holds.  If $N\we P\oplus M$,  then
        there is a commutative square of equivalences
        \begin{equation*}
            \begin{CD}
                \Map_A(N,Q)\otimes_R S    @>>>    \Map_A(P,Q)\otimes_R S\oplus\Map_A(M,Q)\otimes_R S\\
                @VVV                            @VVV\\
                \Map_{A\otimes_R  S}(N\otimes_R  S,Q\otimes_R   S)   @>>>   \Map_{A\otimes_R
                S}(P\otimes_R S,Q\otimes_R S)\oplus\Map_{A\otimes_R S}(M\otimes_R S,Q\otimes_R S).
            \end{CD}
        \end{equation*}
        It  follows  (by  looking  for  instance  at  cofibers   of   the   vertical   maps)
        that~\eqref{eq:10000} holds.
    \end{proof}
\end{lemma}

The following form of the Morita theorem is used frequently to show that certain
$\infty$-categories are categories of modules over some $\EE_1$-ring spectrum.

\begin{theorem}[Morita theory]\label{thm:morita}
    Let $\Cscr$ be a stable presentable $\infty$-category, and let $P$ be an object
    of $\Cscr$. Then, $\Cscr$ is compactly
    generated by $P$ if and only if
    \begin{equation}\label{eq:543}
        \Cscr\xrightarrow{\Map_{\Cscr}(P,-)}\Mod_{\End_{\Cscr}(P)^{\op}}
    \end{equation}
    is an equivalence.
    \begin{proof}
        One direction is the theorem of Schwede and Shipley, in the form found in
        Lurie~\cite{ha}*{Theorem~8.1.2.1}. So, suppose that~\eqref{eq:543} is an
        equivalence. The functor $\Map_{\Cscr}(P,-)$ automatically preserves filtered
        colimits because it is an equivalence,
        so we see that $P$ is compact in $\Cscr$. Since $\Map_{\Cscr}(P,-)$ is conservative,
        it follows from Lemma~\ref{lem:generating} that $\Cscr$ is compactly generated by $P$.
    \end{proof}
\end{theorem}

\subsection{Topologies on affine connective derived schemes}\label{sub:topologies}

Fix an $\Einfinity$-ring $R$, and denote by $\CAlg_R^{\cn}$
the \icat of connective commutative $R$-algebras.
Set $\Aff^{\cn}_R=(\CAlg_R^{\cn})^{\op}$, the \icat of affine connective derived schemes over $\Spec R$.
We make extensive use of \itopoi arising from Grothendieck topologies on $\Aff^{\cn}_R$. For
details    on    the    construction    of    these    \itopoi    see~\cite{htt}*{Chapter~6}
and~\cite{dag7}*{Section~5}. All of these topologies arise from pre-topologies
consisting of special classes of flat morphisms, a notion we now define.

\begin{definition}\label{def:flat}
    A  morphism  $f:S\rightarrow  T$  of  commutative ring spectra   is   called   flat   if
    \begin{equation*}
        \pi_0(f):\pi_0 S\rightarrow \pi_0 T
    \end{equation*}
    is  a  flat   morphism   of   discrete   rings   and   if   $f$   induces   isomorphisms
    \begin{equation*}
        \pi_k S\otimes_{\pi_0 S}\pi_0 T\riso\pi_k T
    \end{equation*}
    for all integers $k$.
\end{definition}

It is useful to use flatness to give a definition of many other properties of morphisms of
$\EE_\infty$-ring spectra.

\begin{definition}
    If $\mathrm{P}$ is a property of flat morphisms of discrete commutative rings, such as
    being faithful or \'etale,
    then a morphism $f:R\rightarrow T$ of commutative rings is said to be $\mathrm{P}$ if
    $f$ is flat in the sense of Definition~\ref{def:flat} and if $\pi_0(f)$ is $\mathrm{P}$.
\end{definition}

The Zariski topology on $\Aff_R^{\cn}$ is the Grothendieck topology generated by Zariski
open covers. Here, a map $\Spec T\rightarrow\Spec S$ is a Zariski open cover if the
associated map on ring spectra $S\rightarrow T$ is flat and induces a Zariski-open cover
$\Spec\pi_0T\rightarrow\Spec\pi_0S$. The associated $\infty$-topos of Zariski sheaves is denoted by $\Shv_R^{\Zar}$.
Similarly, there is an \'etale topology on $\Aff_R^{\cn}$ and an associated \'etale
$\infty$-topos $\Shv_R^{\et}$. A map $\Spec T\rightarrow\Spec S$ is \'etale if $S\rightarrow
T$ is flat and \'etale. Both of these Grothendieck topologies are constructed, explicitly,
via the method of~\cite{dag7}*{Proposition~5.1}. See~\cite{dag7}*{Proposition~5.4} for how
to do this for the flat topology.

\subsection{Tor-amplitude}

Most  of  the  material  below  on  Tor-amplitude  and   perfect   modules   was   developed
in~\cite{sga6}*{Expos\'e I}. We refer also to the exposition
in~\cite{thomason-trobaugh}.  In the simplicial commutative setting, this is treated in To\"en-Vaqui\'e~\cite{toen-vaquie}.
Throughout this section, $R$ is a connective commutative  ring spectrum. We
refer to compact $R$-modules as perfect $R$-modules. This is to agree with the terminology
in the references. Over a scheme $X$, a complex of quasi-coherent $\Oscr_X$-modules is
called perfect if its restriction to any affine subscheme is perfect, or, equivalently,
compact. While the perfect and compact modules agree for affine schemes, on a general scheme
$X$ not every perfect module is compact.

\begin{definition}
    An   $R$-module   $P$   has   Tor-amplitude   contained   in   the   interval    $[a,b]$
    if  for  any  $\pi_0  R$-module  $M$  (any module, not any complex of modules),
    \begin{equation*}
        \Hoh_i(P\otimes_R M)=0
    \end{equation*}
    for $i\notin [a,b]$. If such integers $a,b$ exist, then $P$ is said to have finite Tor-amplitude.
\end{definition}

If $P$ is an $R$-module, then $P$ has Tor-amplitude contained in $[a,b]$ if and only if $P\otimes_R \pi_0
R$  is  a  complex  of  $\pi_0  R$-modules  with  Tor-amplitude  contained  in  $[a,b]$ in
the ordinary sense. Note, however, that our definition  differs  from  that
in~\cite{sga6}*{I  5.2}  simply  in  that  we  work with homology instead of cohomology.

The   next   proposition   is   used   in   the   proof   of    the proposition that follows,
but it seems interesting in its own right.

\begin{proposition}\label{prop:projectives}
    The functor
    \begin{equation*}
        \pi_0:\Ho(\Mod_R^{\proj})\rightarrow\Mod_{\pi_0 R}^{\proj}
    \end{equation*}
    is an equivalence, where the decoration $\proj$ denotes the full subcategory of projective modules.
    \begin{proof}
        This is a special case of~\cite{ha}*{Corollary~8.2.2.19}.
		The  analogous  map  on  free  modules  is  an  equivalence.
        Since projectives are summands of free modules, we deduce
        that the functor $\pi_0$ above is fully faithful.
        
        Let $P$ be a projective $\pi_0 R$-module. Then, there exists a free $\pi_0 R$-module
		$F$ and an idempotent homomorphism $e:F\rightarrow F$ such that  $P$  is  the  image
		of   $e$.    By   definition,   $P$   is    also    the    filtered    colimit    of
        \begin{equation*}
            F\xrightarrow{e} F\xrightarrow{e}F\rightarrow\cdots
        \end{equation*}
        in $\Mod_{\pi_0 R}$. Lift the diagram to a diagram of free $R$-modules $F'$, and let $P'$ be
        the filtered colimit. The $R$-module $P'$ is projective because we can construct a splitting of
        $F'\rightarrow P'$ by mapping $F'$ to each $F'$ in the diagram via the idempotent $e'$. Then, $\pi_0(P')$ is isomorphic to $P$. So, the functor is essentially
        surjective, and hence an equivalence of categories.
    \end{proof}
\end{proposition}

The following proposition provides the technical results needed on perfect
modules over connective $\EE_\infty$-ring spectra. In particular, parts (4)-(7) will be the
key to giving certain inductive proofs about the moduli of objects in module categories in
Section~\ref{sec:moduli}. We emphasize again that it was the insight of To\"en-Vaqui\'e~\cite{toen-vaquie} that suggests this approach to studying perfect objects in the
context of simplicial commutative rings.

\begin{proposition}\label{prop:tor}
    Let $P$ and $Q$ be $R$-modules.
    \begin{enumerate}
        \item[\rm{(1)}]  If  $R$  is   perfect,   then   $P$   has   finite   Tor-amplitude.
        \item[\rm{(2)}]   If $R'$ is a connective commutative $R$-algebra, and if $P$ is an $R$-module with Tor-amplitude contained
            in $[a,b]$, then $P\otimes_R R'$ is an $R'$-module with Tor-amplitude contained in $[a,b]$.
        \item[\rm{(3)}]   If $P$ has Tor-amplitude contained in $[a,b]$ and $Q$ has Tor-amplitude contained in $[c,d]$,
            then   $P\otimes_R   Q$   has   Tor-amplitude    contained    in    $[a+c,b+d]$.
        \item[\rm{(4)}]   If $P$ and $Q$ have Tor-amplitude contained in $[a,b]$, then for any morphism $P\rightarrow Q$,
            the cofiber has Tor-amplitude contained in $[a,b+1]$.
            Dually,    the    fiber    has    Tor-amplitude    contained    in    $[a-1,b]$.
        \item[\rm{(5)}]   If $P$ is a perfect $R$-module with Tor-amplitude contained in $[0,b]$, with $0\leq
            b$,  then   $P$   is   connective,   and   $\pi_0   P=\Hoh_0(P\otimes_R\pi_0R)$.
        \item[\rm{(6)}] If $P$ is  perfect  and  has  Tor-amplitude  contained  in  $[a,a]$,
            then $P$ is equivalent to $\Sigma^{a}M$ for a finitely generated projective $R$-module $M$.
        \item[\rm{(7)}]   If $P$ is perfect and has Tor-amplitude contained in $[a,b]$, then there
            exists a morphism
            \begin{equation*}
                \Sigma^a M\rightarrow P
            \end{equation*}
            such that $M$ is a finitely generated projective $R$-module and the  cofiber  is
            perfect and has Tor-amplitude contained in $[a+1,b]$.
    \end{enumerate}
    \begin{proof}
        Part (1) follows from~\cite{thomason-trobaugh}*{Proposition~2.2.12}
        and~\cite{thomason-trobaugh}*{Proposition~2.3.1.(d)}.    That   the    notions    of
        perfection in Thomason-Trobaugh and Lurie agree is explained
        by~\cite{thomason-trobaugh}*{Theorem~2.4.4}, which is applicable here as the modules
        which appear  in  $\Mod_R$  are  all  quasi-coherent,  and  so  have  quasi-coherent
        homology. Parts
        (2)  and  (3)  are~\cite{sga6}*{Proposition~5.6}.   If  $C$  is   the   cofiber   of
        $P\rightarrow Q$, and if $M$ is a $\pi_0 R$-module, then
        \begin{equation*}
            P\otimes_R M\rightarrow Q\otimes_R M\rightarrow C\otimes_R M
        \end{equation*}
        is a cofiber sequence in $\Mod_{\pi_0 R}$.  The case of a fiber is dual. Thus, part (4) follows immediately  from
        the long exact sequence in homology.

        Consider the $\Tor$ spectral sequence
        \begin{equation*}
            \Eoh^{2}_{p,q}=\Tor_p^{\pi_* R}(\pi_*P,\pi_0 R)_q\Rightarrow\pi_{p+q}(P\otimes_R \pi_0
            R)=\Hoh_{p+q}(P\otimes_R\pi_0R)
        \end{equation*}
        with differentials $d^r_{p,q}$ of degree $(-r,r-1)$ constructed in~\cite{ekmm}.
        If $P$ is a non-trivial perfect $R$-module with Tor-amplitude contained in $[0,b]$, then the abutment
        of the spectral sequence is $0$ when $p+q<0$. We know that $P$ has a bottom homotopy group,
        say $\pi_k$.  That is, $\pi_k P$ is non-zero, and  $\pi_j  P=0$  for  $j<k$.
        Calculating the graded tensor product, we see that $\Eoh^2_{0,k}$ is the coequalizer of
        \begin{equation*}
            \bigoplus_{i+j=k}\pi_iP\otimes_{\pi_0R}\pi_jR\rightrightarrows\pi_kP
        \end{equation*}
        in the category of graded $\pi_0R$-modules. So,
        $\Eoh^2_{0,k}=\pi_kP$ as a $\pi_0R$-module. But, by our hypothesis on $k$, no non-zero
        differential may hit $E^2_{0,k}$. All differentials out are zero for degree reasons. It follows
        that $\pi_{k}(P\otimes_R\pi_0 R)\neq 0$. Therefore, $k\geq 0$, and $P$ is connective. This
        proves the first statement of part (5), and the second statement follows easily from
        the same argument.

        To prove part (6), we may assume that $a=0$.
        By~\cite{thomason-trobaugh}*{Proposition~2.3.1.(d)}, we may assume that $P\otimes_R\pi_0R$ is
        a bounded complex of finitely generated projective $\pi_0R$-modules. Because the kernel of a
        surjective map of finitely generated projective modules is finitely generated
        projective, by induction, the good truncation $\tau_{\geq 0}P\otimes_R\pi_0R\rwe
        P\otimes_R\pi_0R$ is a bounded complex of finitely generated projective
        $\pi_0R$-modules that is concentrated in non-negative degrees.
        We show now that $\pi_0 P$ is a projective $\pi_0R$-module. By part (5),
        $\pi_0P\iso\Hoh_0(P\otimes_R\pi_0R)$. 
        Since the homology is zero above degree $0$, the good truncation $\tau_{\geq
        0}P\otimes_R\pi_0R$ is a
        resolution of the finitely presented $\pi_0R$-module
        $\Hoh_0(P\otimes_R\pi_0R)$ by finitely generated projective $\pi_0R$-modules.
        It suffices to show that $\Hoh_0(P\otimes_R\pi_0R)$ is flat
        by~\cite{matsumura}*{Theorem~7.12}. If $M$ is a $\pi_0R$-module, the Tor spectral
        sequence computing $\Hoh_*(P\otimes_R\pi_0R\otimes_{\pi_0R}M)$ is
        \begin{equation*}
            \Eoh^2_{p,q}=\Tor_p^{\pi_0R}(\Hoh_*(P\otimes_R\pi_0R),M)_q\Rightarrow\Hoh_{p+q}(P\otimes_R\pi_0R\otimes_{\pi_0R}M).
        \end{equation*}
        But, for $q>0$, $\Eoh_{p,q}^2=0$, so that for $p>0$
        \begin{equation*}
            \Tor_p^{\pi_0R}(\Hoh_0(P\otimes_R\pi_0R),M)\iso\Hoh_p(P\otimes_RM)=0,
        \end{equation*}
        by the Tor-amplitude of $P$. Thus, $\Hoh_0(P\otimes_R\pi_0R)$ is flat.

        Thus, by the previous theorem and the connectivity of $P$, there is  a  natural  map
        \begin{equation*}
            Q\rightarrow P,
        \end{equation*}
        where $Q$ is a finitely generated projective $R$-module and $\pi_0 Q\iso\pi_0 P$. It suffices to show that the cofiber
        $C$ of this map is equivalent to zero. The $R$-module $C$ is perfect and has the property that $C\otimes_{R}\pi_0 R$
        is zero. Let $\pi_k C$ be the first non-zero homotopy group of $C$. Then, $\Sigma^{-k}C$ is
        connective with  Tor-amplitude  contained  in  $[0,0]$.    By   part   (5),   $\pi_k
        C=\Hoh_0(\Sigma^{-k}C\otimes_R\pi_0R)=0$, a contradiction. Thus, $C\we 0$.

        To prove  part  (7),  we  assume  that  $a=0$. If $b=0$, the statement follows from
        part (6). Thus, assume that $b>0$, and consider  $P\otimes_R  \pi_0
        R$,  which  is  a  perfect  complex  over   $\pi_0   R$   with   bounded   homology.
        As above, we may assume that $P\otimes_R\pi_0R$ is in fact a bounded complex  of  finitely  generated  projective
        $\pi_0 R$-modules concentrated in non-negative degrees. Thus, there is a
        natural  morphism  of  complexes  $Z_0\rightarrow P\otimes_R\pi_0   R$
        which  induces  a  surjection  in  degree   $0$   homology.    Lift   $Z_0$   to   a
        finitely generated projective $R$-module $M$, by Proposition \ref{prop:projectives}. We can write $Z_0$ as a split
        summand of $\pi_0R^n$, and hence $M$ as a split summand of $R^n$. Since $P$ is connective by
        part (5), the composition
        \begin{equation*}
            \pi_0R^n\rightarrow Z_0\rightarrow P\otimes_R\pi_0R
        \end{equation*}
        lifts to a map $R^n\rightarrow P$. Composing with $M\rightarrow R^n$, we obtain a map
        $M\rightarrow P$ which is a surjection on $\Hoh_0$.  By the long exact  sequence  in
        homology, the cofiber has Tor-amplitude contained in $[1,b]$ (remembering that $b>0$).
        Moreover, the cofiber is  perfect  by  the  two-out-of-three  property  for  perfect
        modules~\cite{thomason-trobaugh}*{Proposition~2.2.13.(b)}.
    \end{proof}
\end{proposition}

\subsection{Vanishing loci}

We show that the complement of the support of a perfect complex on an affine derived scheme
is a quasi-compact open subscheme. Recall that a morphism of schemes $X\rightarrow Y$ is quasi-compact if for every open affine  $\Spec
R$ of $Y$, the pullback $X\times_Y\Spec R$ is  quasi-compact~\cite{ega1}*{Definition  I
6.1.1}. The following result is due to Thomason~\cite{thomason-triangulated} in the ordinary setting
of discrete rings, and to To\"en-Vaqui\'e~\cite{toen-vaquie} for simplicial commutative rings.

\begin{proposition}\label{prop:vanishing}
    Let $R$ be a connective commutative ring spectrum, and let $P$ be a perfect $R$-module.
    The subfunctor $V_P\subseteq\Spec R$ of points $R\rightarrow S$ such that $P\otimes_R S$
    is   quasi-isomorphic   to   zero   is   a   quasi-compact   Zariski   open   immersion.
    \begin{proof}
        If $R$ is discrete, the proposition is~\cite{thomason-triangulated}*{Lemma 3.3.c}.
        To prove the proposition when $R$ is a connective commutative ring spectrum, let $Q=P\otimes_R  \pi_0  R$,
        which is  a  perfect  complex  of  $\pi_0  R$-modules.
        Let $V_Q$ be the quasi-compact Zariski open subscheme of $\Spec\pi_0 R$ specified by the vanishing
        of  $Q$  by  the discrete case. Choose elements   $f_1,\ldots,f_n\in\pi_0   R$   such
        that $V_Q$ is the union of the $\Spec\pi_0 R[1/f_i]$.
        We claim that $V_P$ is  the  union  $V$  of  the  $\Spec  R[1/f_i]$  in  $\Spec  R$.
        But, because $P$ is a perfect $R$-module, given an $S$-point $\Spec S\rightarrow V\subseteq\Spec R$ of $V$,  then
        \begin{equation*}
            \left(P\otimes_R S\right)\otimes_S\pi_0 S\we 0
        \end{equation*}
        if and only if
        \begin{equation*}
            P\otimes_R S\we 0.
        \end{equation*}
        Indeed, $P\otimes_R S$ has a bottom homotopy group, say of degree $k$, and it follows from the proof of
        Proposition~\ref{prop:tor} part (5) that
        \begin{equation*}
            \pi_kP\otimes_RS\iso\Hoh_k\left( \left(P\otimes_R
            S\right)\otimes_S\pi_0S\right).\qedhere
        \end{equation*}
    \end{proof}
\end{proposition}

\section{Module categories and their module categories}\label{sec:modulecategories}

In this section, we examine the algebra of module categories of $\EE_\infty$-ring spectra, viewed
as $\EE_{\infty}$-monoids in the \icat of stable presentable $\infty$-categories.
This leads to an important module-theoretic characterization of Azumaya $R$-algebras for an
$\EE_\infty$-ring spectrum $R$: an $R$-algebra $A$ is Azumaya if and only if $\Mod_A$ is an invertible
$\Mod_R$-module.

\subsection{$R$-linear categories}

In~\cite{htt}*{Chapter 5}, Lurie constructs the \icat $\Prl$ of presentable \icats and colimit
preserving functors. We refer to Lurie's book for the precise definition and properties of
presentable $\infty$-categories. For us, the main points are that a presentable \icat is closed under small limits and colimits and is $\kappa$-compactly generated for some infinite regular
cardinal $\kappa$. Moreover, the \icat $\Prl$ is also closed under small limits and colimits, and there is a symmetric monoidal structure 
on $\Prl$ with unit object the \icat of pointed spaces~\cite{ha}*{Section 6.3}.

A critical fact about $\Prl$ is that if
$R$ is an $\EE_k$-ring spectrum ($1\leq k\leq\infty$), then the \icat of right $R$-modules $\Mod_R$ is an
$\EE_{k-1}$-monoidal stable presentable $\infty$-category with unit $R$ (where
$\infty-1=\infty$). We can equivalently view
$\Mod_R$ as an $\EE_{k-1}$-algebra in
$\Prl$ by~\cite{ha}*{Proposition~8.1.2.6} in this case.
This decrease in coherent commutativity is the analogue of the usual fact that
there is no tensor product of right $A$-modules when $A$ is an associative ring.
Thus, by~\cite{ha}*{Corollary~6.3.5.17}, when $2\leq k\leq\infty$ we may build an \icat $\Cat_R$ of (right) $\Mod_R$-modules in $\Prl$.
In the notation of~\cite{ha},
\begin{equation*}
    \Cat_R=\Mod_{\Mod_R}(\Prl).
\end{equation*}
This \icat is $\EE_{k-2}$-monoidal and is closed under small limits and colimits.
Moreover, the $\EE_{k-2}$-monoidal structure is closed
(see~\cite{ha}*{Remark~6.3.1.17} and the beginning of the next section). The dual of $\Cscr$ is
\begin{equation*}
    \Dual_R\Cscr=\Fun_R^{\Lrm}(\Cscr,\Mod_R),
\end{equation*}
the functor category of left adjoint $R$-linear functors from $\Cscr$ to $\Mod_R$.
When $R$ is the sphere spectrum, $\Cat_R$ is also denoted by $\Prlst$; it is the \icat of
\emph{stable} presentable \icats and colimit preserving functors. Since $\Mod_R$ is stable, we could also define
$\Cat_R$ as $\Mod_{\Mod_R}(\Prlst)$.
We will refer to the objects of $\Cat_R$ as $R$-linear categories. An $R$-linear category is thus a
stable \icat with an enrichment in $\Mod_R$: there are functorial $R$-module mapping spectra
$\Map_{\Cscr}(x,y)$ for $x,y$ in $\Cscr$.

We may also consider the \icat $\Prl_{\st,\omega}$ of
compactly generated stable presentable \icats with morphisms the colimit preserving functors that
preserve compact objects. Then, $\Prl_{\st,\omega}$ inherits a symmetric monoidal
structure from $\Prl_\st$, as one can check by using the proof of \cite{ha}*{6.3.1.14} in the
$\omega$-compactly generated situation. The \icat $\Mod_R$ is again an $\EE_{k-1}$-monoid in
$\Prl_{\st,\omega}$, and so we can consider the \icat $\Cat_{R,\omega}$ of compactly
generated $R$-linear categories and colimit preserving functors that preserve
compact objects. The natural map $\Cat_{R,\omega}\rightarrow\Cat_R$ is an
$\EE_{k-2}$-monoidal map of $\infty$-categories.

There is a natural equivalence
\begin{equation*}
    \Ind:\Cat_{\infty}^{\perf}\leftrightarrows\Prl_{\st,\omega}:(-)^{\omega}
\end{equation*}
of symmetric monoidal $\infty$-categories,  where  $\Cat_{\infty}^{\perf}$  is  the  symmetric
monoidal \icat  of  small
idempotent complete stable \icats and exact functors.
One may also view \icat $\Cat_{\infty}^{\perf}$  as  the  localization  of  the  \icat  of  spectrally
enriched categories $\Cat_{\Sp}$ given by inverting the maps
$\Ascr\rightarrow\Mod_{\Ascr}^{\omega}$ for all (compact)  spectral  categories  $\Ascr$.
For details, see~\cite{bgt1}.
If $R$ is an $\EE_k$-ring, this
equivalence sends the $\EE_{k-1}$-algebra $\Mod_R$
to $\Mod_R^{\omega}$ in $\Cat_{\infty}^{\perf}$.
Thus, it induces an equivalence between $\Cat_{R,\omega}$ and
$\Mod_{\Mod_R^\omega}(\Cat_{\infty}^{\perf})$.

In the rest of this section, we prove some technical results relating algebras and their
modules categories, which we will need later in the paper. While the statements are true for
$\EE_k$-ring spectra with $3\leq k\leq\infty$, for simplicity we treat only
$\EE_\infty$-ring spectra. Fix an $\EE_\infty$-ring $R$.
Let
\[
\Mod_*:\Alg_R\rightarrow(\Cat_R)_{\Mod_R/}
\]
be the symmetric monoidal functor which
sends an $R$-algebra $A$ to the $R$-linear category of right $A$-modules $\Mod_A$ with basepoint $A$.
We abuse notation and write $\Mod_A$ for the object $(\Mod_A,A)$ of $(\Cat_R)_{\Mod_R/}$.
There are analogous functors $\Mod_{*,\omega}:\Alg_R\rightarrow(\Cat_{R,\omega})_{\Mod_R/}$, and we
can forget the basepoint to obtain $\Mod:\Alg_R\rightarrow\Cat_R$ and
$\Mod_{\omega}:\Alg_R\rightarrow\Cat_{R,\omega}$.

There is an adjunction
\begin{equation*}
    \Mod_*:\Alg_R\rightleftarrows \left(\Cat_R\right)_{\Mod_R/}:\End
\end{equation*}
where the right adjoint $\End$ takes a pointed $R$-linear category and sends  it  to
the $R$-algebra of endomorphisms of the distinguished object.

\begin{proposition}
    For an $\EE_\infty$-ring $R$,
    the functors $\Mod_*:\Alg_R\rightarrow\left(\Cat_{R}\right)_{\Mod_R/}$
    and $\Mod_{*,\omega}:\Alg_R\rightarrow\left(\Cat_{R,\omega}\right)_{\Mod_R/}$ are  fully  faithful.
    \begin{proof}
        To check the first statement, for $R$-algebras $A$ and $B$, consider the fiber sequence
        \begin{equation*}
            \map_{\Mod_R/}(\Mod_A,\Mod_B)\rightarrow\map_R(\Mod_A,\Mod_B)\rightarrow\map_R(\Mod_R,\Mod_B).
        \end{equation*}
        Since $\Mod_A$ is dualizable with dual  $\Mod_{A^{\op}}$  and using that the symmetric monoidal
        structure on $\Cat_R$ is closed, we can rewrite the fiber sequence as
        \begin{equation*}
            \map_{\Mod_R/}(\Mod_A,\Mod_B)\rightarrow\Mod_{A^{\op}\otimes_R
            B}^{\eq}\rightarrow\Mod_B^{\eq}.
        \end{equation*}
        The fiber of the map over $B$ is equivalent to the space of $A^{\op}\otimes_R B$-module structures
        compatible   with   the   $B$-module   structure   on   $B$,   which    is    simply
        \begin{equation*}
            \map_{\Alg_R}(A,\End_B(B))\we\map_{\Alg_R}(A,B).
        \end{equation*}
        So, the functor is fully faithful.

        To  check  the  second  statement, simply note  that  there   is   a   pullback   square
        \begin{equation*}
            \begin{CD}
                \map_R^{\omega}(\Mod_A,\Mod_B)  @>>>    \map_R^{\omega}(\Mod_R,\Mod_B)\\
                @VVV                                    @VVV\\
                \map_R(\Mod_A,\Mod_B)           @>>>    \map_R(\Mod_R,\Mod_B).
            \end{CD}
        \end{equation*}
        of mapping spaces, so the fibers are equivalent.
    \end{proof}
\end{proposition}

\begin{corollary}\label{cor:algebrafiber}
    If $A$ is an $R$-algebra, then the fiber over a compact $R$-module $P$ of the forgetful map
    \begin{equation*}
        \map_R^{\omega}(\Mod_A,\Mod_R)\rightarrow\Mod_R^{\eq}
    \end{equation*}
    is naturally equivalent to $\map_{\Alg_R}(A,\End_R(P))$.
\end{corollary}

Despite the fact that $\Mod_R$ is the unit of the symmetric monoidal structure on
$\Cat_{R,\omega}$, it is not formal that $\Mod_R$ is a compact object in $\Cat_{R,\omega}$.
The fact that it is compact is essential in deducing that Azumaya algebras are compact
$R$-algebras (and not just compact as $R$-modules). 

\begin{theorem}
    The unit $\Mod_R$ is a compact object of $\Cat_{R,\omega}$.
    \begin{proof}
        We begin by showing that the $\infty$-category of spectra is compact in
        $\Prl_{\st,\omega}$.
        Equivalently, we must show that the functor
        $\map(\Sp,-):\Prl_{\st,\omega}\rightarrow\Sscr$ preserves filtered colimits. Since
        $\Delta^0$ is a compact object of $\Cat_\infty$, the underlying space functor
        $\Cat_\infty\rightarrow\Sscr$ preserves filtered colimits, and
        we see that it is enough to show that
        $$\Fun^{\Lrm,\omega}(\Sscr,-):\Prl_{\st,\omega}\rightarrow\Cat_\infty$$ preserves filtered
        colimits. By~\cite{htt}*{Proposition~5.5.7.11}, the forgetful functor
        $\Cat_{\infty}^{\Rex(\omega)}\rightarrow\Cat_\infty$ preserves filtered colimits,
        where $\Cat_\infty^{\Rex(\omega)}$ denotes the $\infty$-category of finitely
        cocomplete $\infty$-categories and finite colimit-preserving functors. Recall that
        taking compact objects $(-)^\omega$ identifies $\Prl_{\st,\omega}$ with the full
        subcategory $\Cat_\infty^\perf\subseteq\Cat_\infty^{\Rex(\omega)}$
        consisting of the stable and idempotent-complete objects. Moreover, this inclusion admits a
        left adjoint
        \begin{equation*}
            \mathrm{Stab}(\Ind(-))^{\omega}:\Cat_\infty^{\Rex(\omega)}\rightarrow\Cat_\infty^\perf,
        \end{equation*}
        and the functor
        $(-)^{\omega}:\Prl_{\st,\omega}\rightarrow\Cat_\infty^\perf$ admits a
        left-adjoint $\Ind$ given by ind-completion.

        Let $\colim_i\Cscr_i\we\Cscr$ be a filtered colimit in $\Prl_{\st,\omega}$. It follows
        that the canonical map
        $\colim_i\Cscr_i^\omega\rightarrow\Cscr^\omega$ is an
        idempotent completion, so it is fully faithful and any object $P$ in
        $\Cscr^\omega$ is a retract of an object $Q$ in $\colim_i\Cscr_i^\omega$. In
        particular, there is an idempotent $e\in\pi_0\mathrm{end}(Q)$ such that $P$ is the cofiber
        of
        \begin{equation}\label{eq:234234}
            \bigoplus_{k=0}^{\infty} Q\xrightarrow{1-s(e)}\bigoplus_{k=0}^{\infty} Q,
        \end{equation}
        where $s(e)$ is the map which on the $k$th component maps $Q$ to the $(k+1)$st
        component via $e$.
        Since the colimit $\colim_i\Cscr_i^{\omega}$ is computed in $\Cat_\infty$, it follows that $Q$
        is the image of an object $Q_i$ in $\Cscr_i^\omega$ for some $i$. Write $Q_j$ for
        the image of $Q_i$ in $\Cscr_j$. 
        Because mapping spaces in filtered colimits of $\infty$-categories are given as the
        filtered colimit of the mapping spaces, there is a natural equivalence
        \begin{equation*}
            \colim_{j\geq i}\mathrm{end}(Q_j)\we\mathrm{end}(Q).
        \end{equation*}
        It follows that we may lift $e$ to an idempotent $e_j$ of $Q_j$ for some $j\geq i$.
        Define $P_j$ to be the summand of $Q_j$ split off by this idempotent as
        in~\eqref{eq:234234}. Then, $P_j$ is compact object of $\Cscr_j$ which maps to $P$
        in the colimit. It follows that $\colim_i\Cscr_i^{\omega}\rightarrow\Cscr^\omega$ is
        essentially surjective and hence an equivalence.

        To deduce that, in general, $\Mod_R$ is a compact object of $\Cat_{R,\omega}$, it
        suffices to note that the forgetful functor
        $\Cat_{R,\omega}\we\Mod_{\Mod_R}(\Prl_{\st,\omega})\rightarrow\Prl_{\st,\omega}$
        preserves filtered colimits. This follows from~\cite{ha}*{Corollary 3.4.4.6}, which
        is applicable because $\Prl_{\st,\omega}\we\Cat_\infty^\perf$, as a symmetric
        monoidal $\infty$-categories and the symmetric monoidal structure is closed
        by~\cite{bgt1}*{Theorem 2.14}.
    \end{proof}
\end{theorem}

From the theorem, we deduce an important fact about the endomorphism functor.

\begin{lemma}\label{lem:endfilteredcolimits}
    The  right  adjoint  $\End:(\Cat_{R,\omega})_{\Mod_R/}\rightarrow\Alg_R$ of $\Mod_*$ preserves
    filtered colimits.
    \begin{proof}
        A map $\Mod_R\rightarrow\Cscr$ in $\Cat_{R,\omega}$ classifies a compact object of
        $\Cscr$, i.e., a pointed $R$-linear category. Let $\colim_i\Cscr_i\we\Cscr$ be a
        colimit of pointed compactly generated $R$-linear categories. Let $X_i$ be the image of $R$ in $\Cscr_i$, and let $X$ be the image of
        $R$ in $\Cscr$. Consider the map of $R$-algebras
        \[
        \colim_{i}\End_{\Cscr_i}(X_i)\longrightarrow\End_\Cscr(X).
        \]
        Since the forgetful functors
        \begin{equation*}
            \Alg_R\rightarrow\Mod_R\rightarrow\Sp\xrightarrow{\li\Sigma^n}\Spaces
        \end{equation*}
        preserve filtered colimits, and taken together, they  detect  filtered  colimits  in
        $\Alg_R$, it is enough to show that
        \begin{equation*}
            \colim_i\mathrm{end}_{\Cscr_i^\omega}(X_i)\rightarrow\mathrm{end}_{\Cscr^\omega}(X)
        \end{equation*}
        is an equivalence. This follows because we know that the filtered colimit of pointed
        compactly generated $R$-linear categories agrees with the filtered colimit as
        compactly generated $R$-linear categories,
        with the obvious basepoint, and, by the theorem,
        $\colim_i\Cscr_i^\omega\we\Cscr^\omega$ in $\Cat_\infty$.
    \end{proof}
\end{lemma}

We now prove the important fact that compactness of an $R$-algebra $A$ is detected purely through
the module category of $A$.

\begin{proposition}\label{prop:Acpt}
    Let $A$ be an $R$-algebra. Then, $A$ is compact in $\Alg_R$ if and only if $\Mod_A$ is
    compact in $\Cat_{R,\omega}$.
    \begin{proof}
        Assume first that $A$ is compact in $\Alg_R$, and let $\Cscr$ be a
        filtered  colimit of a diagram $\{\Cscr_i\}_{i\in I}$ in $\Cat_{R,\omega}$.
        Because   $\End$    preserves filtered    colimits by the previous lemma,    it    is    clear    that
        $\Mod_*:\Alg_R\rightarrow(\Cat_{R,\omega})_{\Mod_R/}$ preserves compact objects.
        Using the fact that $\Mod_R$ is a compact
        object in $\Cat_{R,\omega}$, and that $A$ is a compact
        $R$-algebra, for every object $M\in\colim_i \map(\Mod_R,\Cscr_i)\simeq\map(\Mod_R,\Cscr)$,
        which is to say an object $M_i$ of $\Cscr_i$ for $i$ sufficiently large, the    adjunction    provides    a    map    of    fiber    sequences
        \begin{equation}\label{eq:12312}
            \begin{CD}
                \colim_i\map_{\Alg_R}(A,\End_{\Cscr_i}(M_i))   @>>>    \colim_i\map(\Mod_A,\Cscr_i)    @>>>
                \colim_i\map(\Mod_R,\Cscr_i)\\
                @VVV                                       @VVV @VVV\\
                \map_{\Alg_R}(A,\End_{\Cscr}(M))             @>>>    \map(\Mod_A,\Cscr)          @>>>
                \map(\Mod_R,\Cscr)
            \end{CD}
        \end{equation}
        in which the left and right vertical maps are equivalences (note that the top sequence is
        a fiber sequence because filtered colimits commute with finite
        limits, and hence, in particular, fiber sequences
        by~\cite{htt}*{Proposition~5.3.3.3}); it  follows  that  the middle map is an equivalence,
        which shows that $\Mod_A$ is compact in $\Cat_{R,\omega}$.

        Now, assume that $\Mod_A$ is compact in $\Cat_{R,\omega}$. Using~\eqref{eq:12312} and the adjunction
        \[
            \map_{{\Mod_{R}}/}(\Mod_A,\Cscr)\simeq\map_{\Alg_R}(A,\End_{\Cscr}(M)),
        \]
        it is easy to see that $\Mod_A$, with basepoint $A$, is also compact in
        $(\Cat_{R,\omega})_{\Mod_R/}$. Let  $B=\colim  B_i$  be  a
        filtered colimit of $R$-algebras. Then, there are equivalences,
        \begin{align*}
            \colim\map_{\Alg_R}(A,B_i)   &\we    \colim\map_{\Alg_R}(A,\End_{\Mod_{B_i}}(B_i))\\
                                &\we    \colim\map_{\Mod_R/}(\Mod_A,\Mod_{B_i})\\
                                &\we    \map_{\Mod_R/}(\Mod_A,\Mod_B)\\
                                &\we    \map_{\Alg_R}(A,B).
        \end{align*}
        That is, $A$ is a compact object in $\Alg_R$.
    \end{proof}
\end{proposition}

\begin{corollary}
    Compactness is a Morita-invariant property of $R$-algebras.
\end{corollary}

\subsection{Smooth and proper algebras}

If   $\Cscr$   is   an   object    of    $\Cat_{R}$,    then    the    dual    of    $\Cscr$
is the functor category
\begin{equation*}
    \Dual_R\Cscr=\Fun_R^{\Lrm}(\Cscr,\Mod_R)
\end{equation*}
in $\Mod_R$. There is a functorial evaluation map
\begin{equation*}
    \Cscr\otimes_R \Dual_R\Cscr\xrightarrow{\ev}\Mod_R.
\end{equation*}
The $R$-linear category $\Cscr$ is dualizable if there exists
a coevaluation map
\begin{equation*}
    \Mod_R\xrightarrow{\coev}\Dual_R\Cscr\otimes_R\Cscr,
\end{equation*}
which classifies $\Cscr$ as a $\Dual_R\Cscr\otimes_R\Cscr$-module, such that both maps
\begin{equation*}
    \Cscr\xrightarrow{\Cscr\otimes_R\coev}\Cscr\otimes_R\Dual_R\Cscr\otimes_R\Cscr\xrightarrow{\ev\otimes_R\Cscr}\Cscr
\end{equation*}
and
\begin{equation*}
    \Dual_R\Cscr\xrightarrow{\coev\otimes\Dual_R\Cscr}\Dual_R\Cscr\otimes_R\Cscr\otimes_R\Dual_R\Cscr\xrightarrow{\Dual_R\Cscr\otimes\ev}\Dual_R\Cscr
\end{equation*}
are equivalent to the identity.

\begin{lemma}
    An object $\Cscr$ is dualizable in $\Cat_{R,\omega}$ if and only if it is dualizable in
    $\Cat_R$ and the evaluation and
    coevaluation morphisms of its duality data in $\Cat_R$ are morphisms in
    $\Cat_{R,\omega}$.
    \begin{proof}
        Dualizability is detected on the monoidal homotopy category, and the duality data for
        $\Cscr$ in $\Ho(\Cat_{R,\omega})$ must coincide with the duality data in
        $\Ho(\Cat_R)$ by uniqueness.
    \end{proof}
\end{lemma}

\begin{definition}
    A compactly generated $R$-linear category $\Cscr$ is proper if its
    evaluation   map   is   in   $\Cat_{R,\omega}$;
    it is smooth if it is dualizable and its  coevaluation map is in $\Cat_{R,\omega}$.
\end{definition}

If $A$ is an $R$-algebra, then $\Mod_A$ is proper if and only if $A$ is a perfect
$R$-module. Indeed, in this case, the evaluation map is the map
\begin{equation*}
    \Mod_{A\otimes_R A^{\op}}\we\Mod_A\otimes_R\Mod_{A^{\op}}\rightarrow\Mod_R
\end{equation*}
that sends $A\otimes_R A^{\op}$ to $A$. We say in this case that $A$ is a proper $R$-algebra. Similarly,
$\Mod_A$ is smooth if and only if the coevaluation map
$\Mod_R\rightarrow\Mod_{A^{\op}\otimes_R A}$, which sends $R$ to $A$, considered as an
$A^{\op}\otimes_R A$-module, exists and is in
$\Cat_{R,\omega}$. So we see that $\Mod_A$ is smooth if and only
if $A$ is perfect as an $A^{\op}\otimes_R A$-module. Again, we say in this case that $A$ is
a smooth $R$-algebra. In fact, every smooth $R$-linear category is equivalent to a module
category.

\begin{lemma}
    Suppose that $\Cscr$ is a smooth $R$-linear category. Then, $\Cscr\we\Mod_A$ for some
    $R$-algebra $A$.
    \begin{proof}
        See~\cite{toen-derived}*{Lemma 2.6}. The morphism
        \begin{equation*}
            \Mod_R\xrightarrow{\coev}\Dual_R\Cscr\otimes_R\Cscr
        \end{equation*}
        is in $\Cat_{R,\omega}$ by hypothesis. Thus, $R$ is sent by $\coev$ to a compact
        object of $\Dual_R\Cscr\otimes_R\Cscr$. The compact objects of this category are
        the smallest idempotent complete stable subcategory of $\Dual_R\Cscr\otimes_R\Cscr$
        containing the objects of the form
        $\Map_\Cscr(a,-)\otimes_R b$, where $a$ and $b$ are compact objects of $\Cscr$. This
        is because $\Cscr$ is compactly generated, so the compact objects of $\Dual_R\Cscr$
        are precisely the duals of the compact objects of $\Cscr$. We
        can thus write $\coev(R)$ as the result of taking finitely many shifts, cones, and
        summands of $\Map_\Cscr(a_i,-)\otimes_R b_i$, for $i=1,\ldots,n$. The identity map
        \begin{equation*}
            \Cscr\xrightarrow{\Cscr\otimes_R\coev}\Cscr\otimes_R\Dual_R\Cscr\otimes_R\Cscr\xrightarrow{\ev\otimes_R\Cscr}\Cscr
        \end{equation*}
        sends $c\in\Cscr$ to the same diagram built out of $\Map_\Cscr(a_i,c)\otimes_R b_i$.
        It follows that if $\Map_\Cscr(a_i,c)\we 0$ for $i=1,\ldots,n$, then $c\we 0$. Thus,
        the $a_i$ form a set of compact generators for $\Cscr$. Letting
        $A=\End_{\Cscr}(\oplus_i a_i)^{\op}$, it follows that $\Cscr\we\Mod_A$, as desired.
    \end{proof}
\end{lemma}

\begin{definition}
    An $R$-linear category is of finite type if there exists a compact
    $R$-algebra $A$ such that $\Cscr$ is equivalent to $\Mod_A$.
\end{definition}

The condition of being smooth and proper is a strong one for $R$-algebras: it implies
compactness in the $\infty$-category of $R$-algebras.
See~\cite{toen-vaquie}*{Corollary~2.13} for the dg-statement.

\begin{proposition}
    If $\Cscr$ is a smooth and proper $R$-linear category, then $\Cscr$ is of  finite  type.
    \begin{proof}
        Let $A$ be an $R$-algebra such that $\Cscr\we\Mod_A$.
        To show that $A$ is compact as an $R$-algebra, it suffices by Proposition \ref{prop:Acpt} to show that $\Cscr$ is compact in $\Cat_{R,\omega}$.
        To this end, fix a filtered colimit $\Dscr=\colim_{i\in I}\Dscr_i$ in $\Cat_{R,\omega}$; we must show that the natural map
        \begin{equation*}
            \colim_{i\in I}\map_{\Cat_{R,\omega}}(\Cscr,\Dscr_i)\rightarrow\map_{\Cat_{R,\omega}}(\Cscr,\Dscr)
        \end{equation*}
        is an equivalence. The dualizability of $\Cscr$ in $\Cat_{R,\omega}$ gives natural equivalences
        \begin{equation}\label{eq1}
            \map_{\Cat_{R,\omega}}(\Cscr,\Dscr)\we\map_{\Cat_{R,\omega}}(\Mod_{R},D_R\Cscr\otimes_R\Dscr),
        \end{equation}
        and as $\Mod_R$ is compact as a compactly generated $R$-linear category, the result follows from the equivalences
        \[
        \colim_i D_R\Cscr\otimes_R\Dscr_i\we D_R\Cscr\otimes_R\colim_i\Dscr_i\we
        D_R\Cscr\otimes_R\Dscr.\qedhere
        \] 
    \end{proof}
\end{proposition}

The following result is due to To\"en and Vaqui\'e~\cite{toen-vaquie} in the dg-setting, and
the arguments are essentially the same. The result is part of the philosophy of hidden
smoothness due to Kontsevich.

\begin{theorem}
    An $R$-linear category of finite type is smooth.
    \begin{proof}
        It suffices to show that if $A$ is a compact $R$-algebra, then it is perfect as a right $A^{\op}\otimes_R A$-module.
        There is a fiber sequence
        \begin{equation*}
            \Omega_{A/R}\rightarrow A^{\op}\otimes_R A\rightarrow A,
        \end{equation*}
        where $\Omega_{A/R}$ is the $A^{\op}\otimes_R A$-module of differentials (see Lazarev~\cite{lazarev}).
        So, it is enough to show that $\Omega_{A/R}$ is a perfect $A^{\op}\otimes_R
        A$-module when $A$ is a compact $R$-algebra.
        This follows from the adjunction
        \begin{equation*}
            \map_{A^{\op}\otimes_R  A}(\Omega_{A/R},M)\we\map_{(\Alg_R)_{/A}}(A,A\oplus  M),
        \end{equation*}
        together with the fact that, since $A$ is  a compact  $R$-algebra,  $A$  is
        compact in $(\Alg_R)_{/A}$.
    \end{proof}
\end{theorem}

\subsection{Azumaya algebras}

Let  $R$  be   an   $\EE_\infty$-ring   spectrum.    The   following   definition   is   due
to Auslander and Goldman~\cite{auslander-goldman}.  In the derived setting, it
and variations on it have been  considered  by  Lieblich~\cite{lieblich-thesis}, Baker-Lazarev~\cite{baker-lazarev}, To\"en~\cite{toen-derived},
Johnson~\cite{johnson-azumaya}, and Baker-Richter-Szymik~\cite{baker-richter-szymik}.

\begin{definition}
    An $R$-algebra $A$ is an Azumaya $R$-algebra if $A$ is a compact generator of $\Mod_R$ and if
    the   natural   $R$-algebra   map    giving    the    bimodule    structure    on    $A$
    \begin{equation*}
        A\otimes_R A^{\op}\rightarrow\End_{R}(A)
    \end{equation*}
    is an equivalence of $R$-algebras.
\end{definition}

Note that if $A$ is an Azumaya $R$-algebra, then, by definition, $A\otimes_R A^{\op}$ is
Morita equivalent to $R$.
The standard example of an Azumaya algebra is the endomorphism algebra $\End_R(P)$ of a
compact generator of $\Mod_R$. These algebras are not so interesting as they are
already Morita equivalent to $R$. The Brauer group will be the group of Morita equivalence
classes of Azumaya algebras, so these endomorphism algebras will represent the trivial
class. For more examples and various properties, we refer
to~\cite{baker-richter-szymik}. In particular, we will use the fact that if $S$ is an
$\EE_\infty$-$R$-algebra, then $A\otimes_R S$ is Azumaya if $A$ is~\cite{baker-richter-szymik}*{Proposition~1.5}. One main goal of this
paper is to show that if $R$ is a connective commutative ring spectrum, then Azumaya algebras
are \'etale locally Morita equivalent to $R$, which To\"en established in the connective
commutative dg-setting~\cite{toen-derived}.
The first fact we need is the following theorem.

\begin{theorem}[To\"en~\cite{toen-derived}]\label{thm:brk}
    If $R=\Hrm k$, where $k$ is an algebraically closed field, then every Azumaya $R$-algebra is Morita
    equivalent to $R$.
\end{theorem}

We prove now a characterization of Azumaya algebras and smooth and proper algebras. The
corresponding statement for dg-algebras is~\cite{toen-derived}*{Proposition~2.5}.

\begin{theorem}\label{thm:azumayaiffinvertible}
    Let $\Cscr$ be a compactly generated $R$-linear category.
    \begin{enumerate}
        \item[\rm{(1)}] $\Cscr$ is dualizable in $\Cat_{R,\omega}$ if and only if $\Cscr$ is
            equivalent  to   $\Mod_A$   for   a   smooth   and   proper   $R$-algebra   $A$.
        \item[\rm{(2)}] $\Cscr$ is invertible in $\Cat_{R,\omega}$ if and only if $\Cscr$ is
            equivalent to $\Mod_A$ for an Azumaya $R$-algebra $A$.
    \end{enumerate}
    \begin{proof}
        If $A$ is smooth and proper, then $\Mod_A$ is dualizable in $\Cat_{R,\omega}$ since
        the evaluation and coevaluation maps are in $\Cat_{R,\omega}$ by hypothesis. If
        $\Cscr$ is smooth and proper, then $\Cscr\we\Mod_A$ for an $R$-algebra $A$ which is,
        by definition, smooth and proper.

        Suppose that $\Cscr$
        is invertible. Then, it follows that it is dualizable in $\Cat_{R,\omega}$, and thus that
        it is equivalent to $\Mod_A$ where $A$ is a smooth and proper $R$-algebra.   So,  it
        suffices to show that $\Mod_A$ is invertible if and only if  $A$  is  Azumaya.   The
        evaluation map
        \begin{equation*}
            \Mod_{A\otimes_R A^{\op}}\rightarrow\Mod_R
        \end{equation*}
        is an equivalence if and only if $A$ is invertible. This map sends $A\otimes_R
        A^{\op}$ to $A$, and it is contained in $\Cat_{R,\omega}$ if and only if $A$ is a
        compact $R$-module.  The evaluation map is essentially surjective if and only if $A$ is a
        generator   of   $\Mod_R$. Finally,   it   is   fully    faithful    if    and    only    if
        \begin{equation*}
            A\otimes_R A^{\op}\we\End_{A\otimes_R A^{\op}}(A\otimes_R A^{\op})\rightarrow\End_R(A)
        \end{equation*}
        is an equivalence.
    \end{proof}
\end{theorem}

We see that we might define the Brauer space of an $\EE_\infty$-ring spectrum $R$ to
be the grouplike $\EE_{\infty}$-space $\Cat_{R,\omega}^{\times}$. Instead, we will later
give an equivalent definition that generalizes more readily to derived schemes.

\section{Sheaves}\label{sec:sheaves}

We give in this section preliminaries we will need about sheaves of spaces and
$\infty$-categories. In particular, we study smoothness for morphisms of sheaves of spaces,
and we show that under mild hypotheses smooth surjective morphisms admit \'etale local sections.

\subsection{Stacks of algebra and module categories}

Roughly speaking, if $\Xscr$ is an $\infty$-topos and $\Cscr$ is a complete
$\infty$-category, then a $\Cscr$-valued sheaf on $\Xscr$ is a functor
$\Xscr^{\op}\rightarrow\Cscr$ which satisfies descent.

\begin{definition}
    Let $\Cscr$ be a complete $\infty$-category. A $\Cscr$-valued sheaf on $\Xscr$ is a limit-preserving functor
    $\Xscr^{\op}\rightarrow\Cscr$.  The \icat $\Shv_\Cscr(\Xscr)$ is the  full
    subcategory of $\Fun(\Xscr^{\op},\Cscr)$ consisting of the $\Cscr$-valued sheaves on $\Xscr$.
\end{definition}

In the cases we care about, $\Xscr$ will be the $\infty$-topos associated to a Grothendieck
topology on an $\infty$-category $\Ascr$. In this case a $\Cscr$-valued sheaf on $\Xscr$ is
determined by its values on $\Ascr$, because every object in $\Xscr$ is a colimit of
representable functors. Moreover, we will typically be in an even more special situation,
where the Grothendieck topology is given by a pre-topology satisfying the conditions
of~\cite{dag7}*{Propositions 5.1 and 5.7}. In this case, a functor $F:\Ascr^{\op}\rightarrow\Cscr$ is a sheaf if
and only if for every covering morphism $X\rightarrow Y$ in $\Ascr$, the map
\begin{equation*}
    F(Y)\rightarrow\lim_{\Delta}F(X_{\bullet})
\end{equation*}
is an equivalence in $\Dscr$, where $X_{\bullet}$ is the simplicial object associated to
the cover. Similarly, $F$ is a hypercomplete
sheaf, or hypersheaf, if for every hypercovering $V_{\bullet}\rightarrow Y$ in $\Ascr$, the map
\begin{equation*}
    F(Y)\rightarrow\lim_{\Delta}F(V_{\bullet})
\end{equation*}
is an equivalence. See~\cite{dag7}*{Section 5} for details. In particular, Lurie proves that
the collection of faithfully flat morphisms in $(\CAlg_R)^{\op}$ satisfies the necessary
conditions. Thus, the collection of faithfully flat \'etale morphisms
(Section~\ref{sub:topologies}) in $(\CAlg_R^{\cn})^{\op}$ does as well.

In practice, our sheaves will be one of the following three types: sheaves of
$\infty$-groupoids (spaces), which we call sheaves; sheaves of spectra; or, sheaves of
(not necessarily small) $\infty$-categories, which we call stacks. Thus, for instance, a stack on an $\infty$-topos $\Xscr$ is a
limit-preserving functor $\Xscr^{\op}\rightarrow\LCI$. We will also consider sheaves of ring
spectra and stacks of symmetric monoidal $\infty$-categories. A presheaf of symmetric monoidal
$\infty$-categories is a stack if and only if the underlying presheaf of
$\infty$-categories is a stack. Indeed,
the  forgetful  functor  $\CAlg(\LCI)\rightarrow\LCI$ preserves and detects
limits~\cite{ha}*{Corollary~3.2.2.5}.

The conventions spelled out in the previous paragraph might cause some confusion. We have
chosen to emphasize the $\infty$-categorical notion that groupoids are spaces in our
definitions. As a result, we end up saying
``sheaf of Morita equivalences,'' ``classifying sheaf,'' or ``Deligne-Mumford sheaf,''
instead of the more comfortable ``stack of Morita equivalences,'' ``classifying stack,'' and
``Deligne-Mumford stack.'' Our stacks will be sheaves of $\infty$-categories. This approach
is justified by the fact that the three examples just given are actually objects of the
underlying $\infty$-topoi. Since the objects of the $\infty$-topos themselves are sheaves of
spaces, there is no longer any need to have a separate notion of a sheaf of groupoids.

From a stack, we can produce a sheaf of (not necessarily small) spaces as follows. There is a pair of adjoint functors
\begin{equation*}
    i:\Gpd\leftrightarrows\LCI:(-)^{\eq},
\end{equation*}
where the left adjoint $i$ is the natural inclusion, and $(-)^{\eq}$ sends an $\infty$-category $\Cscr$ to its maximal
subgroupoid $\Cscr^{\eq}$.
If $\StM:\Xscr^{\op}\rightarrow\LCI$ is a  stack,  then  the  associated  sheaf
$\StM^{\eq}$ is the composition of $\StM$ with $(-)^{\eq}$, which is a sheaf because
$(-)^{\eq}$ preserves limits.

In the remainder of the section, we will recall some  facts  about  \'etale  (hyper)descent.
Let $R$ be a connective $\Einfinity$-ring, and let $\Shv_R^{\et}$ denote the big \'etale
$\infty$-topos on $R$. Given any commutative $R$-algebra $U$, connective or not, there is a
presheaf $X=\Spec U$ whose values on an $R$-algebra $S$ are given by
\begin{equation*}
    X(S)=\map_{\CAlg_R}(U,S).
\end{equation*}
This presheaf is in fact a sheaf, which says that the \'etale topology on
$\Aff_R^{\cn}$ is subcanonical, though much more is true~\cite{dag7}*{Theorem 5.14}.

\begin{proposition}\label{prop:subcanonical}
    For any commutative $R$-algebra $U$, the presheaf $\Spec U$ is an \'etale hypersheaf.
    \begin{proof}
        Indeed, let $S\rightarrow T^{\bullet}$
        be an \'etale hypercovering. This determines a map $\Nrm(\Delta_{+})\rightarrow\CAlg_R$,
        which is a limit diagram by~\cite{dag7}*{Lemma~5.13}.\qedhere
    \end{proof}
\end{proposition}

Let $\StMod:\left(\Aff^{\cn}_R\right)^{\op}=\CAlg_R^{\cn}\rightarrow\CAlg(\Prl)$ be the presheaf of
symmetric   monoidal   \icats   that   sends   $S$    to    $\Mod_S$.
By~\cite{dag7}*{Theorem~6.1}, this presheaf satisfies descent for \'etale hypercovers. It
follows that we may uniquely extend $\StMod$ to a hyperstack on all of $\Shv_R^{\et}$.
Concretely, when $X$ is an object of $\Shv_R^{\et}$, we let
\begin{equation*}
    \Mod_X=\lim_{\Spec S\rightarrow X}\Mod_S
\end{equation*}
the stable presentable symmetric monoidal $\infty$-category of modules over $X$.
We are actually keeping track of the symmetric monoidal structure on $\Mod_S$, and hence on
$\Mod_X$ by forming the limit in the $\infty$-category $\CAlg(\Prl)$. However, the forgetful
functor $\CAlg(\Prl)\rightarrow\Prl$ preserves limits, so we choose to ignore the
intricacies of symmetric monoidal $\infty$-categories and suppress the symmetric monoidal
structure from the notation.

By composing $\StMod:\Shv_R^{\et}\rightarrow\CAlg(\Prl)$ with the limit-preserving functor
$$\Alg:\CAlg(\Prl)\rightarrow\CAlg(\Prl)$$ that sends a presentable symmetric monoidal
$\infty$-category to $\infty$-category of algebra objects (which is also presentable
by~\cite{ha}*{Corollary 3.2.3.5} and
symmetric monoidal by~\cite{ha}*{Proposition 3.2.4.3 and Example 3.2.4.4}),
we obtain the hyperstack of algebras $\StAlg$ on $\Shv_R^{\et}$. There is a
substack $\StAz$ of Azumaya algebras: an algebra $\Ascr$ over $X$ is Azumaya if its
restriction to any affine scheme is Azumaya.

Recall that if $\Cscr$ is a symmetric monoidal $\infty$-category, then its space of units
$\Pic(\Cscr)$ is the grouplike $\Einfinity$-space consisting of invertible elements of
$\Cscr$ and equivalences. When $\Cscr$ is presentable, then $\Pic(\Cscr)$ is a small space,
as proven in~\cite{abg}*{Theorem~8.9}. Thus, there is a functor
\begin{equation*}
    \Pic:\CAlg(\Prl)\rightarrow\CAlg^{\gp}(\Sscr),
\end{equation*}
where $\CAlg^{\gp}(\Sscr)$ denotes the full subcategory of $\CAlg(\Sscr)$ of grouplike
$\EE_\infty$-spaces.

\begin{proposition}\label{cor:unitsaresheaves}
    If $\StM$ is a hyperstack of presentable symmetric monoidal $\infty$-categories, then
    the presheaf $\Pic(\StM)$ is a hypersheaf.
    \begin{proof}
        By~\cite{abg}*{Theorem~8.10}, $\Pic$ is a right adjoint, so it preserves limits.
    \end{proof}
\end{proposition}

Applying the lemma to the particular stack $\StMod$ on $\Shv_R^{\et}$, we obtain the Picard
sheaf $\ShPic$, and we let $\Shpic$ be the associated sheaf of spectra.

Now, we introduce a stack of $R$-linear categories, $\StCat_R^{\desc}$, which classifies
$R$-linear categories satisfying \'etale hyperdescent.
Let $\StCat_R:(\Aff_R^{\cn})^{\op}=\CAlg_R^{\cn}\rightarrow\LCI$ be the composite functor
\[
\StCat_R:(\Aff_R^{\cn})^{\op}=\CAlg_R^{\cn}\overset{\Mod}{\to}\CAlg(\Prl)\overset{\Mod}{\to}\LCI
\]
whose value at $S$ is the \icat $\Cat_S$ of $S$-linear $\infty$-categories (equivalently, $\Mod_S$-modules in the symmetric monoidal $\infty$-category $\Prl$).

Say that an $R$-linear category $\Cscr$ satisfies \'etale hyperdescent if for each connective commutative $R$-algebra $S$
and each \'etale hypercover $S\to T^\bullet$, the canonical map
\begin{equation*}
    \mathscr{C}\otimes_R    S\rightarrow\lim_{\Delta}\mathscr{C}\otimes_R    T^{\bullet}
\end{equation*}
is an equivalence. We write $\Cat_S^\desc\subseteq\Cat_S$ for the full subcategory of
$\Cat_S$ consisting of the $S$-linear \icats with \'etale hyperdescent and
$\StCat_R^{\desc}\subseteq\StCat_R$ for the full subfunctor of $R$-linear categories with
\'etale hyperdescent.

\begin{example}\label{ex:moddescent}
    If $A$ is an $R$-algebra, then $\Mod_A$ is an $R$-linear category that satisfies \'etale
    hyperdescent.
    Indeed, in this case $\Mod_A$ is dualizable in $\Cat_R$ with dual $\Mod_{A^{\op}}$.
    Therefore, if $S$ is a connective $\Einfinity$-$R$-algebra, then
    \begin{equation*}
        \Mod_A\otimes_R S\we\Dual_R\Dual_R\Mod_A\otimes_R S\we\Fun_R(\Dual_R\Mod_A,\Mod_S)
    \end{equation*}
    in $\Cat_R$. Because functors out of $\Dual_R\Mod_A$ commutes with limits, $\Mod_A$ is an $R$-linear category with hyperdescent.
\end{example}

The important fact about $\StCat_R^{\desc}$ that we need is that it satisfies \'etale
hyperdescent itself.

\begin{proposition}\label{prop:catrstack}
    The   functor   $\StCat_R^{\desc}$   is   an   \'etale   hyperstack    on    $\Aff_R^{\cn}$.
    \begin{proof}
        Lurie proves in~\cite{dag11}*{Theorem~7.5} that the prestack of $R$-linear
        categories satisfying flat hyperdescent is a flat hyperstack.
        The same proof works here.
    \end{proof}
\end{proposition}

\subsection{The cotangent complex and formal smoothness}

We consider notions of smoothness for maps $p:X\rightarrow Y$ of sheaves in $\Shv_R^{\et}$.
References for this material include~\cite{toen-vaquie}, \cite{dag0}, and \cite{dag14}.

\begin{definition}
    Let $p:X\rightarrow Y$ be a map of sheaves. Then, for any point $x\in X(S)$ and any connective $S$-module $M$, the space
    of derivations $\der_{p}(x,M)$ is the fiber of the canonical map
    \[
        X(S\oplus M)\rightarrow X(S)\times_{Y(S)}Y(S\oplus M)
    \]
    over the point corresponding to $x$ and the map $\Spec S\oplus M\rightarrow\Spec S\xrightarrow{x}X\rightarrow Y$,
    where the first map is induced by the map $(id,0):S\rightarrow S\oplus M$. If $Y\we\Spec R$ is a terminal object,
    write $\der_X(-,-)$ for $\der_{p}(-,-)$.
\end{definition}

\begin{definition}
    Let $p:X\rightarrow Y$ be a map of sheaves. An
    object $\Lrm$ of $\Mod_X$ is a relative cotangent complex for  $p$  if  there  exist equivalences
    \begin{equation*}
        \map_{S}(x^*\Lrm,M)\we \der_{p}(x,M)
    \end{equation*}
    natural in $x$ and connective modules $M$.
    When $\Lrm$ exists and is unique up to equivalence then we write $\Lrm_p$ and refer
    to this as \emph{the} cotangent complex of $p$.
    We will often abuse notation and write $\Lrm_{X/Y}$ for $\Lrm_p$ when no confusion will result. When
    $Y\we\Spec R$ is a  terminal  object,  we  write  $\Lrm_X$  in  place  of  $\Lrm_{X/Y}$.
\end{definition}

Note that if $\Lrm$ is a cotangent complex for $p$, then the space of derivations
$\der_p(x,M)$ is never empty.
If $S$ is a ring spectrum, an $S$-module
$M$ is almost connective if it is $k$-connective for some integer $k$. If $X$ is a sheaf, an object $M$ of $\Mod_X$
is almost connective, if its restriction to any $x:\Spec S\rightarrow X$ is almost
connective.

\begin{lemma}
    If $p:X\rightarrow Y$  has  at  least  one  cotangent  complex  $\Lrm$  that  is  almost
    connective,  then  all  cotangent  complexes  are  equivalent,   so   $\Lrm_p$   exists.
    \begin{proof}
        Suppose that $\Lrm$ and $\Lrm'$ are two cotangent complexes for $p$, and suppose
        that $\Lrm$ is almost connective. Let $x:\Spec S\rightarrow X$ be an $S$-point. We
        show that there is an equivalence $x^*\Lrm'\rightarrow x^*\Lrm$, natural in $x$.
        Suppose that $\Sigma^n x^*\Lrm$ is connective. Then, using the chain of equivalences
        \begin{align*}
            \map_S(x^*\Lrm,x^*\Lrm')    &\we\Omega^n\map_S(x^*\Lrm,\Sigma^n x^*\Lrm)\\
                                        &\we\Omega^n\der_p(x,\Sigma^n x^*\Lrm)\\
                                        &\we\Omega^n\map_S(x^*\Lrm',\Sigma^n x^*\Lrm)\\
                                        &\we\map_S(\Sigma^n x^*\Lrm',\Sigma^n x^*\Lrm),
        \end{align*}
        we see there is a unique map $x^*\Lrm'\rightarrow x^*\Lrm$ corresponding to the identity on
        $x^*\Lrm$, which does not depend on $n$, and so is natural in $x$. If there exists
        an integer $k$ such that $\pi_k x^*\Lrm'\rightarrow\pi_k x^*\Lrm$ is not an
        isomorphism, then $\Sigma^k x^*\Lrm'\rightarrow\Sigma^k x^*\Lrm$ is not an
        equivalence. So,
        \begin{equation*}
            \Omega^k\map_S(x^*\Lrm,M)\rightarrow \Omega^k\map_S(x^*\Lrm',M)
        \end{equation*}
        is not an equivalence, which is a contradiction. Thus, $\Lrm'\rightarrow\Lrm$ is an
        equivalence.
    \end{proof}
\end{lemma}

Monomorphisms of sheaves always have cotangent complexes, which vanish.

\begin{lemma}\label{lem:monocotan}
    Let $f:X\rightarrow Y$ be a monomorphism of sheaves. Then, $\Lrm_f\we 0$.
    \begin{proof}
        Suppose that $x:\Spec S\rightarrow X$ is a point and that $M$ is an $S$-module, and consider the diagram
        \begin{equation*}
            \xymatrix{
                X(S\oplus M)    \ar[r] & X(S)\times_{Y(S)}Y(S\oplus M)   \ar[r]\ar[d] &
                Y(S\oplus M)\ar[d]\\
                & X(S) \ar[r] &                          Y(S)},
        \end{equation*}
        in which the square is a pullback. The bottom horizontal arrow is a monomorphism.
        So, the map $X(S)\times_{Y(S)}Y(S\oplus M)\rightarrow Y(S\oplus M)$ is mono. The
        composite $X(S\oplus M)\rightarrow Y(S\oplus M)$ is also mono. So, the map
        $X(S\oplus M)\rightarrow X(S)\times_{Y(S)}Y(S\oplus M)$ is mono, and hence the fibers are
        either empty or contractible. But, the space of derivations
        \begin{equation*}
            \der_f(x,M)
        \end{equation*}
        is the fiber over $x:\Spec S\rightarrow X$ and $\Spec S\oplus M\rightarrow Y$, with
        the latter induced by the composition $\Spec S\oplus M\rightarrow\Spec
        S\xrightarrow{x} X\xrightarrow{f} Y$. It follows that the composite $\Spec S\oplus M\rightarrow X$
        is in the fiber, so it is contractible. Hence, $0$ corepresents derivations.
    \end{proof}
\end{lemma}

The following two lemmas can be proved with
straightforward arguments using only the definition of the space of derivations.

\begin{lemma}\label{lem:cotansequence}
    If $f:X\rightarrow Y$ is a map of sheaves, and if  $\Lrm_X$  and  $\Lrm_Y$  exist,  then
    there is a cofiber sequence
    \begin{equation*}
        f^*\Lrm_Y\rightarrow\Lrm_X\rightarrow\Lrm_f
    \end{equation*}
    in $\Mod_X$. In particular the cotangent complex of $f$ exists.
\end{lemma}

\begin{lemma}\label{lem:colimcotan}
    Let $\{X_i\}$ be a diagram of sheaves in $\Shv_R$ indexed by a simplicial set  $I$,  and
    let $X$ be the limit.  Suppose that the cotangent complex $\Lrm_{X_i}$ exists  for  each
    $i$ in $I$, and write $\Lrm_X$ for the colimit of the diagram $\{\Lrm_{X_i}|_X\}$ in $\Mod_X$.
    If  $\Lrm_X$  is   almost   connective,   then   $\Lrm_X$   is   a   cotangent   complex
    for $X$.
\end{lemma}

The inclusion functor $\tau_{\leq n}\CAlg_R^{\cn}\rightarrow\CAlg_R^{\cn}$ induces a functor
\begin{equation*}
    \tau_{\leq n}^*:\Shv_R^{\et}\we\Shv^{\et}(\CAlg_R^{\cn})\rightarrow\Shv^{\et}(\tau_{\leq n}\CAlg_R^{\cn}).
\end{equation*}
If $S$ is a connective commutative $R$-algebra, then $\tau_{\leq n}^*\Spec S\we\Spec\tau_{\leq n}S$. Indeed, if $T$
is any $n$-truncated connective commutative $R$-algebra, then the natural map
\begin{equation*}
    \map(\tau_{\leq n}S,T)\rightarrow\map(S,T)
\end{equation*}
is an equivalence.

\begin{lemma}\label{lem:truncatedcots}
    If $f:X\rightarrow Y$ is a morphism of sheaves with a cotangent complex $\Lrm_f$ such that $\tau_{\leq n}^*f$ is an equivalence, then $\tau_{\leq n}\Lrm_f\we 0$.
    \begin{proof}
        We sketch the proof. The proof is the same as for the affine case, the details of which can be found
        in~\cite{ha}*{Lemma~8.4.3.17}. In fact, the natural map
        \[
            \tau_{\leq n}(f^*\Lrm_Y)\rightarrow\tau_{\leq n}\Lrm_X
        \]
        is an equivalence. To check this, it is enough to map into an $n$-truncated $\Oscr_X$-module $M$. By the
        universal property of the cotangent complex, we check that the morphism
        \[
        \map_{\CAlg_{R/\Oscr_X}}(\Oscr_X,\Oscr_X\oplus M)\rightarrow\map_{\CAlg_{R/f_*\Oscr_X}}(\Oscr_Y,f_*\Oscr_X\oplus f_*M)
        \]
        is an equivalence, which follows from the fact that
        \[
        \map_{\CAlg_{R/\tau_{\leq n}\Oscr_X}}(\Oscr_X,\tau_{\leq n}\Oscr_X\oplus M)\rightarrow\map_{\CAlg_{R/\tau_{\leq
        n}f_*\Oscr_X}}(\Oscr_Y,\tau_{\leq n}f_*\Oscr_X\oplus f_*M)
        \]
        is an equivalence, since $\tau_{\leq n}\Oscr_Y\rightarrow\tau_{\leq n}f_*\Oscr_X$ is an
        equivalence by hypothesis.
    \end{proof}
\end{lemma}

Let   $R$   be   a   commutative   ring   spectrum.     Then,    the    forgetful    functor
$\CAlg_R\rightarrow\Mod_R$ has a left adjoint
\begin{equation*}
    \Sym_R:\Mod_R\rightarrow\CAlg_R.
\end{equation*}
If $M$ is an $R$-algebra, then $\Sym_R(M)$ is called the symmetric $R$-algebra on $M$.   For
the existence of the functor $\Sym_R$, see~\cite{ha}*{Section~3.1.3}. We can compute the
cotangent complexes of the affine schemes of these symmetric algebras, which provides the
essential step in showing that all maps between connective affine schemes
have cotangent complexes.

\begin{lemma}\label{lem:cotangent}
    Let $M$ be an almost connective $R$-module, and let $S=\Sym_R(M)$. Then the cotangent complex $\Lrm_{\Spec
    S}$ of $\Spec S\rightarrow\Spec R$ exists and is equivalent to the $S$-module $M\otimes_R S$.
    \begin{proof}
        See~\cite{ha}*{Proposition~8.4.3.14}. For any $S$-module $N$, there is a sequence of equivalences
        \begin{equation*}
            \map_{{S}}(M\otimes_R S,N)\we\map_{R}(M,N)\we\map_{R/S}(M,S\oplus N)\we\map_{(\CAlg_R)_{/S}}(S,S\oplus N)
        \end{equation*}
        Thus, $M\otimes_R S$ is an almost connective cotangent complex for  $\Spec  S$, and therefore the unique cotangent complex.
    \end{proof}
\end{lemma}

\begin{proposition}\label{prop:affinecotan}
    If $S\rightarrow T$ is a map of connective commutative $R$-algebras,  then  $\Lrm_{\Spec
    T/\Spec S}$ exists and is connective.
    \begin{proof}
        Indeed, we can write  $T$  as  a  colimit $\colim_i T_i$ of symmetric algebras
        $T_i=\Sym_S {S^{\oplus n_i}}$ so  that  $\Spec  T$  is  a  limit  of  $\Spec  T_i$.
        By Lemma~\ref{lem:colimcotan}, the cotangent complex $\Lrm_T$ is the colimit of  the
        restrictions of $\Lrm_{T_i}$  to  $\Spec  T$.   By  Lemma~\ref{lem:cotangent},  each
        $\Lrm_{T_i}$ is connective. Since colimits of connective $T$-modules are connective,
        $\Lrm_T$ is connective.
    \end{proof}
\end{proposition}

A map of connective commutative $R$-algebras $\phi:\tilde{S}\rightarrow S$ is a nilpotent
thickening if $\pi_0(\phi):\pi_0\tilde{S}\rightarrow\pi_0 S$ is surjective  and  if  the
kernel of $\pi_0(\phi)$ is a nilpotent ideal.
Note that if $S$ is a connective commutative  $R$-algebra,  then  the  maps  $\tau_{\leq
m}S\rightarrow\tau_{\leq n}S$ for $m\geq n$ in the Postnikov tower of $S$ are  nilpotent
thickenings.

\begin{definition}
    A map of sheaves $p:X\rightarrow Y$ is formally smooth if for every nilpotent thickening
    $\tilde{S}\rightarrow S$ the induced map
    \begin{equation}\label{eq:fsmooth}
        X(\tilde{S})\rightarrow X(S)\times_{Y(S)}Y(\tilde{S})
    \end{equation}
    is surjective (that is, surjective on $\pi_0$).
    We say $p:X\rightarrow Y$ is formally \'etale  if  the  maps  in~\eqref{eq:fsmooth}  are
    isomorphisms on $\pi_0$.
\end{definition}

We need the following non-trivial proposition from~\cite{dag7}.

\begin{proposition}[\cite{dag7}*{Proposition~7.26}]\label{prop:smoothprojective}
    If $S$ and $T$ are connective commutative $R$-algebras, then a map $\Spec T\rightarrow\Spec S$ is formally smooth if and only if $\Lrm_{\Spec T/\Spec S}$ is a
    projective $T$-module.
\end{proposition}

To consider the stronger notion of smoothness, we need to consider the notion of compactness
for commutative algebras, and we will need to know later that this notion agrees with the
usual notion of finite presentation for ordinary commutative rings.

\begin{definition}
    A map $S\rightarrow T$ of connective commutative ring spectra is locally of finite
    presentation if $T$ is a compact object of
    $\CAlg_S$~\cite{ha}*{Definition~8.2.5.26}.
\end{definition}

If $T$ is a connective commutative $S$-algebra that is compact in $\CAlg_S$, then it is
compact in $\CAlg_S^{\cn}$. Since truncation preserves compact objects
by~\cite{htt}*{Corollary~5.5.7.4(iii)}, it follows that
$\tau_{\leq 0}T$ is compact in $\tau_{\leq 0}\CAlg_S^{\cn}$.

\begin{lemma}\label{lem:lfp}
    Let $R$ be a discrete commutative ring, and let $S$ be a discrete commutative $R$-algebra. Then,
    $S$ is compact as a discrete commutative $R$-algebra if and only if $S$ is finitely presented.
    \begin{proof}
        Suppose that $S$ is finitely presented, so it can be written as a quotient
        $R[X]/I$, where $X$ is a finite set, and $I$ is a finitely generated ideal. We can
        write $S$ as the pushout
        \begin{equation*}
            \begin{CD}
                \Sym_R I    @>>>    R[X]\\
                @VVV                @VVV\\
                R           @>>>    S
            \end{CD}
        \end{equation*}
        of $R$-algebras. But, this means that if
        $R[Y]\rightarrow\Sym_R I$ exhibits $\Sym_R I$ as a finitely generated $R$-algebra,
        the following square is also a pushout square:
        \begin{equation*}
            \begin{CD}
                R[Y]    @>>>    R[X]\\
                @VVV                @VVV\\
                R           @>>>    S.
            \end{CD}
        \end{equation*}
        Since $R[Y]$, $R[X]$, and $R$ are compact, it follows that $S$ is
        compact as well.

        Now, suppose that $S$ is a compact (discrete) commutative $R$-algebra. Then, $S$ is a
        retract of a finitely presented commutative $R$-algebra $R[X]/I$. Indeed, we can
        write $S$ as a filtered colimit of finitely presented commutative $R$-algebras; by
        compactness, the identity map on $S$ factors through a finite stage. It
        suffices to show that the kernel of $R[X]/I\rightarrow S$ is finitely generated. We
        proceed by noetherian induction. Let $\phi$ be the composition
        \begin{equation*}
            R[X]/I\rightarrow S\rightarrow R[X]/I.
        \end{equation*}
        We may write $I=(p_1(X),\ldots,p_k(X))$, an ideal generated by $k$ polynomials in
        $X$. Let $R_0$ be the subring of $R$ generated over $\ZZ$ by the
        coefficients appearing in the $p_i$ and in the polynomials $\phi(x_i)$ for $x_i\in X$. This is a
        finitely generated commutative $\ZZ$-algebra, so it is in particular noetherian. By
        our choice of $R_0$, we can define an ideal $I_0$ of $R_0[X]$ generated by the same
        polynomials. Moreover, $\phi$ defines a morphism $\phi_0:R_0[X]/I_0\rightarrow
        R_0[X]/I_0$. Let $S_0$ be the image of $\phi_0$, which is a subring of $R_0[X]/I_0$.
        There is an exact sequence of $R_0$-modules
        \begin{equation*}
            0\rightarrow J_0\rightarrow R_0[X]/I_0\rightarrow S_0\rightarrow 0.
        \end{equation*}
        Since $S_0$ is finitely generated and $R_0[X]/I_0$ is noetherian, it follows that
        $S_0$ is finitely presented, and that $J_0$ is finitely generated.
        Tensoring with $R$ over $R_0$, we obtain an exact sequence
        \begin{equation*}
            J_0\otimes_{R_0}R\rightarrow R[X]/I\rightarrow S\rightarrow 0.
        \end{equation*}
        The kernel $J_0$ is finitely generated, and it surjects onto the kernel $J$ of
        $R[X]/I\rightarrow S$. It follows that $J$ is finitely generated, and hence that $S$
        is finitely presented.
    \end{proof}
\end{lemma}

The previous lemma implies that the next definition agrees with the usual definition of
smooth maps between ordinary affine schemes.

\begin{definition}\label{def:smooth}
    A map $\Spec T\rightarrow\Spec S$ is smooth if it is formally smooth  and  $S\rightarrow
    T$ is locally of finite presentation.
\end{definition}

At first glance, the condition that $\Spec T\rightarrow\Spec S$ is surjective in the next
lemma might seem strange. But, we will show in Theorem~\ref{thm:lifts} that this is satisfied if
$\Spec T\rightarrow\Spec S$ is smooth, and if the map $\Spec\pi_0T\rightarrow\Spec\pi_0S$ is
a surjective map of ordinary schemes.

\begin{lemma}\label{lem:2o3lfp}
    Let $R\rightarrow S\rightarrow T$ be maps of connective commutative algebras. If $T$ is
    locally of finite presentation over $R$ and over $S$, and if $\Spec T\rightarrow\Spec S$
    is a surjective map in $\Shv_R^{\et}$, then $S$ is locally of finite presentation over $R$.
    \begin{proof}
        First, the maps $\pi_0R\rightarrow\pi_0S\rightarrow\pi_0T$ satisfy the same
        hypotheses by Lemma~\ref{lem:lfp}. Thus, by~\cite{ega4}*{Proposition I.1.4.3(v)}, $\pi_0S$ is a
        $\pi_0R$-algebra that is locally of finite presentation. Now,
        by~\cite{ha}*{Theorem~8.4.3.18}, it is enough to show that $\Lrm_S=\Lrm_{R/S}$ is a perfect
        $S$-module. There is a fiber sequence
        \begin{equation*}
            \Lrm_S\otimes_S T\rightarrow\Lrm_T\rightarrow\Lrm_{S/T}
        \end{equation*}
        of cotangent complexes. Again, by~\cite{ha}*{Theorem~8.4.3.18}, $\Lrm_T$ and $\Lrm_{S/T}$ are perfect since
        $R\rightarrow T$ and $S\rightarrow T$ are locally of finite presentation. It follows
        that $\Lrm_S\otimes_S T$ is perfect. Since $\Spec T\rightarrow\Spec S$ is
        surjective in $\Shv_R^{\et}$, there are \'etale local sections. Thus, there is a faithfully flat \'etale $S$-algebra
        $P$ and maps $S\rightarrow T\rightarrow P$. Since $\Lrm_S\otimes_S T$ is perfect,
        the $U$-module $\Lrm_S\otimes_S U$ is perfect. But, by faithfully flat descent,
        it follows that $\Lrm_S$ is perfect (to see this, one can either refer forward to
        Lemma~\ref{lem:stm}, or use the fact that an $S$-module is perfect if and only if it is
        dualizable and the fact that dualizability data can be constructed \'etale locally
        by Example~\ref{ex:moddescent}).
    \end{proof}
\end{lemma}

The following example will be used later in the paper.

\begin{example}\label{ex:gln}
    The sheaf of $R$-module endomorphisms of $R^{\oplus n}$ is representable by an affine
    monoid scheme $\Mrm_n$, where
    \begin{equation*}
        \Mrm_n=\Spec\Sym_R\End_R(R^{\oplus n}).
    \end{equation*}
    Given a commutative $R$-algebra $S$, an element of $\Mrm_n(S)$ is a commutative $R$-algebra map
    \begin{equation*}
        \Sym_R\End_R(R^{\oplus n})\rightarrow S.
    \end{equation*}
    But these are equivalent to the $R$-module maps
    \begin{equation*}
        R^{\oplus n^2}\we\End_R(R^{\oplus n})\rightarrow S,
    \end{equation*}
    which by adjunction is the $S$-module
    \begin{equation*}
        S^{\oplus n^2}\we\End_S(S^{\oplus n}).
    \end{equation*}
    Since the space of $S$-module endomorphisms of $S^{\oplus n}$ has a natural monoid structure,
    this shows that $\Mrm_n$ is a monoid scheme.
    We can invert the determinant element of $\pi_0\Mrm_n$, so the sheaf of $S$-module automorphisms of $S^{\oplus n}$ is representable by an affine group scheme $\GL_n$.
    Because the cotangent complex of $\Mrm_n$ at an $S$-point is $S^{\oplus n^2}$, which  is
    a projective $S$-module, the affine schemes $\Mrm_n$  and  $\GL_n$  are  smooth  over  $R$.
\end{example}

\subsection{Geometric sheaves}

Let $R$ be a connective commutative ring spectrum.
The goal of this section is to study certain geometric classes of sheaves in $\Shv_R^{\et}$
built inductively from the representable sheaves by forming smooth quotients.
The notions of $n$-stack here have been studied extensively by Simpson~\cite{simpson},
To\"en and Vezzosi~\cite{hag2}, and Lurie~\cite{dag0}, and we base our approach on theirs.

We define $n$-geometric morphisms and smooth $n$-geometric morphisms inductively as follows.
    \begin{itemize}
        \item A morphism $f:X\rightarrow Y$ in $\Shv_R^{\et}$ is $0$-geometric  if  for  any
            $\Spec S\rightarrow Y$, the fiber product $X\times_Y \Spec S$  is  equivalent
            to   $\coprod_i\Spec   T_i$   for   some   connective    commutative
            $R$-algebras $T_i$.
        \item   A $0$-geometric morphism $f$ is smooth if $X\times_Y \Spec S\rightarrow\Spec S$ is smooth for all $\Spec S\rightarrow Y$.
            To be clear, if $X\times_Y\Spec S\we\coprod_i\Spec T_i$, then this means that
            each morphism $\Spec T_i\rightarrow\Spec S$ is smooth in the sense of
            Definition~\ref{def:smooth}.
        \item A morphism $X\rightarrow Y$ in $\Shv_R^{\et}$  is  $n$-geometric  if  for  any
            $\Spec S\rightarrow Y$, there is a smooth surjective $(n-1)$-geometric  morphism
            $U\rightarrow X\times_Y \Spec S$
            where $U$ is a disjoint union of affines.
        \item An $n$-geometric morphism $f:X\rightarrow Y$ is smooth  if  for  every  $\Spec
            S\rightarrow Y$ we may take the above map $U\rightarrow X\times_Y \Spec S$  such
            that  the  composition  $U\rightarrow  X\times_Y\Spec  S\rightarrow\Spec  S$  is
            a smooth $0$-geometric morphism.
    \end{itemize}
We say that an $n$-geometric morphism $X\rightarrow Y$ is an $n$-submersion  if  it  is  smooth  and
surjective.  If, moreover, $X$ is a disjoint sum of representables, then  we  call  such a
morphism an $n$-atlas. A $1$-geometric sheaf with a Zariski atlas is a derived scheme.
What this means is that a $1$-geometric sheaf $X$ has an atlas $\coprod_i\Spec
T_i\rightarrow X$ which is a $0$-geometric morphism, and, for every point $\Spec
S\rightarrow X$, the pullback $\Spec T_i\times_X\Spec S\rightarrow\Spec S$ is Zariski open.
Similarly, a $1$-geometric sheaf with an \'etale atlas is a Deligne-Mumford stack.

Any $0$-geometric sheaf $X$ is a disjoint union of sheaves $\coprod_{i\in I}\Spec S_i$,
where the $S_i$ are connective commutative $R$-algebras. If $I$ is \emph{finite}, then we call
the sheaf representable. In this case $X=\coprod_{i=1}^{n}\Spec
S_i\we\Spec(S_1\times\cdots\times S_n)$.

A $0$-geometric sheaf is quasi-compact if it is representable, and
a $0$-geometric morphism $f:X\to Y$ is quasi-compact if,
for all $\Spec S\to Y$, the pullback $X\times_Y\Spec S$ is representable.
Inductively, an $n$-geometric sheaf $X$ if quasi-compact if there exists an
$(n-1)$-geometric quasi-compact submersion of the form $\Spec S\to X$, and an $n$-geometric
morphism $f:X\to Y$ is quasi-compact if for each map $\Spec S\to Y$, the fiber
$X\times_Y\Spec S$ is a quasi-compact $n$-geometric sheaf.
Finally, an $n$-geometric morphism $f:X\to Y$ is quasi-separated if the diagonal $X\to
X\times_Y X$ is quasi-compact.

\begin{definition}\label{def:ngeometricprops}
    An $n$-geometric sheaf $X$ is locally of finite presentation over $\Spec R$ if it has an $(n-1)$-atlas
    \begin{equation*}
        \coprod_i\Spec S_i\rightarrow X
    \end{equation*}
    such that each $S_i$ is a connective commutative $R$-algebra that is locally of finite presentation. An
    $n$-geometric morphism $X\rightarrow Y$ is locally of finite presentation if  for  every
    $S$-point of $Y$, $X\times_Y\Spec  S$  is  locally  of  finite  presentation  over
    $\Spec S$. By definition, a smooth $n$-geometric morphism is locally of finite
    presentation.
\end{definition}

We need not only $n$-geometric sheaves, but it is also important to have a good theory of
sheaves that are only locally geometric.
A sheaf $X\rightarrow \Spec S$ is locally geometric if it can be written as a filtered
colimit
\begin{equation*}
    X\we\colim_i X_i,
\end{equation*}
where each sheaf $X_i$ is $n_i$-geometric for some $n_i$ and where the maps
$X_i\rightarrow X$ are monomorphisms. If we can furthermore
take the $X_i$ to be locally of finite presentation, we say that $X$ is locally
geometric and locally of finite presentation.
A morphism $f:X\rightarrow Y$ is locally geometric (locally of finite presentation) if for every $\Spec S\rightarrow Y$,
the pullback $X\times_Y\Spec S$ is locally geometric (locally of finite presentation) over $\Spec S$.

We say that a locally geometric morphism $X\rightarrow Y$  which is locally of finite presentation
is smooth if for every $\Spec S\rightarrow Y$, the pullback $X\times_Y\Spec
S\rightarrow\Spec S$ has a cotangent complex of Tor-amplitude contained in $[-n,0]$ for
some non-negative integer $n$ (depending on $S$).

Note the following easy but important facts.

\begin{lemma}\label{lem:fibersngeometric}
    If $f:X\rightarrow Z$ is an $n$-geometric morphism (resp. smooth $n$-geometric morphism), and if $Y\rightarrow Z$ is any
    morphism, then the pullback $f_{Y}:X\times_Z Y\rightarrow Y$  is  $n$-geometric  (resp.
    $n$-geometric and smooth).
\end{lemma}

\begin{lemma}\label{lem:fiberwise}
    A morphism $f:X\rightarrow Y$ is $n$-geometric if and  only  if  for  every  map  $\Spec
    S\rightarrow Y$, the morphism $X\times_Y\Spec  S\rightarrow\Spec  S$  is  $n$-geometric.
\end{lemma}

The following lemma can be found in~\cite{dag0}.  We include a proof for the reader's convenience.

\begin{lemma}\label{lem:lurie}
    Suppose that $X\xrightarrow{f}Y\xrightarrow{g}Z$ are composable  morphisms  of  sheaves.
    \begin{enumerate}
        \item   If $f$ and $g$ are $n$-geometric (resp. smooth and
            $n$-geometric),  then  $g\circ  f$   is   $n$-geometric   (resp.    smooth   and
            $n$-geometric).
        \item If $f$ is an $n$-submersion and $g\circ  f$  is  $(n+1)$-geometric,  then  $g$
            is $(n+1)$-geometric.
        \item  If  $g\circ  f$  is  $n$-geometric  and  $g$   is   $(n+1)$-geometric,   then
            $f$ is $n$-geometric.
    \end{enumerate}
            \begin{proof}
                We prove (1) by  induction  on  $n$.   Using  Lemma~\ref{lem:fiberwise},  it
                suffices to suppose that $Z$ is representable. Assume that $n=0$. Then, the fact that $g$ is
                $0$-geometric implies that $Y$ is representable, and the fact  that  $f$  is
                $0$-geometric then implies that $Z$ is representable. Evidently any morphism of
                representables is $0$-representable, and compositions of smooth morphisms of
                representables are smooth.
                Now assume the statement (1) for $(n-1)$-geometric morphisms. Since we assume
                that $Z$  is  representable,  it  suffices  to  find  an  $(n-1)$-submersion
                $U\rightarrow X$ where $U$ is a disjoint  union  of  representables.   Since
                $g$ is $n$-representable, there is an  $(n-1)$-submersion  $V\rightarrow  Y$
                where  $V$  is  a  sum  of  representables.   Constructing   the   pullback
                $X\times_Y V$, we know  by  the  $n$-geometricity  of  $f$  and  the  formal
                representability  of  $V$  that   there   is   an   $(n-1)$-submersion
                $U\rightarrow   X\times_Y   V$   with   $U$   a   sum   of   representables.
                This is summarized in the diagram
                \begin{center}
                    \begin{tikzpicture}[description/.style={fill=white,inner sep=2pt}]
                            \matrix   (m)   [matrix   of   math    nodes,    row    sep=3em,
                            column  sep=2.5em,   text   height=1.5ex,   text   depth=0.25ex]
                            {     U & X\times_Y V & V &  \\
                                      & X & Y   &   Z\\};
                            \path[->,font=\scriptsize]
                            (m-1-1)  edge  node[auto]  {$  (n-1)\textrm{-sub.}  $}   (m-1-2)
                            (m-1-2) edge node[above] {$ $} (m-1-3)
                            (m-1-2)  edge  node[left]  {$  (n-1)\textrm{-sub.}  $}   (m-2-2)
                            (m-1-3)  edge  node[right]  {$  (n-1)\textrm{-sub.}  $}  (m-2-3)
                            (m-2-3) edge node[below] {$ g $} (m-2-4)
                            (m-2-2) edge node[below] {$ f $} (m-2-3);
                    \end{tikzpicture}.
                \end{center}
                Since, by the induction hypothesis, composition of  $(n-1)$-submersions  are
                $(n-1)$-submersions, the inductive step follows.

                To prove (2), it is again enough to assume that $Z$ is representable. Then,
                there  is  an   $n$-atlas   $u:U\rightarrow   X$   since   $g\circ   f$   is
                $(n+1)$-geometric.  Since $f$  is  an  $n$-submersion,  the  composition
                $f\circ u$ is an $n$-atlas by part (1).  Hence,  $g$  is  $(n+1)$-geometric.

                To prove (3), suppose that $p:\Spec S\rightarrow Y$ is a point of $Y$. Consider the
                diagram
                \begin{center}
                    \begin{tikzpicture}[description/.style={fill=white,inner sep=2pt}]
                            \matrix   (m)   [matrix   of   math    nodes,    row    sep=3em,
                            column  sep=2.5em,   text   height=1.5ex,   text   depth=0.25ex]
                            {    &    X\times_Y    \Spec    S    &    \Spec    S    &    V\\
                                  U & X\times_Z\Spec  S  &  Y\times_Z\Spec  S  &  \Spec  S\\
                                      & X & Y   &   Z\\};
                            \path[->,font=\scriptsize]
                            (m-1-2) edge node[auto] {$ $} (m-1-3)
                            (m-1-2) edge node[left] {$ $} (m-2-2)
                            (m-1-3) edge node[right] {$ $} (m-2-3)
                            (m-2-3) edge node[below] {$  $} (m-2-4)
                            (m-2-2) edge node[below] {$  $} (m-2-3)
                            (m-2-2) edge node[] {$ $} (m-3-2)
                            (m-2-3) edge node[] {$ $} (m-3-3)
                            (m-2-1) edge  node[below]  {$  (n-1)\mathrm{-atlas}  $}  (m-2-2)
                            (m-1-4)  edge  node[auto]   {$   n\mathrm{-atlas}   $}   (m-2-3)
                            (m-2-1) edge[dashed] node[] {$ $} (m-1-4)
                            (m-2-4)   edge   node[right]   {$   g\circ    p    $}    (m-3-4)
                            (m-2-4) edge node[auto] {$ p $} (m-3-3)
                            (m-3-2) edge node[below] {$ f $} (m-3-3)
                            (m-3-3) edge node[below] {$ g $} (m-3-4);
                    \end{tikzpicture},
                \end{center}
                where the squares are all pullback squares, $U\rightarrow X\times_Z\Spec S$
                is an $(n-1)$-atlas (or the identity map if $n=0$) and $V\rightarrow Y\times_Z\Spec S$ is an $n$-atlas.
                Since $V\rightarrow Y\times_Z\Spec S$ is surjective, up to refining $U$,  we
                may assume  that  the  composite  $U\rightarrow  Y\times_Z\Spec  S$  factors
                through $V$.  The map $U\rightarrow V$ is thus $0$-geometric.  By part  (1),
                the map $U\rightarrow Y\times_Z\Spec S$ is $n$-geometric.  By part  (2),  it
                follows   that   $X\times_Z\Spec   S\rightarrow   Y\times_Z\Spec    S$    is
                $n$-geometric.    Therefore,   $X\times_Y\Spec   S\rightarrow\Spec   S$   is
                $n$-geometric, and we conclude  by  Lemma~\ref{lem:fiberwise}  that  $f$  is
                $n$-geometric.
            \end{proof}
\end{lemma}

\begin{remark}\label{rem:quasicompact}
    The previous lemma goes through as stated with the additional assumptions and
    conclusions of quasi-compactness.
\end{remark}

\begin{lemma}\label{lem:atlaseslfp}
    Suppose that $X$ is an $n$-geometric sheaf that is locally of finite presentation.
    If $U=\coprod_i\Spec T_i\xrightarrow{p} X$ is \emph{any} atlas, then each $T_i$ is locally
    of finitely presentation over $R$.
    \begin{proof}
        Let $V=\coprod_i\Spec S_i\xrightarrow{q} X$ be an atlas where each $S_i$ is locally of
        finite presentation over $R$. Since $V\rightarrow X$ is a surjection of sheaves, we
        may assume, possibly by refining $U$, that there is a factorization of $p$ as $U\rightarrow
        V\xrightarrow{q} X$. Now, consider the fiber-product $U\times_X V$, which is a
        smooth $(n-1)$-geometric sheaf over either $U$ or $V$. Let $W=\coprod_i\Spec
        P_i\rightarrow U\times_X V$ be an atlas. Since $U\times_X V\rightarrow U$ is
        surjective, we may arrange indices so that the composition $W\rightarrow U$ is a
        coproduct of smooth surjections of the form $\Spec P_i\rightarrow\Spec T_i$. Assume
        also that in the map $U\rightarrow V$ we have $\Spec T_i\rightarrow\Spec S_i$. Then,
        there is a composition of commutative ring maps $S_i\rightarrow T_i\rightarrow P_i$.
        The composite is locally of finite presentation since it is smooth, the map
        $T_i\rightarrow P_i$ is locally of finite presentation for the same reason, and by
        construction the map $\Spec P_i\rightarrow\Spec T_i$ is surjective. Thus, the
        conditions of Lemma~\ref{lem:2o3lfp} are satisfied. It follows that $T_i$ is locally
        of finite presentation over $S_i$. Since $S_i$ is locally of finite presentation
        over $R$, it follows that $T_i$ is locally of finite presentation over $R$.
    \end{proof}
\end{lemma}

Now we prove an analogue of~\cite{dag0}*{Principle~5.3.5}.

\begin{lemma}\label{lem:screwing}
    Suppose that $\Prm$ is a property of sheaves. Suppose that every disjoint union of
    affines has property $\Prm$, and suppose that whenever $U\rightarrow X$ is a surjective
    morphism of sheaves such that $U_X^k=U\times_X\cdots\times_X U$ has property $\Prm$ for
    all $k\geq 0$, then $X$ has property $\Prm$. Then, all $n$-geometric sheaves have property $\Prm$.
    \begin{proof}
        Let $X$ be an $n$-geometric sheaf. Then, there exists a smooth $(n-1)$-geometric
        surjection $U\rightarrow X$, where $U$ is a disjoint union of affines. Each product
        $U_X^k$ is $(n-1)$-geometric. So, it suffices to observe that the statement follows
        by induction.
    \end{proof}
\end{lemma}

\begin{lemma}\label{lem:stacksdescent}
    Let $X\rightarrow Y$ be a surjection of sheaves.
    Suppose  that  $X$  and  $X\times_Y  X$  are   $n$-geometric   stacks   and   that   the
    projections  $X\times_Y  X\rightarrow  X$  are   $n$-geometric and smooth.    Then,   $Y$   is   an
    $(n+1)$-geometric stack. If, in addition, $X$ is quasi-compact and $X\rightarrow Y$ is a
    quasi-compact morphism, then $Y$ is quasi-compact. Finally, if $X$ is locally of finite presentation, then so is $Y$.
    \begin{proof}
        Let
        \begin{equation*}
            \coprod_i\Spec S_i\rightarrow X
        \end{equation*}
        be an atlas for $X$. For any $i$ and $j$, we have a diagram of pullbacks
        \begin{equation*}
            \begin{CD}
                \Spec S_i\times_Y\Spec S_j  @>>>    \Spec S_i\times_Y X   @>>>    \Spec S_i\\
                @VVV                                @VVV                    @VVV\\
                X\times_Y \Spec S_j         @>>>    X\times_Y X         @>>>    X\\
                @VVV    @VVV    @VVV\\
                \Spec S_j   @>>>    X   @>>>    Y.
            \end{CD}
        \end{equation*}
        We will show that the composite $\coprod_i\Spec S_i\rightarrow X\rightarrow Y$
        is an $n$-submersion. The surjectivity follows by hypothesis. To check that $\Spec
        S_i\rightarrow X$ is smooth and $n$-geometric, it is enough to check on the fiber of
        $\Spec S_j$ for all $j$. Since $X\times_Y X\rightarrow X$ and $\Spec S_i\rightarrow
        X$ are smooth and $n$-geometric, it follows that $\Spec S_i\times_Y\Spec
        S_j\rightarrow X\times_Y\Spec S_j$ and $X\times_Y\Spec S_j\rightarrow S_j$ are
        smooth and $n$-geometric. Therefore, the composite is as well, which completes the
        proof of the first statement. To prove the second statement, note that we can take
        $\coprod_i\Spec S_i$ to be a finite disjoint union, since $X$ is quasi-compact.
        Then, since $\coprod_i\Spec S\rightarrow X$ and $X\rightarrow Y$ are quasi-compact,
        it follows that the composition is quasi-compact by Remark~\ref{rem:quasicompact}.
        The third statement
        is immediate, since $\coprod_i\Spec S_i$ can be chosen so that each $S_i$ is locally
        of finite presentation.
    \end{proof}
\end{lemma}

\begin{lemma}\label{lem:the}
    Suppose  that  $f:X\rightarrow  Y$  is  a  morphism  of  sheaves   where   $X$   is   an
    $n$-geometric sheaf and the diagonal $Y\rightarrow  Y\times_{\Spec  R}Y$  is  $n$-geometric.   Then,  $f$  is
    $n$-geometric. Moreover, if $f$ is smooth and surjective, then $Y$ is an $(n+1)$-geometric stack.
    \begin{proof}
        Since $X$ is $n$-geometric, there is an $(n-1)$-submersion
        $\coprod_{i}\Spec T_i\rightarrow X$.  Suppose that  $\Spec  S\rightarrow  Y$  is
        arbitrary.   Form  the  fiber  products  $X\times_Y  \Spec  S$  and  $\coprod_i\Spec
        T_i\times_Y\Spec S$, and note that the map
        \begin{equation*}
            \coprod_i\Spec T_i\times_Y \Spec S\rightarrow X\times_Y\Spec S
        \end{equation*}
        is an $(n-1)$-submersion. The map
        \begin{equation*}
            \coprod_i\Spec  T_i\times_Y\Spec   S\rightarrow\coprod_i\Spec   T_i\otimes_R   S
        \end{equation*}
        is $n$-geometric because it is the  pullback  of  $\coprod_i\Spec  T_i\times_{\Spec
        R}\Spec S\rightarrow Y\times_{\Spec R}  Y$  along  the  diagonal  map  $Y\rightarrow
        Y\times{\Spec  R}$.   Therefore  $\coprod_i\Spec  T_i\times_Y\Spec  S$   admits   an
        $(n-1)$-submersion from a disjoint union of affines $\coprod_{i,j}\Spec U_{ij}$. We
        obtain the following diagram
        \begin{center}
            \begin{tikzpicture}[description/.style={fill=white,inner sep=2pt}]
                    \matrix (m) [matrix of math nodes, row sep=3em,
                    column sep=2.5em, text height=1.5ex, text depth=0.25ex]
                    {   \coprod\Spec   U_{ij}   &   \coprod    \Spec    T_i\times_Y\Spec{S}&
                    X\times_Y\Spec S & \Spec S \\
                          \coprod\Spec T_i\otimes_R  S  &  \coprod\Spec  T_i  &  X  &  Y\\};
                    \path[->,font=\scriptsize]
                    (m-1-1)  edge  node[auto]  {$   (n-1)\textrm{-submersion}   $}   (m-1-2)
                    (m-1-2)  edge  node[above]  {$  (n-1)\textrm{-submersion}   $}   (m-1-3)
                    (m-1-3) edge node[below left] {$  $} (m-1-4)
                    (m-1-2)   edge   node[left]   {$   n\textrm{-geometric}    $}    (m-2-1)
                    (m-1-2) edge node[] {$ $} (m-2-2)
                    (m-1-3) edge node[] {$ $} (m-2-3)
                    (m-1-4) edge node[] {$ $} (m-2-4)
                    (m-2-2)  edge  node[below]  {$  (n-1)\textrm{-submersion}   $}   (m-2-3)
                    (m-2-3) edge node[above] {$ f $} (m-2-4);
            \end{tikzpicture}.
        \end{center}
        By Lemma~\ref{lem:lurie}, the composition of the top two horizontal maps  is  also
        an  $(n-1)$-submersion  from  a  disjoint  union  of  affines,   establishing   that
        $f$ is $n$-geometric. The second claim is clear.
    \end{proof}
\end{lemma}

\begin{lemma}\label{lem:loops}
    If $X$ is $n$-geometric, and if $p:\Spec  S\rightarrow  X$  is  a  point  of  $X$,  then
    $\Omega_p X=\Spec S\times_X \Spec S\rightarrow\Spec S$ is an $(n-1)$-geometric morphism.
    The projection $X^{S^{m}}\rightarrow X$ induced by choosing a point in the $m$-sphere is
    an $(n-m)$-geometric map.
    \begin{proof}
        We use the equivalent description of $\Omega_p X$ as the pullback  in  the  diagram
        \begin{equation*}
            \begin{CD}
                \Omega_p X  @>>>    \Spec S\times_{\Spec R}\Spec S\\
                @VVV                @V(p,p)VV\\
                X           @>>>    X\times_{\Spec R}X.
            \end{CD}
        \end{equation*}
        Since the diagonal of $X$  is  $(n-1)$-geometric,  it  follows  that  the  composite
        \begin{equation*}
            \Omega_p   X\rightarrow\Spec   S\times_{\Spec   R}\Spec   S\rightarrow\Spec    S
        \end{equation*}
        is also $(n-1)$-geometric. To prove the second statement, it suffices to note that the fiber of the  projection
        map over a point $p:\Spec S\rightarrow  X$  is  the  $m$-fold  iterated  loop  space
        $\Omega_p^m  X$.   So,  this  follows  from   the   first   part   of   the   lemma.
    \end{proof}
\end{lemma}

\begin{example}\label{ex:ngroups}
    If $G$ is a smooth $n$-geometric stack of groups (\emph{i.e.}, grouplike
    $A_{\infty}$-spaces), then $\Brm G$ is a pointed
    smooth $(n+1)$-geometric stack. Indeed, the loop space of $\Brm G$ at the canonical point is
    just $G$.  Therefore, the point $\Spec R\rightarrow\Brm G$  is  an  $n$-submersion.   Using
    Lemma~\ref{lem:lurie} part (2), the claim follows.
\end{example}

For the next two lemmas, fix a base sheaf $Z$ in $\Shv_R^{\et}$, and consider the
$\infty$-topos $\Shv^{\et}_{/Z}$ of objects over $Z$.
Let $\Shv^n_{/Z}$ be the full subcategory of $\Shv^{\et}_{/Z}$ consisting of the $n$-geometric
morphisms $Y\rightarrow Z$.

\begin{lemma}\label{lem:limitsnstacks}
    The full subcategory $\Shv^n_{/Z}$ of $\Shv^{\et}_{/Z}$ is closed under finite limits in $\Shv^{\et}_{/Z}$.
    \begin{proof}
        As $\Shv^n_{/Z}$ has a terminal object which agrees  with  the  terminal  object  of
        $\Shv^{\et}_{/Z}$, it is enough to check the case of pullbacks. Suppose that $X\rightarrow Y$
        and $W\rightarrow Y$ are two morphisms in $\Shv^n_{/Z}$.  In  order  to  check  that
        $X\times_Y W$ is in $\Shv^n_{/Z}$ it suffices to check  that  $X\times_Y  W$  is  in
        $\Shv^n_{/Y}$ since $Y\rightarrow Z$ is $n$-geometric.  Moreover, we  can  obviously
        reduce to the case that $Y\we\Spec S$ is representable. Then, $X$ and $W$ are $n$-geometric
        stacks over $S$. Taking an $n$-atlas $U\rightarrow X$ and an $n$-atlas $V\rightarrow
        W$, the fiber product $U\times_Y V$ is an $n$-atlas for $X\times_Y W$ by  using  the
        stability of geometricity and smoothness under pullbacks.
    \end{proof}
\end{lemma}

\begin{lemma}\label{lem:limits0}
    The full subcategory of $\Shv^{\et}_{/Z}$ consisting of quasi-compact $0$-geometric sheaves over $Z$ is closed under
    \emph{all} limits in $\Shv^{\et}_{/Z}$.
    \begin{proof}
        It suffices to note that arbitrary limits of representables are representable, since
        the   \icat   of   connective   commutative   $R$-algebras   has    all    colimits.
    \end{proof}
\end{lemma}

\begin{lemma}\label{lem:nlimits}
    A finite limit of $n$-geometric morphisms (locally of finite presentation) is $n$-geometric
    (and locally of finite presentation).
    \begin{proof}
        The proof is by induction on $n$. The base case $n=0$ simply follows because finite
        limits of representable sheaves are representable, and finite limits distribute over coproducts.
        Suppose the lemma is true for $k$-geometric
        sheaves for all $k<n$, and let $f_i:X_i\rightarrow Y_i$ be a finite diagram of
        $n$-geometric morphisms (locally of finite presentation). Let $f:X\rightarrow Y$ be
        the limit. Let $\Spec S\rightarrow Y$ be an $S$-point. Then, we may construct an
        atlas for the pullback $X\times_Y\Spec S$ as the (finite) limit of a compatible
        family of atlases for the pullbacks $X_i\times_{Y_i}\Spec S$. The morphism from this
        atlas to $X\times_Y\Spec S$ is $(n-1)$-geometric by the inductive hypothesis. It is
        also clear that it is a submersion. If the maps are locally of finite presentation,
        then the atlases over each $X_i\times_{Y_i}\Spec S$ may be chosen to be locally of
        finite presentation, and hence their (finite) limit is again locally of finite
        presentation.
    \end{proof}
\end{lemma}

\begin{lemma}\label{lem:cosimplicial}
    Let $X^\bullet$ be a cosimplicial diagram of quasi-separated $n$-geometric sheaves over
    $Z$. Then, the limit $X=\lim_{\Delta}X^\bullet$ is $n$-geometric over $Z$.
    \begin{proof}
        By \cite{goerss-jardine}*{Proposition VII.1.7}, there  is  push-out  diagram  for  any  $m$
        \begin{equation*}
            \begin{CD}
                \Delta^\bullet_m\times\partial\Delta^m    @>>>    sk_{m-1}\Delta^\bullet\\
                @VVV                                        @VVV\\
                \Delta^\bullet_m\times\Delta^m      @>>>    sk_m\Delta^\bullet
            \end{CD}
        \end{equation*}
        of cosimplicial spaces. Given $X^{\bullet}$ we obtain a pullback diagram of sheaves
        \begin{equation*}
            \begin{CD}
                \map(sk_m\Delta^\bullet,X^\bullet)  @>>>
                \map(\Delta^\bullet_m\times\Delta^m,X^\bullet)\we X^m\\
                @VVV    @VVV\\
                \map(sk_{m-1}\Delta^\bullet,X^\bullet) @>>>
                \map(\Delta^{\bullet}_m\times\partial\Delta^m,X^\bullet)\we (X^m)^{S^{m-1}}.
            \end{CD}
        \end{equation*}
        By   Lemma~\ref{lem:loops},   the   map    $(X^m)^{S^{m-1}}\rightarrow    X^m$    is
        $(n-m-1)$-geometric.  Since $X^m$ is $n$-geometric, Lemma~\ref{lem:lurie}  part  (3)
        implies that the left-hand  vertical  maps  above  are  $(n-m-2)$-geometric.   Thus,
        if $m\geq n-2$, we see that the  left-hand  vertical  maps  above  are  $0$-geometric.
        Moreover, by hypothesis, each diagonal $X^m\to X^m\times_Z X^m$ is quasi-compact, so, pulling back, we see that each of the maps
        \[
        (X^m)^{S^{m-1}}\to (X^m)^{S^{m-2}}\to\cdots\to (X^m)^{S^1}\to X^m
        \]
        is quasi-compact, so the composite $(X^m)^{S^{m-1}}\to X^m$ and the section $X^m\to
        (X^m)^{S^{m-1}}$ are as well by Remark~\ref{rem:quasicompact}.
        We conclude that the left-hand vertical map is $0$-geometric and quasi-compact.
        As we have equivalences
        \begin{equation*}
            \lim_{\Delta}X^{\bullet}\we
            \lim_m\map(sk_m\Delta^{\bullet},X^\bullet)\we\lim_{m\geq
            n-2}\map(sk_m\Delta^{\bullet},X^{\bullet}),
        \end{equation*}
        we see that $\lim_{\Delta}X^{\bullet}$ is a limit of quasi-compact
        $0$-geometric morphisms over $\map(sk_{n-2}\Delta^\bullet,X^{\bullet})$,  which,  as
        as a finite limit of $n$-geometric sheaves  over  $Z$,  is  $n$-geometric  over  $Z$.
        Hence, by Lemma~\ref{lem:limits0}, the limit is $n$-geometric.
    \end{proof}
\end{lemma}

\begin{proposition}\label{prop:retractsnstacks}
    If $Y$ is a retract of a sheaf $X$ over $Z$, and if $X$ is quasi-separated and
    $n$-geometric over $Z$, then $Y$ is $n$-geometric over $Z$.
    \begin{proof}
        We refer to \cite{htt}*{Section~4.4.5} for details about retracts in \icats. In particular, any
        retract   in   $\Shv^{\et}_{/Z}$   is   given   as   the   limit   of    a    diagram
        $\tilde{X}:\Idem\rightarrow\Shv^{\et}_{/Z}$, where $\Idem$ is an \icat with only  one
        object $*$ and with finitely many simplices in each
        degree.  Let $X=\tilde{X}(*)$.  It follows that the  cosimplicial  replacement  (see
        \cite{bousfield-kan}*{XI 5.1}) of $p$ is a cosimplicial
        sheaf  $X^{\bullet}$  which  in  degree  $k$  is  a  finite  product  of  copies  of
        $X$. Thus, if $p$ takes values in quasi-separated $n$-geometric sheaves over $Z$,
        then each $X^{k}$ is still quasi-separated and $n$-geometric.
        By Lemma~\ref{lem:cosimplicial}, the retract of $X$  classified  by  $\tilde{X}$  is
        $n$-geometric.
    \end{proof}
\end{proposition}

\begin{lemma}\label{lem:retractslfpnstacks}
    If $X$ is a quasi-separated $n$-geometric sheaf that is locally of finite presentation, and if $Y$ is a
    retract of $X$, then $Y$ is locally of finite presentation.
    \begin{proof}
        By the previous lemma, $Y$ is itself $n$-geometric. Let $U=\coprod_i\Spec
        T_i\rightarrow Y$ be an atlas, and let $V=\coprod_i\Spec S_i\rightarrow X\times_Y U$
        be an atlas for the fiber product. Since the composition $V\rightarrow X\times_Y
        U\rightarrow X$ is an $(n-1)$-geometric submersion, it follows that $V$ is an atlas
        for $X$. By Lemma~\ref{lem:atlaseslfp}, each $S_i$ is locally of finite presentation
        over $R$. Taking the pullback of $X\times_Y U\rightarrow X$ over $Y\rightarrow X$, we get
        $U\rightarrow X\times_Y U$, since $Y\rightarrow X\rightarrow Y$ is the identity. Possibly by refining $U$, we can assume that
        $U\rightarrow X\times_Y U$ factors through the surjection $V\rightarrow X\times_Y
        U$. We thus have shown that each $T_i$ is a retract of $S_j$ for some
        $j$. Since $S_j$ is locally of finite presentation over $R$, it follows that $T_i$
        is as well.
    \end{proof}
\end{lemma}

We now prove in two lemmas that images of smooth $n$-geometric morphisms are Zariski open.
This is a generalization of the fact that images of smooth maps of schemes are Zariski open.
Restricting a sheaf in $\Shv_R^{\et}$ to discrete connective commutative rings induces a
geometric morphism of $\infty$-topoi $\pi_0^*:\Shv_R^{\et}\rightarrow\Shv_{\pi_0R}^{\et}$.
Note that $\pi_0^*\Spec S$ is equivalent to $\Spec\pi_0S$.

\begin{lemma}\label{lem:truncationsubobjects}
    Let $S$ be a connective commutative $R$-algebra. Then, a subobject $Z$ of $\Spec S$ is
    Zariski open if and only if $\pi_0^*Z$ is Zariski open in $\Spec\pi_0S$.
    \begin{proof}
        The necessity is trivial. So, suppose that $\pi_0^*Z$ is Zariski open.
        Because $\pi_0^*$ admits a left adjoint, we see that
        $\pi_0^*$ preserves $(-1)$-truncated objects and finite limits. Thus, $\pi_0^*$
        preserves subobjects, so $\pi_0^*Z$ is a subobject of $\Spec\pi_0S$. Let $F$ be the
        set of $f\in\pi_0S$ such that $\Spec S[f^{-1}]\rightarrow\Spec S$ factors through
        $Z$. Note that $f\in F$ if and only if
        $\Spec\pi_0S[f^{-1}]\rightarrow\Spec\pi_0S$ factors through $\pi_0^*Z$.
        By construction, there is a monomorphism over $\Spec S$
        \begin{equation*}
            W:=\bigcup_{f\in F}\Spec S[f^{-1}]\rightarrow Z.
        \end{equation*}
        Since $\pi_0^*Z$ is Zariski open, it follows that $\pi_0^*W=\pi_0^*Z$.
        The counit map of the adjunction
        \begin{equation*}
            \pi_0^*:\Shv_R^{\et}\rightleftarrows\Shv_{\pi_0R}^{\et}:\pi_{0*}
        \end{equation*}
        gives a map $Z\rightarrow\pi_{0*}\pi_0^*Z=\pi_{0*}\pi_0^*W$. Now, we can recover $W$
        from $\pi_{0*}\pi^*W$ as the pullback
        \begin{equation*}
            \begin{CD}
                W   @>>>    \Spec S\\
                @VVV        @VVV\\
                \pi_{0*}\pi_0^*W  @>>>    \pi_{0*}\pi_0^*\Spec S.
            \end{CD}
        \end{equation*}
        Indeed, since $W$ is open, it is a union of $\Spec S[f_i^{-1}]$.
        This is clear when $W$ is a basic open $\Spec S[f^{-1}]$, and the general case
        follows from the fact that $\pi_0^*$ induces an equivalence between the small
        Zariski site of $\Spec S$ and
        (the nerve of) the small Zariski site of $\Spec\pi_0 S$. Thus, there are maps $W\rightarrow Z$ and
        \begin{equation*}
            Z\rightarrow \pi_{0*}\pi_0^*Z\times_{\pi_{0*}\Spec\pi_0S}\Spec S\rwe W
        \end{equation*}
         over $\Spec S$. Since $W$ and $Z$ are subobjects of $\Spec S$, it follows that they are equivalent. Thus, $Z$ is Zariski open.
    \end{proof}
\end{lemma}

The image of a map $f:X\rightarrow Y$ of sheaves is defined as the epi-mono factorization
$X\twoheadrightarrow\im(f)\rightarrowtail Y$. In particular, the morphism
$\im(f)\rightarrowtail Y$ is a monomorphism.

\begin{lemma}\label{lem:smoothimages}
    Let $f:X\rightarrow Y$ be a smooth $n$-geometric morphism. Then, the map
    $\im(f)\rightarrow Y$ is a Zariski open immersion.
    \begin{proof}
        We may assume without loss of generality that $Y=\Spec S$ for some connective
        commutative $R$-algebra $S$. Then, by hypothesis, there is a smooth
        $(n-1)$-geometric chart
        \begin{equation*}
            \coprod_i\Spec T_i\rightarrow X
        \end{equation*}
        such that the compositions $g_i:\Spec T_i\rightarrow\Spec S$ are smooth, and thus have
        cotangent complexes $\Lrm_{g_i}$ which are projective.
        By~\cite{ha}*{Proposition~8.4.3.9}, $\pi_0\Lrm_{g_i}$ is the cotangent complex of
        $\pi_0^*(g_i):\Spec\pi_0T_i\rightarrow\Spec\pi_0S$, and it is projective. Since
        $g_i$ is locally of finite presentation, by Lemma~\ref{lem:lfp}, $\pi_0^*(g_i)$
        is locally of finite presentation, and hence smooth. Since smooth morphisms
        of discrete schemes are flat by~\cite{ega4}*{Theorem 17.5.1}, the image of $\pi_0^*(g_i)$ has an open image in
        $\Spec\pi_0S$. By the previous lemma, it follows that the image of $g_i$ in $\Spec
        S$ is open.
    \end{proof}
\end{lemma}

\subsection{Cotangent complexes of smooth morphisms}

In this section, we show that $n$-geometric morphisms have cotangent complexes, and we give a
criterion for an $n$-geometric morphism to be smooth in terms of formal smoothness and the
cotangent complex.

Let  $S$  be  a  connective  commutative  $R$-algebra,  and  let  $M$  be  a   connective
$S$-module.  Then, a small extension of $S$ by $M$ over $R$ is  a  connective  commutative
$R$-algebra $\tilde{S}$ together with an $R$-algebra section $d$ of $S\oplus \Sigma M\rightarrow S$
such that $\tilde{S}$ is equivalent to the pullback
\begin{equation*}
    \begin{CD}
        \tilde{S}   @>>>    S\\
        @VVV                @V(\id,0) VV\\
        S           @>d>>   S\oplus \Sigma M.
    \end{CD}
\end{equation*}
The $\infty$-category of small extensions $\CAlg_{R/S}^{\sm}$ is the full subcategory of $\CAlg_{R/S}^{\cn}$ spanned
by small extensions of $S$ over $R$.

\begin{lemma}
    There is a natural equivalence
    $\CAlg_{R/S}^{\sm}\we\left(\tau_{>0}\Mod_S\right)_{\Lrm_{R/S}/}$. The composite
    \begin{equation*}
        \CAlg_{R/S}^{\sm}\we\left(\tau_{>0}\Mod_S\right)_{\Lrm_{R/S}/}\rightarrow\tau_{>0}\Mod_S
    \end{equation*}
    is given by taking the cofiber $\Sigma M$ of $\tilde{S}\rightarrow S$.
    \begin{proof}
        By adjunction, the space of $R$-algebra sections of $S\oplus \Sigma M\rightarrow S$ is
        equivalent to the space of $S$-module maps $\Lrm_{R/S}\rightarrow\Sigma M$.
    \end{proof}
\end{lemma}

The previous lemma allows to compute the cotangent complex of a small extension.

\begin{lemma}\label{lem:cotansmall}
    Let $\tilde{S}\rightarrow S$ be a small extension of $S$ by $M$. Then, the cotangent complex $\Lrm_i$ of $i:\Spec
    S\rightarrow\Spec\tilde{S}$ is naturally equivalent to $\Sigma M$.
    \begin{proof}
        By the previous lemma, it suffices to show that $\tilde{S}$ is an initial object of
        $\CAlg_{\tilde{S}/S}^{\sm}$. But, it is easy to check that $\tilde{S}$ is a small extension of $S$ over
        $\tilde{S}$, and as $\tilde{S}\rightarrow S$ is the initial object of
        $\CAlg_{\tilde{S}/S}^{\cn}$, it follows that it is an initial object of
        $\CAlg_{\tilde{S}/S}^{\sm}$.
    \end{proof}
\end{lemma}

A sheaf $X$ is infinitesimally  cohesive  if  for  all  $R$-algebras  $S$  and  all  small
extensions $\tilde{S}\we S\times_{S\oplus \Sigma M} S$ of $S$ by an $S$-module $M$ the natural map
\begin{equation*}
    X(\tilde{S})\rightarrow X(S)\times_{X(S\oplus \Sigma M)} X(S)
\end{equation*}
is an equivalence.

\begin{lemma}\label{lem:lifting}
    Let $X$ be an infinitesimally cohesive sheaf with a cotangent complex $\Lrm_X$, let $u:\Spec S\rightarrow X$ be a point, and let $\tilde{S}\rightarrow
    S$ be a small extension of $S$ by $M$ classified by a class $x\in\pi_0\map_S(\Lrm_{\Spec S},\Sigma M)$. Then, $u$ extends to
    $\tilde{u}:\Spec\tilde{S}\rightarrow X$ if and only if the image of $x$ vanishes under the map
    $\pi_0\map_S(\Lrm_{\Spec S},\Sigma M)\rightarrow\pi_0\map_S(u^*\Lrm_X,\Sigma M)$ induced by $u^*\Lrm_X\rightarrow\Lrm_{\Spec S}$.
    \begin{proof}
        Since $X$ is infinitesimally cohesive, there is a cartesian square
        \[
            \begin{CD}
                X_u(\tilde{S})    @>>>    *\\
                @VVV                    @V(u,0) VV\\
                *    @>\alpha >>         X_u(S\oplus \Sigma M),
            \end{CD}
        \]
        where $X_u(\tilde{S})$ and $X_u(S\oplus \Sigma M)$ are the fibers of $X(\tilde{S})\rightarrow X(S)$ and $X(S\oplus
        \Sigma M)\rightarrow X(S)$ over $u$, and $\alpha$ is induced by $\Spec d:\Spec S\oplus M\rightarrow\Spec S$. By definition of the cotangent complex, $X_u(S\oplus\Sigma
        M)\we\map_S(u^*\Lrm_X,\Sigma M)$. So, $X_u(\tilde{S})$ is non-empty if and only if the point $\alpha$ of
        $\map_S(u^*\Lrm_X,\Sigma M)$ is $0$. But, $d$ is classified by $x$, so that $\alpha$ is the image of $x$ in
        $\map_S(u^*\Lrm_X,\Sigma M)$, as claimed.
    \end{proof}
\end{lemma}

A sheaf $X$ is nilcomplete  if  for  any  connective  commutative  $R$-algebra  $S$  the
canonical map
\begin{equation*}
    X(S)\rightarrow\lim_n X(\tau_{\leq n}S)
\end{equation*}
is an equivalence.
If  $T$   is   any   commutative   $R$-algebra,   then   $X=\Spec   T$   is   nilcomplete.
Indeed, if $S$ is a connective commutative $R$-algebra, then
\begin{equation*}
    X(S)=\map_{\CAlg_R}(T,S)\we\lim_n\map_{\CAlg_R}(T,\tau_{\leq
    n}S)=\lim_nX(\tau_{\leq n}S).
\end{equation*}
A map of sheaves $p:X\rightarrow Y$ is nilcomplete if
for all connective commutative $R$-algebras $S$ and all  $S$-points  of  $Y$  the  fiber
product $X\times_Y\Spec S$ is nilcomplete.

\begin{remark}\label{rem:smallsite}
    Suppose that $S\rwe\lim_\alpha S_\alpha$ is a limit diagram of connective
    commutative $R$-algebras such that each map $S\rightarrow S_\alpha$ induces an
    isomorphism on $\pi_0$. In this case, the underlying small \'etale $\infty$-topoi of $S$
    and each $S_\alpha$ are equivalent. Given a sheaf $X$ in $\Shv_R^{\et}$, let $X_S$
    (resp. $X_{S_\alpha}$) denote the restriction of $X$ to the small \'etale
    site of $S$. Thus, for instance, the space of global sections $X_{S_\alpha}(S)$ is equivalent
    to $X(S_\alpha)$. In order for $X(S)\rightarrow \lim_\alpha X(S_\alpha)$ to be an
    equivalence, it suffices to show that $X_S$ is equivalent to $\lim_\alpha X_{S_\alpha}$.
\end{remark}

As we now show, all $n$-geometric morphisms have cotangent complexes, and we use this to
show that the property of smoothness for $n$-geometric morphisms can be detected via a
Tor-amplitude condition on the cotangent complex. The proof of the next proposition is a mix
of several proofs in~\cite{dag0}, particularly those of Propositions 5.1.5 and 5.3.7.

\begin{proposition}\label{prop:smoothperfect}
    An $n$-geometric morphism $f:X\rightarrow Y$ is infinitesimally cohesive, nilcomplete,
    and has a $(-n)$-connective cotangent complex $\Lrm_{f}$.
    If $f$ is smooth, then $\Lrm_f$ is perfect of Tor-amplitude contained in $[-n,0]$.
    Finally, if $f$ is smooth, then it is formally smooth.
    \begin{proof}
        We prove the proposition by induction on $n$. We may assume moreover that $Y=\Spec
        S$, and prove the statements for $X$ and $\Lrm_X$. When $X$ is a disjoint union of
        affines, it is automatically infinitesimally cohesive and nilcomplete, since maps
        out of a commutative ring commute with limits. The other statements in the base case $n=0$ follow from
        Propositions~\ref{prop:affinecotan} and~\ref{prop:smoothprojective}.
        Thus, suppose the proposition  is  true  for
        $k$-geometric morphisms with $k<n$. Then, since $X$  is  $n$-geometric, we fix an  $(n-1)$-submersion
        $p:U\rightarrow  X$  where  $U$  is  a  disjoint   union   of   affines
        $\coprod_i\Spec T_i$. Write $p_i$ for the composition $\Spec T_i\rightarrow
        U\rightarrow X$.

        To prove the statements about infinitesimal cohesiveness and nilcompleteness, we
        apply Lemma \ref{lem:screwing} and use Remark~\ref{rem:smallsite}. Let $X$ be an
        $n$-geometric sheaf, and let $U\rightarrow X$ be a surjection of sheaves. Let
        $S\rwe \lim_\alpha S_{\alpha}$ be a limit diagram of connective commutative
        $R$-algebras such that each map $S\rightarrow S_{\alpha}$ induces an isomorphism on
        $\pi_0$. Consider the simplicial object obtained by taking iterated fiber products
        of the map
        \begin{equation*}
            \lim_\alpha U_{S_\alpha}\rightarrow\lim_\alpha X_{S_\alpha}.
        \end{equation*}
        By using identifications of the form
        \begin{equation*}
            U_{S_\alpha}\times_{X_{S_\alpha}}U_{S_{\alpha}}\we(U\times_X U)_{S_\alpha},
        \end{equation*}
        the simplicial objected induces a $(-1)$-truncated map from the geometric realization $|\lim_\alpha
        U_{S_\alpha,\bullet}|$ to $\lim_\alpha X_{S_\alpha}$. We obtain a commutative
        diagram
        \begin{equation*}
            \xymatrix{
                X_S \ar[r] &   \lim_\alpha X_{S_\alpha}\\
                |U_{S,\bullet}|\ar[u]\ar[r] &  |\lim_\alpha U_{S_\alpha,\bullet}|\ar[u]}
        \end{equation*}
        where the bottom map is an equivalence, the left vertical map is an equivalence, and
        the right vertical map is $(-1)$-truncated. To show the top map is
        an equivalence, it is enough to show that for any \'etale $S$-algebra $T$ the map $\lim_\alpha
        U(T_\alpha)\rightarrow\lim_\alpha X(T_\alpha)$ is surjective, where
        $T_\alpha=S_\alpha\otimes_S T$.

        \emph{Infinitesimal cohesiveness}. We specialize the above considerations to where
        $S$ is a small extension of $S_0$ by $M$
        \begin{equation*}
            \begin{CD}
                S   @>>>    S_0\\
                @VVV        @VVV\\
                S_0 @>>>    S_0\oplus\Sigma M.
            \end{CD}
        \end{equation*}
        Set $S_1=S_0\oplus\Sigma M$. We want to show that for any \'etale $S$-algebra $T$,
        the natural map
        \begin{equation*}
            \lim_i U(T_i)\rightarrow\lim_i X(T_i)
        \end{equation*}
        is surjective, where $T_i=S_i\otimes_S T$. It suffices to prove this when the
        value $\Spec(T_1)\rightarrow X$ at the terminal object lifts to $U$.
        To show that the map on limits is surjective, it suffices to show that for any
        point $x_0$ of $X(T_0)$ that maps to $x_1$ in $X(T_1)$, and for any lift of $x_1$ to
        $y_1$ in $U(T_1)$, there exists $y_0$ in $U(T_0)$ mapping to $y_1$ and $x_0$ (for
        either of the maps $T_0\rightarrow T_1$). This surjectivity follows from the fact
        that the cotangent complex of $U\rightarrow X$ exists, which is due to the inductive
        step of the proof. The surjectivity follows from Lemma~\ref{lem:lifting} because the
        cotangent complex of $U$ over $X$ is perfect and its dual is connective.

        \emph{Nilcompleteness}. The proof of nilcompleteness is similar to that of
        infinitesimal cohesiveness, and is left to the reader.

        \emph{Existence}. Fix a $T$-point $x:\Spec T\rightarrow X$. Since $p:U\rightarrow X$
        is surjective, we can assume that $x$ factors through $p$. Let $y:\Spec T\rightarrow
        U$ be such a factorization. Then, there is a natural morphism
        \begin{equation*}
            F:\der_{f\circ p}(y,M)\rightarrow\der_f(x,M).
        \end{equation*}
        The $0$-point of $\der_f(x,M)$ is the point in the fiber of $X(T\oplus M)\rightarrow
        X(T)\times_{Y(T)}Y(T\oplus M)$ corresponding to $\Spec T\oplus M\rightarrow\Spec
        T\rightarrow X$. The fiber over $0$ of $F$ is naturally equivalent to $\der_p(y,M)$.
        Thus, we have a natural fiber sequence
        \begin{equation*}
            \map_T(\Lrm_p,M)\rightarrow\map_T(\Lrm_{f\circ p},M)\rightarrow\der_f(x,M).
        \end{equation*}
        By the formal smoothness of the smooth $(n-1)$-geometric morphism $p$, the map of spaces $F$ is surjective. It
        follows that we can identify $\der_f(x,M)$ with the fiber of the delooped map
        \begin{equation*}
            \Brm\map_T(\Lrm_p,M)\rightarrow\Brm\map_T(\Lrm_{f\circ p},M).
        \end{equation*}
        Therefore, the fiber of $\Lrm_{f\circ p}\rightarrow\Lrm_p$ is a cotangent complex
        for $f$. The connectivity statement is immediate.

        \emph{Tor-amplitude}. Now, suppose that $X\rightarrow\Spec S$ is
        smooth.  Then, we may assume that $\Spec T_i$ is  smooth  over  $\Spec  S$; in
        particular  $\Spec T_i$ is locally  finitely  presented  and  $\Lrm_{T_i}$   is
        finitely generated projective. By
        Lemma~\ref{lem:cotansequence}, there is a cofiber sequence
        \begin{equation*}
            p_i^*\Lrm_X\rightarrow\Lrm_{\Spec T_i}\rightarrow\Lrm_{\Spec T_i/X}.
        \end{equation*}
        By the inductive hypothesis, $\Lrm_{\Spec T_i/X}$ is perfect and has Tor-amplitude contained in
        $[-n+1,0]$.  Therefore, $\Lrm_X$ is  perfect  and  has  Tor-amplitude  contained  in $[-n,0]$.

        \emph{Formal smoothness}. Let $\Kscr$ be the class of nilpotent thickenings
        $u:\tilde{T}\rightarrow T$ that satisfy the left lifting property with respect to
        $f$. Since $f$ has a cotangent complex, $\Kscr$ contains all trivial square-zero
        extensions $T\oplus M\rightarrow T$. To see that $\Kscr$ contains all small
        extensions of $T$ by $M$, we use the fact that
        \begin{equation*}
            X(\tilde{T})\we X(T)\otimes_{X(T\oplus\Sigma M)}X(T).
        \end{equation*}
        Therefore, to check that the projection
        \begin{equation*}
            X(\tilde{T})\we X(T)\otimes_{X(T\oplus\Sigma M)}X(T)\rightarrow X(T)
        \end{equation*}
        is surjective, it suffices to note that
        \begin{equation*}
            \pi_0X(T)\times_{\pi_0X(T\oplus\Sigma M)}\pi_0
            X(T)\rightarrow\pi_0X(T)\times_{\pi_0X(T)}\pi_0X(T)=\pi_0X(T)
        \end{equation*}
        is surjective, because the map of pullback diagrams admits a section induced by the
        inclusion $\pi_0X(T)\rightarrow\pi_0X(T\oplus\Sigma M)$. Finally, that $\Kscr$
        contains all nilpotent thickenings follows from the method of the proof
        of~\cite{dag7}*{Proposition 7.26}, which decomposes such a thickening as a limit of
        small extensions.
    \end{proof}
\end{proposition}

The fact that smooth $n$-geometric morphisms have perfect cotangent complexes with
Tor-amplitude contained in $[-n,0]$ characterizes smooth morphisms, at least if we include
the assumption that the morphism be locally of finite presentation.

\begin{proposition}
    An $n$-geometric morphism $f:X\rightarrow Y$ is smooth if and only if it is locally of
    finite presentation and $\Lrm_f$ is a perfect complex with Tor-amplitude contained in $[-n,0]$.
    \begin{proof}
        We may assume that $Y=\Spec S$. Let $U=\coprod_i\Spec T_i\rightarrow X$ be a smooth
        $(n-1)$-submersion onto $X$. Then, each composition $\Spec T_i\rightarrow\Spec S$ is
        smooth, and hence locally of finite presentation. Therefore, $f$ is locally of
        finite presentation. The fact that $\Lrm_f$ is perfect with Tor-amplitude contained in $[-n,0]$
        follows from the previous proposition. Suppose now that $f$ is $n$-geometric,
        locally of finite presentation, and that $\Lrm_f$ has Tor-amplitude contained in
        $[-n,0]$. Take a chart $p:U=\coprod_i\Spec T_i\rightarrow X$ where the $\Spec T_i$ are all
        locally of finite presentation over $\Spec S$. The pullback sequence
        \begin{equation*}
            p^*\Lrm_X\rightarrow\Lrm_U\rightarrow\Lrm_p
        \end{equation*}
        of cotangent complexes, together with the facts that $p^*\Lrm_X$ and $\Lrm_U$ are
        perfect complexes with Tor-amplitude contained in $[-n,0]$ and $[-n+1,0]$, respectively, shows that
        $\Lrm_U$ is perfect with Tor-amplitude contained in $[-n,0]$. But, since $U$ is a
        disjoint union of affines, $\Lrm_U$ is connective. Thus, $\Lrm_U$ is
        equivalent to a finitely generated projective module, so that each $\Spec
        T_i\rightarrow\Spec S$ is smooth, as desired.
    \end{proof}
\end{proposition}

\subsection{Etale-local sections of smooth geometric morphisms}

The theorem in this section says that smooth morphisms that are
surjective on geometric points are in fact surjections of \'etale sheaves.

\begin{theorem}\label{thm:lifts}
    If $p:X\rightarrow Y$ is a smooth locally geometric morphism that is surjective on
    geometric points,
    then for every $S$-point $\Spec S\rightarrow
    Y$ there exists an \'etale cover $\Spec T\rightarrow\Spec  S$  and  a  $T$-point  $\Spec
    T\rightarrow X$ such that
    \begin{equation*}
        \begin{CD}
            \Spec T @>>>    X\\
            @VVV            @VVV\\
            \Spec S @>>>    Y
        \end{CD}
    \end{equation*}
    commutes.
    \begin{proof}
        We may assume without loss of generality that $Y=\Spec S$, and it then suffices to
        prove that $X\rightarrow\Spec S$ has \'etale local sections. Write $X$ as a filtered
        colimit
        \begin{equation*}
            X\we\colim_i X_i
        \end{equation*}
        of $n_i$-geometric sheaves, such that each $X_i\rightarrow X$ is a monomorphism. By
        Lemma~\ref{lem:monocotan}, the cotangent complex $\Lrm_{X_i/X}$ vanishes. Suppose
        that the cotangent complex of $X$ has Tor-amplitude contained in $[-n,0]$. Then, for
        $i$ sufficiently large, it follows that $X_i\rightarrow\Spec S$ is a smooth
        $n_i$-geometric morphism. Since $\Spec S$ is quasi-compact, and since the image of
        $X_i$ in $\Spec S$ is open by Lemma~\ref{lem:smoothimages}, it follows that for some $i$, $X_i\rightarrow\Spec S$ is
        a smooth $n_i$-geometric morphism that is surjective on geometric points.
        There  exists  an
        $(n-1)$-submersion $U\rightarrow X_i$ such that $U$ is the disjoint  union  of  smooth
        affine $S$-schemes. Let
        $\pi_0^*U\rightarrow\Spec\pi_0 S$ be the
        associated map of ordinary schemes.  By hypothesis,  this  morphism  is  smooth  and
        is surjective on geometric points.  By~\cite{ega4}*{Corollaires IV.17.16.2 et IV.17.16.3(ii)}, there exists an \'etale
        cover  $\Spec\pi_0  T\rightarrow\Spec\pi_0  S$  and  a   factorization   $\Spec\pi_0
        T\rightarrow\pi_0^*U\rightarrow\Spec\pi_0 S$.
        The  \'etale  map  $\pi_0  S\rightarrow\pi_0  T$  determines  a  unique   connective
        commutative $S$-algebra $T$ by~\cite{ha}*{Theorem~8.5.0.6}. We would like to like  to  lift
        the   $\pi_0   T$-point   of   $\pi_0^*U$   to    a    $T$-point    of    $U$.
        Since  $U$  is  a  disjoint  union  of  affines,  it  is  nilcomplete.    Therefore,
        $U(T)\we\lim_n U(\tau_{\leq  n}T)$,  so  it  suffices  to  show  that  $U(\tau_{\leq
        n}T)\rightarrow U(\tau_{\leq n-1}T)$ is surjective. This is true since $U$ is formally
        smooth and $\tau_{\leq n}T\rightarrow \tau_{\leq n-1}T$ is a  nilpotent  thickening.
    \end{proof}
\end{theorem}

\section{Moduli of objects in linear \icats}\label{sec:moduli}

We study moduli spaces of objects in $R$-linear categories. This extends the work of
To\"en-Vaqui\'e~\cite{toen-vaquie} to the setting of commutative ring spectra. We
give some results on local moduli, which form the basis of an important geometricity
statement for global moduli sheaves. As a corollary, we show in the final section
that if $A$ is an Azumaya algebra over $R$, the sheaf of Morita equivalences from $A$ to $R$
is smooth over $\Spec R$, and hence has \'etale local sections.

\subsection{Local moduli}

In this section we prove the geometricity of the sheaf corepresented by a free commutative
$R$-algebra $\Sym_R(P)$ where $P$ is a perfect $R$-module, and we show that  when  $P$  has
Tor-amplitude contained in $[-n,0]$, then this sheaf is smooth. These facts are non-trivial precisely
because $\Sym_R(P)$ is not necessarily connective.  This turns out to be the  main  step  in
understanding the geometricity of more general moduli problems.

Let   $\ShProj_R$   denote    the    sheaf    of    finite    rank    projective    modules.

\begin{proposition}
    The sheaf $\ShProj_R$ is equivalent to $\coprod_n\B\GL_n$. In particular, $\ShProj_R$ is
    locally $1$-geometric and locally of finite presentation.
    \begin{proof}
        A   projective   module   is   locally   free   by    Theorem~\ref{prop:projectives}.
        Hence, the sheaf  of  projective  rank  $n$  modules  is  equivalent  to  $\B\GL_n$.
        This sheaf has a $0$-atlas $\Spec R\rightarrow\B\GL_n$, which shows that it is
        $1$-geometric and locally of finite presentation.
    \end{proof}
\end{proposition}

\begin{theorem}\label{thm:localmoduli}
    Let $P$ be a perfect $R$-module with Tor-amplitude contained in $[a,b]$ with $a\leq  0$.
    Then, the sheaf $\Spec\Sym_R(P)$
    is  a quasi-compact and quasi-separated $(-a)$-geometric  stack  that is locally  of   finite   presentation.    Moreover,   the
    cotangent complex of $\Spec\Sym_R(P)$ at an $S$-point $x:\Spec S\rightarrow\Spec\Sym_R(P)$ is the $S$-module
    \begin{equation*}
        \Lrm_{\Spec\Sym_R(P),x}\we P\otimes_R S.
    \end{equation*}
    Therefore, if $b\leq 0$, $\Spec\Sym_R(P)$ is smooth.
    \begin{proof}
        We prove all of the statements except for quasi-separatedness by induction on $-a$.
        If  $a=0$,  then  $P$  is  connective  by
        Proposition~\ref{prop:tor}, so
        that $\Sym_R(P)$ is connective as well.  Thus,  $\Spec\Sym_R(P)$  is $0$-geometric
        and quasi-compact. It is locally of finite presentation since the $R$-module $P$ is
        perfect.
        Now, assume that $P$  has  Tor-amplitude  contained  in  $[a,b]$  where  $a<0$.   By
        Proposition~\ref{prop:tor}, we can write $P$ as the
        fiber  of  some  map   $Q\rightarrow\Sigma^{a+1}N$,   where   $Q$   is   a   perfect
        $R$-module with Tor-amplitude contained in $[a+1,b]$ and $N$ is a finitely generated
        projective $R$-module. The fiber sequence induces a
        fiber sequence of sheaves
        \begin{equation*}
            \Spec\Sym_R(\Sigma^{a+1}N)\rightarrow\Spec\Sym_R(Q)\rightarrow\Spec\Sym_R(P)\rightarrow\Spec\Sym_R(\Sigma^aN),
        \end{equation*}
        where,   inductively,   both   $\Spec\Sym_R(\Sigma^{a+1}N)$   and   $\Spec\Sym_R(Q)$
        are  $(-a-1)$-geometric  stacks  that are locally  of  finite   presentation. The map
        \begin{equation}\label{eq:asdf}
            \Spec\Sym_R(Q)\rightarrow\Spec\Sym_R(P)
        \end{equation}
        is surjective, because if $S$ is a connective commutative $R$-algebra, we get a
        fiber sequence of spaces
        \begin{equation*}
            \map_R(\Sigma^{a+1}N,S)\rightarrow\map_R(Q,S)\rightarrow\map_R(P,S)\rightarrow\map_R(\Sigma^aN,S)\we\Brm\map_R(\Sigma^{a+1}N,S),
        \end{equation*}
        which shows that $\map_R(Q,S)$ is the total space of a principal bundle over
        $\map_R(P,S)$. The map~\eqref{eq:asdf} is also quasi-compact, since the fiber
        $\Spec\Sym_R(\Sigma^{a+1}N)$ is quasi-compact. Note   that
        \begin{align*}
            \Spec\Sym_R(Q)\times_{\Spec\Sym_R(P)}\Spec\Sym_R(Q)
            &\we\Spec\left(\Sym_R(Q)\otimes_{\Sym_R(P)}\Sym_R(Q)\right)\\
            &\we\Spec\Sym_R(Q\oplus_P Q)
        \end{align*}
        Using that the natural  map  given  by  an  inclusion  followed  by  the  codiagonal
        \begin{equation*}
            Q\rightarrow Q\oplus_P Q\rightarrow Q
        \end{equation*}
        is  the  identity,  it  follows  that  $Q\oplus_P   Q\we   Q\oplus   \Sigma^{a+1}N$.
        Therefore,
        \begin{equation*}
            \Spec\Sym_R(Q\oplus_P Q)\we\Spec\Sym_R(Q)\times_{\Spec
            R}\Spec\Sym_R(\Sigma^{a+1}N).
        \end{equation*}
        It follows that the projection  $\Spec\Sym_R(Q\oplus_P  Q)\rightarrow\Spec\Sym_R(Q)$
        is the pullback of a $(-a-1)$-geometric morphism, and so is itself
        $(-a-1)$-geometric. The projection is smooth because $\Spec\Sym_R(\Sigma^{a+1}N)$ is
        smooth. Therefore, by all of the statements of Lemma~\ref{lem:stacksdescent}, $\Spec\Sym_R(P)$ is a
        quasi-compact $(-a)$-geometric stack that is locally of finite
        presentation. Finally, by Lemma~\ref{lem:cotangent}, the cotangent complex of $\Spec\Sym_R(P)$  is
        $P\otimes_R\Sym_R(P)$, so $\Spec\Sym_R(P)$ is smooth
        by~\ref{prop:smoothperfect} if $b\leq 0$.

        It remains to show that $\Spec\Sym_R(P)$ is quasi-separated. Let $Q$ be the cofiber
        of the diagonal morphism $P\rightarrow P\oplus P$. Then, the fiber of the diagonal
        morphism $$\Spec\Sym_R(P)\rightarrow\Spec\Sym_R(P)\times_{\Spec
        R}\Spec\Sym_R(P)\we\Spec\Sym_R(P\oplus P)$$ is equivalent to $\Spec\Sym_R(Q)$, which
        is quasi-compact by the first part of the proof.
    \end{proof}
\end{theorem}

\begin{remark}
    If $P$ is a perfect $R$-module with Tor-amplitude contained in $[a,b]$ and $a\geq 0$, then
    it also has Tor amplitude  contained  in  $[0,b]$,  and  the  proposition  implies  that
    $\Spec\Sym_R(P)$  is  a  $0$-geometric  stack.
\end{remark}

\subsection{The moduli sheaf of objects}

In this section, we apply the study of  local  moduli  above  to  global  moduli  sheaves  of
objects. The main theorems in this section, Theorems~\ref{thm:shm} and~\ref{thm:moduli}, are generalizations of results
of~\cite{toen-vaquie} to connective $\EE_\infty$-ring spectra.

Let $\Cscr$ be a compactly generated $R$-linear $\infty$-category. Define a functor
$$\StM_\Cscr:(\Aff_{R}^{\cn})^{\op}\rightarrow\LCI$$ whose value at $\Spec S$ is the full subcategory of
$\Dual_R\Cscr\otimes_R\Mod_S\we\Fun_R^{\Lrm}(\Cscr,\Mod_S)$
consisting   of   those   objects $f$ such that for every compact object $x$ of
$\Cscr$, the value $f(x)$ is a compact object of $\Mod_S$. Put another way, we
can define $\StM_\Cscr$ as the pullback
\begin{equation*}
    \begin{CD}
        \StM_{\Cscr}    @>>>    \Fun_R^{\Lrm}(\Cscr,\StMod)\\
        @VVV                    @VVV\\
        \prod_{x\in\pi_0\Cscr^{\omega}}\StMod^{\omega}  @>>> \prod_{x\in\pi_0\Cscr^{\omega}}\StMod,
    \end{CD}
\end{equation*}
where $\StMod^{\omega}$ is the functor of compact objects
$\StMod^{\omega}:(\Aff_R^{\cn})^{\op}\rightarrow\LCI$ given by
\begin{equation*}
    \StMod^{\omega}(\Spec S)=\Mod_S^{\omega}.
\end{equation*}

\begin{lemma}\label{lem:stm}
    For any compactly generated $R$-linear \icat  $\Cscr$,  the  functor  $\StM_\Cscr$
    satisfies \'etale hyperdescent.
    \begin{proof}
        It is clear that $\Fun_R^{\Lrm}(\Cscr,\StMod)$ satisfies \'etale hyperdescent since $\StMod$ is an \'etale hyperstack.
        Moreover, we claim that the functor of compact objects $\StMod^{\omega}$ also
        satisfies \'etale hyperdescent. It suffices to check that
        $\Mod_S^{\perf}\rightarrow\lim_{\Delta}\Mod_{T^{\bullet}}^{\perf}$ is an
        equivalence for any \'etale hypercover $S\rightarrow T^{\bullet}$. But this follows from the commutative diagram
        \begin{equation*}
            \begin{CD}
                \Mod_S^{\omega}  @>>>    \lim_{\Delta}\Mod_{T^{\bullet}}^{\omega}\\
                @VVV                    @VVV\\
                \Mod_S          @>>>    \lim_{\Delta}\Mod_{T^{\bullet}}.
            \end{CD}
        \end{equation*}
        The vertical arrows are fully faithful, and the bottom arrow is an equivalence. It
        follows that the top arrow is fully faithful. It is also essentially surjective for
        the following reason. The compact objects of $\Mod_S$ are precisely the dualizable
        ones. But, the property of being dualizable can be checked \'etale locally.
        Thus, $\StMod^{\omega}$ satisfies \'etale hyperdescent.
        Now, by the pullback definition of $\StM_{\Cscr}$ above, it follows that
        $\StM_\Cscr$ satisfies \'etale hyperdescent.
    \end{proof}
\end{lemma}

Because it satisfies \'etale hyperdescent, the functor $\StM_\Cscr$ extends uniquely to a
limit preserving functor $\Shv_R^\et\rightarrow\LCI$. We abuse notation and write
$\StM_\Cscr$ for the resulting stack.
Let $\ShM_\Cscr=\StM_\Cscr^{\eq}$ be the sheaf of equivalences in $\StM_\Cscr$. We call
this   the   moduli   sheaf   (or    moduli    space)    of    objects    in    $\Cscr$.
It is a sheaf of small spaces because $\Cscr^\omega$ is small.
If $\Cscr$ is $\Mod_A$ for some $R$-algebra $A$, we  also  write  $\Mscr_A$  for
$\Mscr_{\Mod_A}$. This sheaf classifies left $A$-module structures on perfect $S$-modules.

\begin{definition}
    Let $\ShM_R^{[a,b]}$ be the  full  sub-sheaf  of  $\ShM_R$  whose  $S$-points  consist
    of perfect $S$-modules with Tor-amplitude contained in $[a,b]$.  Note  that  this  makes
    sense since Tor-amplitude is stable  under  base-change  by  Proposition~\ref{prop:tor}.
    By the same proposition, there is an equivalence
    \begin{equation*}
        \colim_{a\leq b}\ShM_R^{[a,b]}\rwe\ShM_R,
    \end{equation*}
    and each map $\ShM_R^{[a,b]}\rightarrow\ShM_R$ is a monomorphism.
\end{definition}

\begin{theorem}\label{thm:shm}
    The sheaf $\ShM_R$ is locally  geometric  and  locally  of  finite  presentation.
    \begin{proof}
        By  the  definition  of  local  geometricity,  it  suffices  to   show   that   each
        $\ShM_R^{[a,b]}$ is $(n+1)$-geometric and locally of finite presentation where $n=b-a$.
        To begin, we show that each diagonal map
        \begin{equation}\label{eq:diagonal401}
            \ShM_R^{[a,b]}\rightarrow\ShM_R^{[a,b]}\times_{\Spec R}\ShM_R^{[a,b]}
        \end{equation}
        is   $(b-a)$-geometric   and locally   of   finite    presentation.
        Given a map from $\Spec S$ to the product classifying two  perfect  $S$-modules  $P$
        and  $Q$,  the  pullback  along  the  diagonal  is  equivalent   to   $\Eq(P,Q)$,
        the sheaf over $\Spec S$ classifying equivalences between $P$ and $Q$.
        This is a  Zariski  open  subsheaf  of  $\Spec\Sym_S(P\otimes_S  Q^{\vee})$.   Since
        $P\otimes_S   Q^{\vee}$   has   Tor-amplitude   contained   in    $[a-b,b-a]$,    by
        Theorem~\ref{thm:localmoduli}, $\Eq(P,Q)\rightarrow\Spec S$ is $(b-a)$-geometric. Therefore,
        the   diagonal map~\eqref{eq:diagonal401}  is   $(b-a)$-geometric,   as   desired.

        We now proceed by induction on $n=b-a$.  When $n=0$, $a$-fold suspension induces  an
        equivalence
        \begin{equation*}
            \ShProj_R\rightarrow\ShM_R^{[a,a]}.
        \end{equation*}
        By Corollary~\ref{ex:gln}, $\StProj_R$ is $1$-geometric and locally of finite presentation.
        Now, suppose that $\ShM_R^{[a+1,b]}$ is $(b-a)$-geometric  and  locally  of  finite
        presentation.  The general outline of the induction is as follows.  We  construct  a
        sheaf $U$ that admits a $0$-geometric smooth map to
        $\ShM_R^{[a+1,b]}\times_{\Spec R}\ShM_R^{[a+1,a+1]}$ and use this to show that $U$ is a
        $(b-a)$-geometric sheaf locally of finite presentation.
        Then, we show that $U$ submerses onto $\ShM_R^{[a,b]}$.
        By Lemma~\ref{lem:the}, this is enough to conclude that
        $\ShM_R^{[a,b]}$ is $(b-a+1)$-representable and  locally  of  finite  presentation.

        Let   $U$   be   defined   as   the   pullback   of    sheaves
        \begin{equation*}
            \begin{CD}
                U   @>>>    \Fun(\Delta^1,\StMod_R^{\omega})^{\eq}\\
                @Vp VV        @VVV\\
                \ShM_R^{[a+1,b]}\times_{\Spec R}\ShM_R^{[a+1,a+1]}    @>>>
            \Fun(\partial\Delta^1,\StMod_R^{\omega})^{\eq}.
            \end{CD}
        \end{equation*}
        Suppose that $\Spec S\rightarrow\ShM_R^{[a+1,b]}\times_{\Spec R}\ShM_R^{[a+1,a+1]}$ is
        a point classifying a perfect $S$-module $Q$ of Tor-amplitude contained in $[a+1,b]$
        and a perfect $S$-module $\Sigma^{a+1}M$ of Tor-amplitude contained in $[a+1,a+1]$. The fiber of
        $p$ at this point is simply the local moduli sheaf
        \begin{equation*}
            \Spec\Sym_S(Q\otimes_S \Sigma^{-a-1}M).
        \end{equation*}
        As $Q\otimes_S \Sigma^{-a-1}M$ has Tor-amplitude contained in $[0,b-a-1]$, it follows that
        this local moduli sheaf is $0$-geometric and
        locally of finite presentation because $\Sym_S Q\otimes_S \Sigma^{-a-1}M$ is a compact
        commutative $S$-algebra (because $Q\otimes_S \Sigma^{-a-1}M$ is compact).
        Thus, $p$ is $0$-geometric and locally of finite presentation. Moreover,
        $\ShM_R^{[a+1,b]}\times_{\Spec R}\ShM_R^{[a+1,a+1]}$ is a $(b-a)$-geometric sheaf locally of
        finite presentation by the inductive hypothesis. So, $U$ is $(b-a)$-geometric by
        Lemma~\ref{lem:lurie}, and it is locally of finite presentation.

        Let $q:U\rightarrow\ShM_R^{[a,b]}$ be the map that sends an object of  $U$  to  the
        fiber of the map it classifies in $\Fun(\Delta^1,\StMod_R^{\omega})^{\eq}$, which  has
        the asserted Tor-amplitude by part~(5) of Lemma~\ref{prop:tor}. Since $U$ is $(b-a)$-geometric
        and because the diagonal of $\ShM_R^{[a,b]}$ is $(b-a)$-geometric, it follows  from
        Lemma~\ref{lem:the} that $q$ is $(b-a)$-geometric.  If we prove that $q$  is  smooth
        and surjective, it will follow that $\ShM_R^{[a,b]}$ is $(b-a+1)$-geometric by
        Lemma~\ref{lem:the}.

        The surjectivity of $q$ is immediate from part~(7) of Proposition~\ref{prop:tor}. To
        check smoothness, consider a point $\Spec S\rightarrow\ShM^{[a,b]}$, which classifies
        a compact $S$-module $P$ of Tor-amplitude contained in $[a,b]$. Let $Z$ be the
        fiber product of this map with $U\rightarrow\ShM^{[a,b]}$. The $T$-points of the sheaf $Z$ may be
        described as ways of writing $P\otimes_S T$ as a fiber of a map
        $Q\rightarrow\Sigma^{a+1}M$ where $M$ is a
        finitely generated projective $T$-module, and $Q$ is a $T$-module with Tor-amplitude contained in $[a+1,b]$.
        Possibly after passing to a Zariski cover of $\Spec T$, we may assume that $M\we T^{\oplus r}$ is finitely
        generated and free. In other words, the sheaf $Z$ decomposes as $$Z\we\coprod_{r\geq 0}
        Z_r,$$ where $Z_r$ classifies maps $\Sigma^{a}S^{\oplus r}\rightarrow P$ with
        cofiber having Tor-amplitude contained in $[a+1,b]$. Since
        $\Spec\Sym_S(\Sigma^a(P^{\vee})^{\oplus r})$ classifies all maps $\Sigma^a S^{\oplus
        r}\rightarrow P$, we see that $Z_r$ consists of the points of
        $\Spec\Sym_S(\Sigma^a (P^{\vee})^{\oplus r})$
        classifying maps $\Sigma^a S^{\oplus r}\rightarrow P$ that are surjective on
        $\pi_a$. Since $\pi_aP$ is finitely generated by Proposition~\ref{prop:compactness},
        this surjectivity condition is open, because the vanishing of the cokernel of $\pi_0S^{\oplus r}\rightarrow\pi_aP$ can
        be detected on fields by Nakayama's lemma.
        As the perfect module $(P^{\vee})^{\oplus r}$ has Tor-amplitude contained in
        $[a-b,0]$, $\Spec\Sym_S(\Sigma^a(P^{\vee})^{\oplus r})$ is smooth by
        Theorem~\ref{thm:localmoduli}. Thus, $Z_r$ is
        smooth, and hence the morphism $U\rightarrow\ShM^{[a,b]}$ is smooth, which completes
        the proof.
    \end{proof}
\end{theorem}

To analyze $\ShM_A$ for other rings $A$,
we use sub-sheaves $\ShM_A^{[a,b]}$ of $\ShM_A$ induced by the corresponding
sub-sheaves of $\ShM_R$: specifically, define
$\ShM_A^{[a,b]}$ to be the pullback in
\begin{equation*}
    \begin{CD}
        \ShM_A^{[a,b]}   @>\pi^{[a,b]}>>    \ShM_R^{[a,b]}\\
        @VVV                    @VVV\\
        \ShM_A           @>\pi >>    \ShM_R.
    \end{CD}
\end{equation*}
Since the filtration $\{\ShM_R^{[-a,a]}\}_{a\geq 0}$ exhausts $\ShM_R$, it follows that
$\{\ShM_A^{[-a,a]}\}_{a\geq 0}$ exhausts $\ShM_A$.

\begin{proposition}
    Let $\Mod_A$ be an $R$-linear category of finite  type,  so  that  $A$  is  a  compact
    $R$-algebra. Then, the natural morphism
    $\pi:\ShM_A^{[a,b]}\rightarrow\ShM_R^{[a,b]}$ is $(b-a)$-geometric and locally of finite presentation.
    \begin{proof}
        It is easy to see using Corollary~\ref{cor:algebrafiber} that the space of $T$-points of the fiber of $\pi^{[a,b]}$ at a point $\Spec
        S\rightarrow\ShM_R^{[a,b]}$ classifying a perfect $S$-module $P$ is equivalent to
        \begin{equation*}
            \map_{\Alg_T}(A\otimes_R T,\End_T(P\otimes_S T)).
        \end{equation*}
        We will write $\Shmap(A\otimes_R S,\End_S(P))$ for the resulting sheaf over $\Spec S$.
        Now, since $A\otimes_R S$ is of finite presentation as an $S$-algebra, $A\otimes_R S$ is a retract of a finite colimit of
        the  free  $S$-algebra  $S\langle  t  \rangle$.   It  follows from
        Lemma~\ref{lem:nlimits}, Proposition~\ref{prop:retractsnstacks}, and
        Lemma~\ref{lem:retractslfpnstacks} that  to  prove  that
        $\Shmap(A\otimes_R S,\End_S(P))$ is $(b-a)$-geometric and locally of finite
        presentation, it suffices to show that $\Shmap(S\langle t\rangle,\End_S(P))$
        is a quasi-separated $(b-a)$-geometric sheaf that is locally of finite presentation. But,
        \begin{equation*}
            \Shmap_{\Alg_S}(S\langle t \rangle,\End_S(P))\we\Shmap_{\Mod_S}(S,\End_S(P))\we\Spec\Sym_S(\End_S(P)^{\vee}).
        \end{equation*}
        Since $\End_S(P)^{\vee}\we P^{\vee}\otimes_S P$ is  perfect  and  has  Tor-amplitude
        contained in $[a-b,b-a]$,  it follows from  Theorem~\ref{thm:localmoduli}  that  the
        fiber is $(b-a)$-geometric, quasi-separated, and locally of finite presentation.
    \end{proof}
\end{proposition}

Given the proposition, it is now straightforward to prove the following theorem.

\begin{theorem}\label{thm:moduli}
    Let   $A$   be   a compact   $R$-algebra.    Then,    the    stack
    $\ShM_A$   is   locally    geometric    and    locally    of    finite    presentation,
    and the functor $\pi:\ShM_A\rightarrow\ShM_R$ is locally geometric and locally of finite
    presentation.
    \begin{proof}
        By the previous proposition, $\ShM_A^{[a,b]}\rightarrow\ShM_R^{[a,b]}$ is
        $(b-a)$-geometric and locally of finite presentation. Since $\ShM_R^{[a,b]}$ is
        $(b-a+1)$-geometric and locally of finite presentation, it follows from
        Lemma~\ref{lem:lurie}, that $\ShM_A^{[a,b]}$ is also $(b-a+1)$-geometric and locally
        of finite presentation. It follows that $\ShM_A$ is locally geometric and locally of
        finite presentation. To prove the second statement, let $\Spec S\rightarrow\ShM_R$
        classify a perfect $S$-module $P$, which has Tor-amplitude contained in some
        $[a,b]$. The fiber of $\pi$ over this point is equivalent to the fiber of
        $\pi^{[a,b]}$ over $P$, which by the previous proposition is $(b-a)$-geometric and
        locally of finite presentation. 
    \end{proof}
\end{theorem}

Note that, in the proof, the fiber is not only locally geometric,
but $(b-a)$-geomtric. However, the geometricity of the fibers changes from point to
point.

\begin{corollary}\label{cor:cotangent}
    Let $A$ be  a  compact  $R$-algebra,  and  let  $\Spec  S\rightarrow\ShM_A$
    classify a perfect $S$-module $P$ with a left $A$-module structure. Then, the cotangent complex of
    $\ShM_A$ at the point $P$ is the $S$-module
    \begin{equation*}
        \Lrm_{\ShM_A,P}\we\Sigma^{-1}\End_{A^{\op}\otimes_R S}(P)^{\vee}.
    \end{equation*}
    \begin{proof}
        By Lemma~\ref{lem:monocotan} and Proposition~\ref{prop:smoothperfect}, the cotangent complex $\Lrm_{\ShM_A}$
        exists. Consider  the  loop sheaf $\Omega_{P}\ShM_A$.   By  Lemma~\ref{lem:colimcotan},  the
        cotangent  complex  of  this  sheaf  at  the  basepoint  $*$  is   simply
        \begin{equation*}
            \Lrm_{\Omega_{P}\ShM_A,*}\we\Sigma\,\Lrm_{\ShM_A,P}.
        \end{equation*}
        Thus, it suffices to compute $\Lrm_{\Omega_{P}\ShM_A,*}$. The stack $\Omega_{P}\ShM_A$
        is described by
        \begin{align*}
            \Omega_{P}\ShM_A(T)\we\aut_{A^{\op}\otimes_R T}(P\otimes_S
            T)&\subseteq\Omega^{\infty}\End_{A^{\op}\otimes_R T}(P\otimes_S T)\\
            &\we\map_S\left(S,\End_{A^{\op}\otimes_R T}(P\otimes_S T)\right)\\
            &\we\map_S\left(S,\End_{A^{\op}\otimes_R S}(P)\otimes_S T\right)\\
            &\we\map_S(\End_{A^{\op}\otimes_R S}(P)^{\vee},T)\\
            &\we\map\left(\Spec  T,\Spec\Sym_S(\End_{A^{\op}\otimes_R  S}(P)^{\vee})\right),
        \end{align*}
        where   the   equivalence   between   the   second   and   third   lines    is    by
        Lemma~\ref{lem:functoriality}.
        The inclusion is Zariski open, and hence the computation of Theorem~\ref{thm:localmoduli} says that
        \begin{equation*}
            \Lrm_{\Omega_{P}\ShM_A,*}\we \End_{A^{\op}\otimes_R S}(P)^{\vee},
        \end{equation*}
        which completes the proof.
    \end{proof}
\end{corollary}

\subsection{\'Etale local triviality of Azumaya algebras}

Let $R$ be a connective $\EE_\infty$-ring spectrum, and let $A$ be an Azumaya $R$-algebra.
We prove now that the stack of Morita equivalences from $A$ to $R$ is smooth and surjective
over $\Spec R$. As a corollary, we obtain one of the major theorems of the paper: the \'etale
local triviality of Azumaya algebras.

Since $A$ is compact as an $R$-module, the space $\ShM_A(S)$ of $S$-points is equivalent to
the space
\begin{equation*}
    \ShM_A(S)\we\Fun_{\Cat_{S,\omega}}\left(\Mod_{A\otimes_R S},\Mod_S\right)^{\eq}
\end{equation*}
of (compact object preserving) functors between compactly-generated $S$-linear categories.
We define the full subsheaf $\ShMor_A\subseteq\ShM_A$ of Morita equivalences from $A$ to $R$
by restricting the $S$-points to the full subspace of $\ShM_A(S)$ consisting of the equivalences $\Mod_{A\otimes_RS}\we\Mod_S$.

\begin{proposition}\label{prop:moritasmoothness}
    Suppose that $R$ is a connective $\EE_\infty$-ring  and  that  $A$  is  an
    Azumaya algebra. The sheaf $\ShMor_A\rightarrow\Spec R$ of Morita equivalences is locally
    geometric and smooth.
    \begin{proof}
        We show that $\ShMor_A\subseteq\ShM_A$ is quasi-compact and Zariski open. Fix  an  $S$-valued  point  of
        $\ShM_A$ classifying an $A^{\op}\otimes_R S$-module $P$ which is compact as an
        $S$-module. The bimodule $P$ defines an adjoint pair of functors
        \begin{equation*}
            -\otimes_AP:\Mod_{A\otimes_RS}\rightleftarrows\Mod_S:\Map_S(P,-).
        \end{equation*}
        To show that $\ShMor_A\subseteq\ShM_A$ is open, it suffices to check that the subsheaf of points of $\Spec S$ on which the unit
        \begin{equation*}
            \eta(A):A\rightarrow\Map_S(P,A\otimes_A P)
        \end{equation*}
        and counit
        \begin{equation*}
            \epsilon(S):\Map_S(P,S)\otimes_A P\rightarrow S
        \end{equation*}
        morphisms  are  equivalences is open in $\Spec S$,  since  the  unit  and   counit
        transformations   are
        equivalences   if   and   only   if   they   are   equivalences    on    generators.
        As    $A$    is    a    perfect    $S$-module    by    the    Azumaya    hypothesis,
        Proposition~\ref{prop:vanishing} implies that
        the cofibers of these maps vanish on quasi-compact Zariski open subschemes of $\Spec
        S$. Taking the intersection of these two open subschemes yields the desired quasi-compact Zariski  open  subscheme
        of   $\Spec   S$   on   which    $P$ defines a Morita equivalence.
        It follows that $\ShMor_A$ is locally geometric and locally of finite  presentation.
        
        Given a point $P:\Spec  S\rightarrow\ShMor_A$,  the  cotangent  complex  at  $P$  is
        \begin{equation*}
            \Lrm_{\ShMor_A,P}\we\Sigma^{-1}\End_{A\otimes_R S}(P)\we \Sigma^{-1}S,
        \end{equation*}
        a perfect $S$-module with Tor-amplitude contained in $[-1,-1]$. Thus, by definition,
        $\ShMor_A$ is a smooth locally geometric sheaf.
    \end{proof}
\end{proposition}

The following theorem is a generalization of~\cite{toen-derived}*{Proposition~2.14} to
connective $\EE_\infty$-ring spectra.

\begin{theorem}\label{thm:etalelocaltriviality}
    Let $R$ be a connective $\EE_\infty$ ring spectrum, and let $A$ be an Azumaya
    $R$-algebra. Then, there is a faithfully flat \'etale $R$-algebra $S$ such that
    $A\otimes_R S$ is Morita equivalent to $S$.
    \begin{proof}
        The theorem follows immediately from the previous proposition, Theorem~\ref{thm:brk}, and
        Theorem~\ref{thm:lifts}.
    \end{proof}
\end{theorem}

\begin{proposition}
    If $R$ is a connective $\EE_\infty$-ring spectrum and $A$ is an Azumaya $R$-algebra,
    then the sheaf of Morita equivalences $\ShMor_A$ is a $\ShPic$-torsor. In particular, it
    is $1$-geometric and smooth.
    \begin{proof}
        The action of $\ShPic$ on $\ShMor_A$ is simply by tensoring $A^{\op}$-modules with
        line bundles. \'Etale-locally, $\ShMor_A$ is equivalent to the space of equivalences
        $\Mod_S\we\Mod_S$ by the theorem. This is precisely $\ShPic$ over $\Spec S$.
    \end{proof}
\end{proposition}

\section{Gluing generators}\label{sec:gluing}

The main result of the previous section shows that Azumaya algebras over connective
$\EE_\infty$-rings are \'etale locally trivial. In this section, we want to show that certain \'etale cohomological
information on derived schemes $X$ can be given by Azumaya algebras. In other words we want
to prove that ``$\Br(X)=\Br'(X)$'' in various cases. This is established
once we prove the following theorem.

\begin{theorem}[Local-global principle]
    Let $\Cscr$ be an $S$-linear category, and suppose that $S\rightarrow T$ is an \'etale
    cover such that $\Cscr\otimes_S T$ has a compact generator. Then, $\Cscr$ is has a compact generator.
\end{theorem}

In fact, we prove a version of this theorem for quasi-compact and quasi-separated
derived schemes.
The result we prove
expands on~\cite{dag11}*{Theorem~6.1}, which is about a similar statement for the
property of being compactly generated.

The proof breaks up into several parts. First, we prove a local-global statement for Zariski covers. This is used in
two ways: to reduce the problem from schemes to affine schemes and to help prove Ninsevich descent. Second, we prove a local-global
statement for \'etale covers, following~\cite{dag11}*{Section~6}. The main insight there is to use the fact that a presheaf is
an \'etale sheaf if and only if it is a sheaf for the Nisnevich and finite \'etale topologies. Then, by using a theorem
of Morel and Voevodsky (see~\cite{dag11}*{Theorem~2.9}), we can reduce the proof of Nisnevich local-global principle to certain special squares, which we
analyze directly using techniques that, essentially, go back to Thomason~\cite{thomason-trobaugh} and
B\"okstedt-Neeman~\cite{bokstedt-neeman}. Proof of a local-global
principle for finite \'etale covers is subsumed in a more general statement for finite and
flat covers, which is purely $\infty$-categorical.

This theorem has applications to the module theory of perfect stacks, as is developed in~\cite{bzfn} and~\cite{dag11}, and is
related to questions about when derived categories are compactly generated, which have been
studied by Thomason~\cite{thomason-trobaugh},
Neeman~\cite{neeman-1992} and~\cite{neeman-1996},
B\"okstedt-Neeman~\cite{bokstedt-neeman}, and Bondal-van den
Bergh~\cite{bondal-vandenbergh}.

With two major exceptions, the outline of the proof
is already contained in~\cite{toen-derived}. First, our proof differs significantly from To\"en's
when it comes to the \'etale local-global principle. Since To\"en works with simplicial commutative rings, he is able to
use some concrete constructions based on work of Gabber~\cite{gabber} to reduce to the finite \'etale case. These constructions, which
involve algebras of invariants of symmetric groups acting on polynomial rings and quotient
algebras, simply fail in the case of $\EE_\infty$-ring spectra.
Thus, we use Lurie's idea of using the Morel-Voevodsky result to prove the \'etale
local-global principle. Second, we
cannot prove the fppf version contained in~\cite{toen-derived}. Because of the lack of
$\EE_\infty$-structures on quotient rings, we do not know how to show that the stack of
quasi-sections used in the proof of~\cite{toen-derived}*{Proposition~4.13}. Hence, we work
everywhere in the \'etale topology. For the cases of interest to us, this restriction is not
a problem.

\subsection{Azumaya algebras and Brauer classes over sheaves}

Fix a connective $\Einfinity$-ring spectrum $R$.
In Section~\ref{sec:sheaves}, we introduced  the  \'etale hyperstacks $\StAlg$, $\StAlg^{\Az}$, and
$\StCat_R^{\desc}$. Let $\ShAlg$, $\ShAz$, and $\ShP$ be the associated underlying (large) \'etale
hypersheaves.
There are natural maps $\ShAlg\rightarrow\ShP$ and $\ShAz\rightarrow\ShP$. Let $\ShMr$ and
$\ShBr$ be the \'etale hypersheafifications of the images of these maps.
To be precise,
for every connective commutative $R$-algebra $S$, there is a map $\ShAlg(S)\rightarrow\ShP(S)$,
and the image of this map is full subspace of $\ShP(S)$ consisting of those points of
$\ShP(S)$ in the image of $\ShAlg(S)$. As $S$ varies, these images determine a presheaf of
spaces, and $\ShMr$ is the \'etale sheafification of this presheaf. The story for $\ShBr$ is
similar, but with $\ShAz$ in place of $\ShAlg$.

\begin{definition}
    Let $X$ be an object of $\Shv_R^{\et}$.
    \begin{enumerate}
        \item[(1)]  A quasi-coherent algebra over $X$ is a morphism $X\rightarrow\ShAlg$.
        \item[(2)]  An  Azumaya  algebra  over  $X$  is  a   morphism   $X\rightarrow\ShAz$.
        \item[(3)]  A   Morita   class   over    $X$    is    a    morphism
            $X\rightarrow\ShMr$.
        \item[(4)]  A Brauer class over $X$ is a morphism $X\rightarrow\ShBr$.
        \item[(5)]  A linear category with descent over $X$ is a morphism $X\rightarrow\ShP$.
    \end{enumerate}
    Note that of the above classifying spaces, only $\ShAz$ and $\ShBr$ are actually sheaves
    of \emph{small} spaces.
\end{definition}

A Brauer class over $X$ is thus a linear category over $X$ which is \'etale locally
equivalent to modules over some Azumaya algebra.  The rest of this section will
prove that every Brauer class (resp. Morita class) over $X$ lifts to an Azumaya algebra
(resp. algebra) when $X$ is a quasi-compact and quasi-separated derived scheme.

If $\alpha:\Spec S\rightarrow\ShP$ is  a linear category with descent  over  $\Spec  S$,  let
$\Mod_S^{\alpha}$ denote the $\Mod_S$-module classified by $\alpha$ by the Yoneda lemma.
The following construction is studied extensively in~\cite{bzfn} and~\cite{dag11}.

\begin{definition}
    Let $\alpha:X\rightarrow\ShP$ be a linear category with descent over  $X$. Then,  the  \icat  of
    $\alpha$-twisted $X$-modules is
    \begin{equation*}
        \Mod_X^{\alpha}=\lim_{f:\Spec S\rightarrow X}\Mod_S^{\alpha\circ f}.
    \end{equation*}
    This limit exists and $\Mod_X^{\alpha}$ is stable and presentable, because $\Prl_{\st}$ is closed under limits.
\end{definition}

We can describe this construction as the right Kan extension of $(\Aff_{/X}^{\cn})^{\op}\simeq\CAlg_R^{\cn}\rightarrow\Prl$, the functor
which  sends  $f:\Spec  S\rightarrow  X$  to  $\Mod_S^{\alpha\circ  f}$,  evaluated   at
$X$. As a particularly important case, let $\alpha:X\rightarrow\ShP$ be the linear category with descent over $X$ which sends a point $x:\Spec S\rightarrow
X$ to $\Mod_S^{x^*\alpha}=\Mod_S$. Then, by definition,
$\Mod_X^{\alpha}$ is simply $\Mod_X$, which is an $\Einfinity$-algebra in $\Prl$. When $X$
is a (non-derived) scheme, then the homotopy category of $\Mod_X$ recovers the usual
derived category $D_\mathrm{qc}(X)$ of complexes of $\Oscr_X$-modules with quasi-coherent
cohomology sheaves.

The construction of the stable $\infty$-category of $\alpha$-twisted modules commutes
with colimits.

\begin{lemma}\label{lem:limitstwisted}
    Let $I\rightarrow\Shv_R^{\et}$ be a small diagram of sheaves $X_i$ with colimit $X$, let
    $\alpha:X\rightarrow\ShP$ be a linear category with descent over  $X$,  and  let  $\alpha_i$  be
    the restriction of $X$ to $X_i$. Then, the canonical map
    \begin{equation*}
        \Mod_X^{\alpha}\rightarrow\lim_I\Mod_{X_i}^{\alpha_i}
    \end{equation*}
    is an equivalence.
    \begin{proof}
        This follows from our definition of
        $\Mod_{X}^{\alpha}$ as a right Kan extension.
    \end{proof}
\end{lemma}

To attack our main theorem, the local-global principle, we require some additional
terminology.

\begin{definition}\label{def:perfectobject}
    Let $\alpha:X\rightarrow\ShP$ be a linear category with descent over $X$.
    \begin{enumerate}
        \item   An object $P$ of $\Mod_X^{\alpha}$ is called perfect if for every point $x:\Spec S\rightarrow X$, $x^*
            P$ is a compact object of $\Mod_S^{x^*\alpha}$.
        \item   An object  $P$ of $\Mod_{X}^{\alpha}$ is a perfect generator if for every point
            $x:\Spec S\rightarrow X$, the pullback $x^*P$ is a compact generator of $\Mod_{S}^{x^*\alpha}$.
        \item   An object $P$ is a global generator of $\Mod_X^{\alpha}$ if it is a compact
            generator and a perfect generator.
    \end{enumerate}
\end{definition}

Note that while perfect objects are preserved automatically by any pullback induced by a
map $\pi:X\rightarrow Y$ in $\Shv_X^{\et}$, it is not the case that compact objects are
preserved by pullbacks. For instance, if $X$ is not quasi-compact over the base $\Spec R$, then
$\Mod_R\rightarrow\Mod_X$ sends $R$ to $\Oscr_X$, which is perfect but might not be compact.
It is for this reason why perfect objects play such an important role. However, in most
cases of interest, it is possible to show that the perfect and compact objects do coincide.
See, for example,~\cite{bzfn}*{Section 3}.

When $X$ is affine, the next lemma shows that there is no difference between the notions of
compact generators and perfect generators of $\Mod_X^{\alpha}$. In particular, every perfect
generator is automatically a global generator.

\begin{lemma}\label{lem:compactsonaffines}
    If $\alpha:\Spec S\rightarrow\ShP$ is a linear category with descent, then an object $P$ of
    $\Mod_S^\alpha$ is a compact generator if and only if it is a perfect generator.
    \begin{proof}
        If $P$ is a perfect generator, then $P$
        is a compact generator of $\Mod_S^{\alpha}$ by definition. So, suppose that $P$ is a
        compact generator of $\Mod_S^{\alpha}$. We must show that for any $f:\Spec T\to\Spec S$ is a connective
        $\Einfinity$-$S$-algebra, then $P\otimes_S T$ is a compact generator of
        $\Mod_T^{\alpha}$. There is a commutative diagram of equivalences
        \begin{equation*}
            \begin{CD}
                \Mod_S^{\alpha}\otimes_S\Mod_T  @>>>    \Mod_T^{f^*\alpha}\\
                @V\Map(P,-) VV                                @VVV\\
                \Mod_{\End(P)^{\op}}\otimes_S\Mod_T @>>>    \Mod_{\End(P)^{\op}\otimes_S T},
            \end{CD}
        \end{equation*}
        in which the
        left-hand equivalence is Morita Theory~\ref{thm:morita} and the right-hand
        equivalence is induced from the other three.
        By commutativity, the object $P\otimes_S T$ in the upper-left corner is sent to
        $\End(P)^{\op}\otimes_S T$ in the lower-right corner, which is indeed a compact
        generator.
    \end{proof}
\end{lemma}

The following lemma will be used below to detect when an object is a compact generator of
$\Mod_S^{\alpha}$ by passing to $\Mod_T^{\alpha}$ for an \'etale cover $S\rightarrow T$.

\begin{lemma}\label{lem:flatdetection}
    If $S\rightarrow T$ is an \'etale cover, and if $\alpha:\Spec S\rightarrow\ShP$ is an
    linear category with descent, then a \emph{compact} object $P$ of $\Mod_S^\alpha$ is a compact generator of
    $\Mod_S^{\alpha}$ if and only if $P\otimes_S T$ is a compact generator of $\Mod_T^\alpha$.
    \begin{proof}
        One direction is clear: if $P$ is a compact generator of $\Mod_S^{\alpha}$, then by
        the lemma above, it is a perfect generator, so that $P\otimes_S T$ is a compact
        generator of $\Mod_T^{\alpha}$. So, suppose that $P$ is a compact object of
        $\Mod_S^{\alpha}$ such that $P\otimes_S T$ is a compact generator of
        $\Mod_T^{\alpha}$. Let $A=\End_S(P)^{\op}$, and let $\Mod_A$ be the stable
        $\infty$-category of $A$-modules. Write $T^\bullet$ for the cosimplicial commutative
        $S$-algebra associated to the cover $S\rightarrow T$. Consider commutative diagram
        \begin{equation*}
            \begin{CD}
                \Mod_S^{\alpha} @>>>    \lim_{\Delta}\Mod_{T^{\bullet}}^{\alpha}\\
                @V\Map(P,-)VV           @V\Map(P\otimes_S T^{\bullet},-) VV\\
                \Mod_A          @>>>    \lim_{\Delta}\Mod_{A\otimes_S T^{\bullet}}.
            \end{CD}
        \end{equation*}
        The horizontal maps are equivalences since both $\Mod_S^{\alpha}$ and $\Mod_A$
        satisfy \'etale descent, the latter by Example~\ref{ex:moddescent}. On the other hand, since
        $P\otimes_S T$ is a compact generator of $\Mod_T^{\alpha}$, the right vertical map
        is an equivalence, since it is the limit of a level-wise equivalence of simplicial
        $\infty$-categories. It follows that the left vertical map is an equivalence. In
        particular, if $\Map(P,M)\we 0$, then $M\we 0$ in $\Mod_S^{\alpha}$. Thus, $P$ is a
        compact generator of $\Mod_S^{\alpha}$.
    \end{proof}
\end{lemma}

As we now see, the linear categories with descent over $X$ that possess perfect generators
are exactly those which are $\infty$-categories of modules for quasi-coherent algebras over $X$. Our strategy in
proving the local-global principle is then to construct perfect generators.

\begin{proposition}\label{prop:factorization}
    A linear category with descent $\alpha:X\rightarrow\ShP$ factors through $\ShAlg\rightarrow\ShP$ if
    and only if $\Mod_X^{\alpha}$ possesses a perfect generator.
    \begin{proof}
        Suppose that $\alpha:X\rightarrow\ShP$ factors as
        \begin{equation*}
            X\xrightarrow{\Ascr}\ShAlg\rightarrow\ShP.
        \end{equation*}
        Then, there is an algebra object $A$ in $\Mod_X^{\alpha}$, which restricts to an $S$-algebra
        $A\otimes S$ over each affine $\Spec S\rightarrow X$, and which is a compact
        generator of the $S$-linear category $\Mod_S^\alpha\we\Mod_{A\otimes S}$. Hence, $A$ is a perfect generator. Now, suppose
        that $P$ is a perfect generator of $\Mod_X^{\alpha}$.
        By hypothesis,
        for any point $x:\Spec S\rightarrow X$, the object $P$ of
        $\Mod_S^{x^*\alpha}$ induces an equivalence
        \begin{equation*}
            \Map(P,-):\Mod_S^{x^*\alpha}\rightarrow\Mod_{\End(P)^{\op}\otimes S}.
        \end{equation*}
        In other words, we obtain a natural
        equivalence of functors
        \begin{equation*}
            \Mod_{\Spec -/X}^{\alpha}\rightarrow\Mod_{\End(P)^{\op}\otimes -}.
        \end{equation*}
        Therefore, $\End(P)^{\op}$ classifies a lift of $\alpha$ through $\ShAlg\rightarrow\ShP$.
    \end{proof}
\end{proposition}

\subsection{The Zariski local-global principle}

There is a long history to the arguments in this section. On the one hand, the ideas about lifting
compact objects along localizations goes back to Thomason-Trobaugh~\cite{thomason-trobaugh}
and Neeman~\cite{neeman-1992}*{Theorem~2.1}. On the other hand, the arguments about Zariski gluing appeared
in B\"okstedt-Neeman~\cite{bokstedt-neeman}*{Section 6}, in an argument about derived categories of
quasi-coherent sheaves. They were further used in~\cite{neeman-1996}*{Proposition~2.5}
and~\cite{bondal-vandenbergh}*{Theorem~3.1.1}
before being used by To\"en~\cite{toen-derived}*{Proposition~4.9} for module categories over quasi-coherent sheaves of algebras.

Given a colimit-preserving functor $F:\Cscr\rightarrow\Dscr$ of stable presentable $\infty$-categories, the
kernel of $F$ is full subcategory of $\Cscr$ consisting of those objects which become equivalent to $0$ in $\Dscr$. Since the $\infty$-category of stable presentable $\infty$-categories has limits which are computed in $\widehat{\Cat}_\infty$, we see that the kernel of $F$ is stable, presentable, and equipped with a colimit-preserving inclusion into $\Cscr$.

In this section, when $U$ is a quasi-compact open subscheme of a derived scheme $X$, we will
write $\Mod_{X,Z}^{\alpha}$ for the kernel of
$\Mod_X^{\alpha}\rightarrow\Mod_U^{\alpha}$, where $Z$ is the complement of $U$ in $X$.
Of course, this complement will usually not exist as a derived scheme, but only as a closed subspace of $X$.

The following proposition, which appears essentially in~\cite{bokstedt-neeman},
is one of the major points of ``derived'' geometry in our proof
that $\Br(X)=\Br'(X)$. The generator $K$ in the proof is truly a derived object, and thus produces,
even in the case of ordinary schemes, a derived Azumaya algebra.

\begin{proposition}[\cite{bokstedt-neeman}*{Proposition 6.1} and \cite{toen-derived}*{Proposition 3.9}]\label{prop:localization}
    Let $j:U\subset X=\Spec S$ be a quasi-compact open subscheme with complement $Z$, and let
    $\alpha:X\rightarrow\ShP$ be a $S$-linear category such that $\Mod_X^{\alpha}$ has a
    compact generator $P$. Then, the restriction functor
    $j^*:\Mod_X^\alpha\rightarrow\Mod_U^\alpha$ is a localization whose kernel
    $\Mod_{X,Z}^{\alpha}$ is generated by a single compact object $L$ in $\Mod_X^\alpha$.
    \begin{proof}
        Note that under these hypotheses, it is enough to treat the special case in which
        $\alpha$ classifies $\Mod_S$. Indeed, in this case, we have a localization sequence
        \begin{equation*}
            \Mod_{X,Z}\rightarrow\Mod_X\rightarrow\Mod_U.
        \end{equation*}
        Since $\Mod_X^{\alpha}$ is dualizable (it admits a compact generator), tensoring with $\Mod_X^{\alpha}$ 
        preserves limits, and we obtain the localization sequence
        \begin{equation*}
            \Mod_{X,Z}^{\alpha}\rightarrow\Mod_X^{\alpha}\rightarrow\Mod_U^{\alpha}.
        \end{equation*}
        To complete the proof, it suffices to show that
        $\Mod_{X,Z}$ has a compact generator.
        Write
        \begin{equation*}
            U=\bigcup_{i=1}^r \Spec S[f_i^{-1}],
        \end{equation*}
        and let $K_i$ be the cone of $S\xrightarrow{f_i} S$. Then,
        $K=K_1\otimes_S\cdots\otimes_S K_r$ is a compact object of
        $\Mod_X^{\alpha}$, and $j^*L\we 0$. We claim that $K$ is a compact generator of the
        kernel of $\Mod_{X,Z}$. Suppose that $\Map(K,M)\we 0$ and $j^*(M)\we 0$. Then,
        \begin{equation*}
            \Map(K_1,\Map(K_2\otimes_S\cdots\otimes_S K_r,M)\we\Map(K,M)\we 0,
        \end{equation*}
        where we are using the fact that $K_i$ is self-dual up to a shift. It follows that $f_1$ acts
        invertibly on $\Map(K_2\otimes_S\cdots\otimes_S K_r,M)$, so that
        \begin{align*}
            \Map(K_2,\otimes_S\cdots\otimes_S K_r,M)    &\we \Map(K_2,\otimes_S\cdots\otimes_S K_r,M)\otimes_S S[f_1^{-1}]\\
                                                        &\we \Map\left(K_2,\otimes_S\cdots\otimes_S K_r,M\otimes_S S[f_1^{-1}]\right)\\
                                                        &\we 0,
        \end{align*}
        where the last equivalence follows from the fact that $j^*M\we 0$ and that $\Spec
        S[f_1^{-1}]$ is contained in $U$. By induction, it follows that
        \begin{equation*}
            \Map(K_r,M)\we 0,
        \end{equation*}
        and thus that
        \begin{equation*}
            M\we M\otimes_S S[f_r^{-1}]\we 0.
        \end{equation*}
        Therefore, $K$ is a compact generator of $\Mod_{X,Z}$.
    \end{proof}
\end{proposition}

We also need the following $K$-theoretic characterization, due to~\cite{thomason-trobaugh}
in the case of schemes and~\cite{neeman-1992} more generally, of when an object lifts through a
localization. Recall that if $\Cscr$ is a compactly generated stable $\infty$-category, then
$\K_0(\Cscr)$ is the Grothendieck group of the compact objects of $\Cscr$. That is, it is
the free abelian group on the \emph{set} of compact objects of $\Cscr$, modulo the relation
$[M]=[L]+[N]$ whenever there is a cofiber sequence $L\rightarrow M\rightarrow N$.
Note that $\K_0(\Cscr)$ depends only on the triangulated homotopy category $\Ho(\Cscr)$.

\begin{proposition}\label{prop:decomposition}
    Let $\alpha:X\rightarrow\ShP$ be a linear category such that $\Mod_X^{\alpha}$ is
    compactly generated, where $X$ is a derived scheme over $R$ which can be embedded as a
    quasi-compact open subscheme of an affine $\Spec S\in\Aff_R^{\cn}$, and let $U\subseteq X$ be a
    quasi-compact open subscheme. Then,
    a compact object $P$ of $\Mod_U^{\alpha}$ lifts to $\Mod_X^{\alpha}$ if and only if it is in the image of
    $\K_0(\Mod_X^{\alpha})\rightarrow\K_0(\Mod_U^{\alpha})$.
    \begin{proof}
        This follows from Neeman's localization theorem~\cite{neeman-1992}*{Theorem~2.1} and
        its corollary~\cite{neeman-1992}*{Corollary~0.9}. The only thing to check is that
        $\Mod_{X,Z}^{\alpha}$ is compactly generated by a set of objects that are compact in
        $\Mod_X^{\alpha}$. For this, we refer to the beginning of the proof of
        Lemma~\ref{lem:closedgluing}, which shows that the inclusion
        $\Mod_{X,Z}^{\alpha}\rightarrow\Mod_X^{\alpha}$ preserves compact objects.
    \end{proof}
\end{proposition}

We are now ready to state and prove our Zariski local-global principle, which is a
generalization of the arguments of~\cite{bokstedt-neeman}*{Section 6} and the
theorems~\cite{bondal-vandenbergh}*{Theorem 3.1.1}
and~\cite{toen-derived}*{Proposition~4.9}.

\begin{theorem}\label{thm:zariskigluing}
    Let $X$ be a quasi-compact, quasi-separated derived scheme over $R$, and let
    $\alpha:X\rightarrow\ShP$ be a linear category with descent over $X$. If there exists Zariski cover
    $f:\Spec S\rightarrow X$ such that $\Mod_S^{f^*\alpha}$ has a compact generator, then
    there exists a global generator of $\Mod_X^{\alpha}$.
    \begin{proof}
        The proof is by induction on $n$, the number of affines in an open cover of $X$ over which
        there are compact generators. If $\Mod_X^{\alpha}$ has a compact generator when $X=\Spec S$,
        then it has a global generator by Lemma~\ref{lem:compactsonaffines}. Now, assume that for all quasi-compact,
        quasi-separated derived schemes $Y$ and all $\beta:Y\rightarrow\ShP$, if $$\coprod_{i=1}^n\Spec
        S_i\xrightarrow{\coprod f_i} Y$$ is a Zariski cover such that $\Mod_{S_i}^{f_i^*\beta}$
        has a compact generator for $i=1,\ldots, n$, then $\Mod_Y^{\beta}$ has a global
        generator. Let $$\coprod_{i=1}^{n+1}\Spec T_i\xrightarrow{\coprod
        g_i} X$$ be a Zariski cover such that each $\Mod_{T_i}^{g_i^*\alpha}$ has a compact generator.
        The proof will be complete if we produce a global generator of $\Mod_X^{\alpha}$.

        Let $Y$ be the union of $\Spec T_i$, $i=1,\ldots,n$ in $X$, let $Z=\Spec
        T_{n+1}$, and let $W=Y\cap Z$. So, there is a pushout square of sheaves
        \begin{equation*}
            \begin{CD}
                W   @>>>    Z\\
                @VVV        @VVV\\
                Y   @>>>    X.
            \end{CD}
        \end{equation*}
        By Lemma~\ref{lem:limitstwisted}, it follows that
        \begin{equation}\label{eq:cartesianmodules}
            \begin{CD}
                \Mod_X^{\alpha} @>>>    \Mod_Z^{\alpha}\\
                @VVV                    @VVV\\
                \Mod_Y^{\alpha} @>>>    \Mod_W^{\alpha}
            \end{CD}
        \end{equation}
        is a pullback square of stable presentable $\infty$-categories. By the induction hypothesis, there exists a
        global generator $P_Y$ of $\Mod_Y^{\alpha}$. The restriction of
        $P_Y\oplus\Sigma P_Y$ to $W$ lifts to a compact object of $\Mod_Z^{\alpha}$ by
        Proposition~\ref{prop:decomposition}. Since $P_Y\oplus\Sigma P_Y$ is also a compact
        generator of $\Mod_Y^{\alpha}$, we can assume in fact that the restriction $P_W$ of $P_Y$ to $W$ lifts to a
        compact object $P_Z$ of $\Mod_Z^{\alpha}$. The cartesian square~\eqref{eq:cartesianmodules} says
        that there is an object $P_X$ of $\Mod_X^{\alpha}$ 
        that restricts to $P_W$, $P_Y$, and $P_Z$ over $W$, $Y$, and $Z$, respectively. The
        object $P_X$ is in fact compact, because for any $M_X$ in $\Mod_X^{\alpha}$, the
        mapping space $\map(P_X,M_X)$ is computed as the pullback
        \begin{equation*}
            \begin{CD}
                \map_X(P_X,M_X) @>>>    \map_Z(P_Z,M_Z)\\
                @VVV                    @VVV\\
                \map_Y(P_Y,M_Y) @>>>    \map_W(P_W,M_W).
            \end{CD}
        \end{equation*}
        Since finite limits commute with filtered colimits, since the restriction functors
        preserve colimits, and since $P_Z$, $P_Y$, and $P_W$ are compact, it follows that
        $P_X$ is compact.

        Because $Z$ is affine and $W\subset Z$ is quasi-compact, by
        Proposition~\ref{prop:localization}, the restriction functor
        \begin{equation*}
            \Mod_Z^{\alpha}\rightarrow\Mod_W^{\alpha}
        \end{equation*}
        is a localization, which kills exactly the stable subcategory of $\Mod_Z^{\alpha}$
        generated
        by a compact object $Q_Z$ of $\Mod_Z^{\alpha}$. We may lift $Q_Z$ to an
        object $Q_X$ of $\Mod_X^{\alpha}$ that restricts to $0$ over $Y$
        using~\eqref{eq:cartesianmodules}. The object $Q_X$ is compact for the same reason
        that $P_X$ is compact. Then, $L_X=P_X\oplus Q_X$ is a compact object of
        $\Mod_X^{\alpha}$, which we claim is a global generator of $\Mod_X^{\alpha}$.

        Suppose that $M_X$ is an object of $\Mod_X^{\alpha}$ such that $\map_X(L_X,M_X)\we 0$. Then,
        $\map_X(P_X,M_X)\we 0$ and $\map_X(Q_X,M_X)\we 0$. For any $N$ in $\Mod_X^{\alpha}$, we have a cartesian square
        \begin{equation*}
            \begin{CD}
                \map_X(N_X,M_X)   @>>>    \map_Z(N_Z,M_Z)\\
                @VVV                    @VVV\\
                \map_Y(N_Y,M_Y)\we 0  @>>>    \map_W(N_W,M_W).
            \end{CD}
        \end{equation*}
        When $N_X=Q_X$, the bottom mapping spaces are trivial
        so that $0\we\map_X(Q_X,M_X)\we\map_Z(Q_Z,M_Z)$. It follows that $M_Z$ is supported
        on $W$. On the other hand, $0\we\map_X(P_X,M_X)\we\map_Y(P_Y,M_Y)$ since in that
        case, the right-hand vertical map is an equivalence as $M_Z$ is supported on $W$. As
        $P_Y$ is a compact generator of $\Mod_Y^{\alpha}$, the restriction of $M$ to $U$ is
        trivial. But, the support of $M_Z$ is contained in $W\subset U$, so $M$ is trivial.
        Therefore, $L$ is a compact generator of $\Mod_X^{\alpha}$.

        To prove that $L_X$ is a perfect generator of $\Mod_X^{\alpha}$, it
        suffices to show that $L_Y$ is a perfect generator of $\Mod_Y^{\alpha}$ and that $L_Z$ is a compact
        generator of $\Mod_Z^{\alpha}$ (since $Z$ is affine). Indeed, given any affine
        $V=\Spec S$ mapping into $X$, we can intersect it with the affine hypercover
        determined by the $T_i$. Write $S\rightarrow T$ for this map. By hypothesis,
        $L\otimes T$ is a compact generator for $\Mod_T^{\alpha}$. By
        Lemma~\ref{lem:flatdetection}, it
        follows that $L\otimes S$ is a compact generator of $\Mod_S^{\alpha}$.

        That $L_Y$ is a global generator of $\Mod_Y^{\alpha}$ is trivial,
        since $Q_Y\we 0$ and so $L_Y\we P_Y$ was chosen to be a global generator of $\Mod_Y^{\alpha}$.
        If $M$ is an object of $\Mod_Z^{\alpha}$ such that $\map_Z(L_Z,M)\we 0$,
        then $\map_Z(Q_Z,M)\we 0$ so that $M$ is supported on $W$. Thus,
        \begin{equation*}
            0\we\map_Z(P_Z,M)\we\map_W(P_W,M)\we\map_Y(P_Y,M).
        \end{equation*}
        But, $P_Y$ is a global generator of $\Mod_Y^{\alpha}$, and
        $\Mod_W^{\alpha}\rightarrow\Mod_Y^{\alpha}$ is fully faithful. Thus, $M\we 0$.
    \end{proof}
\end{theorem}

\subsection{The \'etale local-global principle}

In this section, we adapt an idea of Lurie to show that for $R$-linear $\infty$-categories, the property of
having a compact generator is local for the \'etale topology. The context of this section is slightly different from
that of the rest of Section~\ref{sec:gluing}: we do not require that our $R$-linear
categories to satisfy \'etale hyperdescent. As
every $R$-linear $\infty$-category satisfies \'etale descent by~\cite{dag11}*{Theorem~5.4} (not \'etale
hyperdescent), this is a natural hypothesis to drop when considering \'etale covers. So, instead of studying morphisms
$X\rightarrow\ShP$, we instead fix an $R$-linear category $\Cscr$. If $S$ is a commutative $R$-algebra, we write $\Mod_S(\Cscr)$
for the $\infty$-category of $S$-modules in $\Cscr$. In particular, $\Mod_R(\Cscr)\we\Cscr$,
and more generally $\Mod_S(\Cscr)\we\Cscr\otimes_R S$. For a general \'etale sheaf $X$,
we define
\[
    \Mod_X(\Cscr)=\lim_{\Spec S\rightarrow X}\Mod_S(\Cscr).
\]
If $\Cscr$ is in fact a linear category with \'etale hyperdescent arising from a map $\alpha:\Spec R\rightarrow\ShP$, then these definitions
agree with our definitions of $\Mod_X^{\alpha}$ above.

\begin{lemma}
    Let $F:\Cscr\rightleftarrows\Dscr:G$ be a pair of adjoint functors between stable
    presentable $\infty$-categories such that the right adjoint $G$ is conservative and preserves
    filtered colimits. If $P$ is a compact generator of $\Cscr$, then $F(P)$ is a compact generator of $\Dscr$.
    \begin{proof}
        Since $G$ preserves filtered colimits,
        $F$ preserves compact objects, so that $F(P)$ is compact. Suppose that $M$ is an
        object of $\Dscr$ such that $\Map_{\Dscr}(F(P),M)\we 0$. Then,
        $\Map_{\Cscr}(P,G(M))\we 0$. Since $P$ is a compact generator of $\Cscr$, this
        implies that $G(M)\we 0$. The conservativity of $G$ implies that $M\we 0$, so
        that $F(P)$ is a compact generator of $\Dscr$.
    \end{proof}
\end{lemma}

Following Lurie, we let $\Test_{\pi_0 R}$ be the category of (non-derived) $\pi_0 R$-schemes $X$ which admit a quasi-compact open
immersion $X\rightarrow\Spec\pi_0 S$, where $\pi_0 S$ is an \'etale $\pi_0 R$-algebra.
There is a Grothendieck topology on $\Test_{\pi_0 R}$ that extends the Ninsevich topology~\cite{dag11}*{Proposition~2.7}. Lurie
proves~\cite{dag11}*{Theorem~2.9} a version of the theorem of Morel and Voevodsky which says that for a presheaf $F$ on
$\Test_{\pi_0 R}$, being a Nisnevich sheaf is equivalent to satisfying affine Nisnevich
excision. Recall that $F$ satisfies affine Nisnevich excision if $F(\emptyset)$ is contractible and
for all affine morphisms $X'\rightarrow X$ and quasi-compact open subschemes $U\subseteq X$ such that $X'-U'\rightarrow X-U$ is an
isomorphism, where $U'=X'\times_X U$, the diagram
\begin{equation*}
    \begin{CD}
        F(X)    @>>>    F(X')\\
        @VVV            @VVV\\
        F(U)    @>>>    F(U')
    \end{CD}
\end{equation*}
is a pullback square of spaces.

Let $\CAlg_R^{\et}$ denote the $\infty$-category of \'etale $R$-algebras.
There is a fully faithful embedding $\CAlg_R^{\et}\rightarrow \Nrm(\Test_{\pi_0R}^{\op})$ given by sending $S$ to
$\Spec\pi_0 S$. Given an $R$-linear category
$\Cscr$, we extend the construction that sends an \'etale $R$-algebra $S$ to $\Mod_S(\Cscr)$
to $\Test_{\pi_0 R}$ by
right Kan extension. In other words, if $X$ is an object of $\Test_{\pi_0 R}$,
\begin{equation*}
    \Mod_X(\Cscr)=\lim_{\Spec\pi_0 S\rightarrow X}\Mod_S(\Cscr),
\end{equation*}
where the limit runs over all \'etale $R$-algebras $S$ and all maps $\Spec\pi_0 S\rightarrow X$.

If $j:U\subseteq X$ is a quasi-compact open immersion in $\Test_{\pi_0 R}$ with complement
$Z$, viewed as a $\pi_0 R$-scheme with its reduced scheme structure, then we let
$\Mod_{X,Z}(\Cscr)$ be the full subcategory of $\Mod_X(\Cscr)$ consisting of those objects $M$ such that $j^*M\we 0$ in
$\Mod_U(\Cscr)$. Roughly speaking, these are the quasi-coherent $\Oscr_X$-modules in $\Cscr$ with
support contained in $Z$.

\begin{lemma}\label{lem:closedgluing}
    Let $X$ be an object of $\Test_{\pi_0 R}$, and let $j:U\rightarrow X$ be a quasi-compact open immersion with complement $Z$.
    If there exists a compact object $Q$ in $\Mod_X(\Cscr)$ such that $\Mod_U(\Cscr)$ is generated by $j^* Q$
    and if $\Mod_{X,Z}(\Cscr)$ has a compact generator $P$, then $i_!P\oplus Q$ is a compact
    generator of $\Mod_X(\Cscr)$, where $i_!$ is the inclusion functor from
    $\Mod_{X,Z}(\Cscr)$ into $\Mod_X(\Cscr)$.
    \begin{proof}
        Since $Q$ is compact by hypothesis, to show $i_!P\oplus Q$ is compact, we must show
        that $i_!P$ is compact. In fact, we show that $i_!$ preserves compact objects. To
        see this, consider the right adjoint $i^!$ of $i_!$, which is defined as the fiber
        of the natural unit natural transformation
        \begin{equation*}
            i^!\rightarrow\id_{\Mod_X(\Cscr)}\rightarrow j_*j^*,
        \end{equation*}
        where $j_*$ is the right adjoint of $j^*$. The functor $j^*$, being a left adjoint,
        preserves small colimits. By~\cite{dag11}*{Proposition~5.15}, the functor $j_*$
        preserves small colimits as well (this is where the quasi-compact hypothesis is
        used). Since $i^!$ is defined via a finite limit diagram,
        it follows that $i^!$ preserves filtered colimits, and hence that $i_!$ preserves
        compact objects. Hence, $i_!P\oplus Q$ is a compact object of $\Mod_X(\Cscr)$.
        Suppose now that $M$ is an object of $\Mod_X(\Cscr)$ such that
        \begin{equation*}
            \Map_X(i_!P\oplus Q,M)\we 0.
        \end{equation*}
        Then, $\Map_Z(P,i^!M)\we\Map_X(i_!P,M)\we 0$. Since $P$ is a compact generator of
        $\Mod_{X,Z}(\Cscr)$, it follows that $i^!M\we 0$. Hence the unit map $M\rightarrow
        j_*j^*M$ is an equivalence. At the same time,
        \begin{equation*}
            \Map_U(j^*Q,j^*M)\we\Map_X(Q,j_*j^*M)\we 0.
        \end{equation*}
        Thus, $j^*M\we 0$, since $j^*Q$ is a compact generator of $\Mod_U(\Cscr)$. Since
        $M\we j_*j^*M$, we see that $M\we 0$. Therefore, $i_!P\oplus Q$ is indeed a compact
        generator of $\Mod_X(\Cscr)$.
    \end{proof}
\end{lemma}

\begin{lemma}\label{lem:nisgluing}
    Let $\Cscr$ be an $R$-linear category.
    Let $f:X'\rightarrow X$ be a morphism in $\Test_{\pi_0 R}$ where $X'$ is affine. Suppose that $\Mod_X(\Cscr)$ is compactly
    generated, and suppose that there exists a quasi-compact
    open subset $U\subseteq X$ with complement $Z$ such that $f|_{Z'}:Z'\rightarrow Z$ is an equivalence, where
    $Z'=Z\times_X X'$, and such that $\Mod_{X'}(\Cscr)$ and $\Mod_U(\Cscr)$ possess compact generators $P$ and $Q$.
    Then, $\Mod_X(\Cscr)$ has a compact generator.
    \begin{proof}
        We verify the conditions of Lemma~\ref{lem:closedgluing}. Because $\Mod(\Cscr)$ is a Nisnevich sheaf, there is a
        cartesian square of $R$-linear $\infty$-categories
        \begin{equation*}
            \begin{CD}
                \Mod_X(\Cscr)   @>>>    \Mod_{X'}(\Cscr)\\
                @VVV                    @VVV\\
                \Mod_U(\Cscr)   @>>>    \Mod_{U'}(\Cscr).
            \end{CD}
        \end{equation*}
        Taking the fibers of the vertical maps induces an equivalence $\Mod_{X,Z}(\Cscr)\we\Mod_{X',Z'}(\Cscr)$. By
        Proposition~\ref{prop:localization}, the fact that $\Mod_{X'}(\Cscr)$ has a compact generator implies that
        $\Mod_{X',Z'}(\Cscr)$ has a compact generator, and hence $\Mod_{X,Z}(\Cscr)$ has a compact generator.
        To finish the proof, we show that $\Mod_U(\Cscr)$ has a compact generator which is the restriction of a compact
        object over $X$. But, by Proposition~\ref{prop:decomposition}, $Q\oplus\Sigma Q$ is
        the restriction of a compact object of $X$. It clearly generates $\Mod_U(\Cscr)$.
    \end{proof}
\end{lemma}

Let $\Cscr$ be an $R$-linear $\infty$-category, and let $\chi_{\Cscr}$ be the presheaf on
$\CAlg_R^{\et}$ defined by
\begin{equation*}
    \chi_{\Cscr}(S)=\begin{cases}
        {*} &   \text{if $\Mod_S(\Cscr)$ has a compact generator,}\\
        \emptyset & \text{otherwise.}
    \end{cases}
\end{equation*}
The presheaf $\chi_{\Cscr}$ extends to a presheaf $\chi'_{\Cscr}$ on $\Test_{\pi_0 R}$ by right Kan extension. By definition, if $X$ is an
object of $\Test_{\pi_0 R}$, then $\chi'_{\Cscr}(X)$ is contractible if and only if $\chi_{\Cscr}(S)$ is
non-empty for all $R$-algebras $S$ and all $\Spec\pi_0 S\rightarrow X$.

\begin{lemma}\label{lem:localglobalfiniteflat}
    Suppose that $\Cscr$ is an $R$-linear $\infty$-category and that $R\rightarrow S$ is a finite
    faithfully flat cover. Then, $\Mod_S(\Cscr)$ has a compact generator if and only if $\Mod_R(\Cscr)$ does.
    \begin{proof}
        If $\Mod_R(\Cscr)$ has a compact generator, then it is a perfect generator by
        Lemma~\ref{lem:compactsonaffines}, so that $\Mod_S(\Cscr)$ has a compact generator.
        Suppose that $\Mod_S(\Cscr)$ has a compact generator $P$. The functor
        $\pi_*:\Mod_S\rightarrow\Mod_R$ has a right adjoint, which is given explicitly by
        $\pi^!(M)=\Map_R(S,M)$. Since $S$ is a finite and flat $R$-module, it follows that
        $\pi^!$ preserves filtered colimits. Therefore,
        $\pi_*:\Mod_S(\Cscr)\rightarrow\Mod_R(\Cscr)$ has a continuous right adjoint given
        by tensoring $\Cscr$ with $\pi^!:\Mod_R\rightarrow\Mod_S$. We abuse notation and
        write $\pi^!$ for this right adjoint as well. It follows immediately that $\pi_*$
        preserves compact objects so that $\pi_*(P)$ is a compact object of $\Mod_R(\Cscr)$.
        To show that $\pi_*(P)$ is a compact generator of $\Mod_R(\Cscr)$, suppose that
        $\Map_R(\pi_*(P),M)\we 0$. Using the adjunction, we get that $\Map_S(P,\pi^!(M))\we
        0$. Therefore, $\pi^!(M)\we 0$. In general, the functor $\pi_*$ is conservative.
        But, $\pi_*\pi^!(M)\we S^{\vee}\otimes_R M$, so that $\pi_*\pi^!$ is conservative by
        the faithful flatness of $S$. Therefore, $\pi^!$ is conservative. Thus, $M\we 0$, so
        that $\pi_*(P)$ is a compact generator of $\Mod_R(\Cscr)$.
    \end{proof}
\end{lemma}

Now, we come to the \'etale local-global principle. The idea of the proof is due to
Lurie~\cite{dag11}*{Section 6}.

\begin{theorem}\label{thm:etalegluing}
    If $\Cscr$ is an $R$-linear $\infty$-category, then $\chi_{\Cscr}$ is an \'etale sheaf.
    \begin{proof}
        By~\cite{dag11}*{Theorems~2.9, 3.7}, it suffices to show that $\chi_{\Cscr}$ satisfies
        finite \'etale descent, and that $\chi_{\Cscr}'$ satisfies affine Nisnevich excision. Finite \'etale
        descent follows from Lemma~\ref{lem:localglobalfiniteflat}. To show that
        $\chi_{\Cscr}'$
        satisfies affine Nisnevich excision, suppose that $f:X'\rightarrow X$ is an affine morphism in
        $\Test_R$, that $U\subseteq X$ is a quasi-compact open subset such that
        $X'-U'\we X-U$, where $U'=X'\times_X U$, and that $\chi'_{\Cscr}(X')$ and $\chi'_{\Cscr}(U)$ are non-empty.
        Note that by~\cite{dag11}*{Proposition 6.12 and Lemma~6.17} all of the stable
        presentable $\infty$-categories that appear in proof are compactly generated. This is
        important because we will use Lemma~\ref{lem:nisgluing}.
        To show that $\chi'_{\Cscr}(X)$ is non-empty, let $\Spec
        S\rightarrow X$ be a point of $X$. Pull back the affine elementary Nisnevich square via this
        map, to obtain
        \begin{equation*}
            \begin{CD}
                X'\times_X U\times_X\Spec S @>>>    X'\times_X\Spec S\\
                @VVV                                @VVV\\
                U\times_X\Spec S            @>>>    \Spec S.
            \end{CD}
        \end{equation*}
        By our hypotheses, $X'\times_X\Spec S$ is affine, so that
        $\chi_{\Cscr}(X'\times_X\Spec S)\we\chi_{\Cscr}'(X'\times_X\Spec S)$
        is contractible, as we see by using the map $X'\times_X\Spec S\rightarrow X'$.
        By Lemma~\ref{lem:nisgluing}, to complete the proof, it suffices to show
        that $\Mod_{U\times_X\Spec S}(\Cscr)$ has a compact
        generator. By hypothesis, we know that $\chi'_{\Cscr}(U\times_X\Spec S)$ is non-empty. As
        $U$ is quasi-compact in $X$, we may write $U\times_X\Spec S$ as a union of Zariski open subschemes.
        \begin{equation*}
            U\times_X\Spec S=\bigcup_{i=1}^n\Spec S_i.
        \end{equation*}
        Since $\chi'_{\Cscr}(U\times_X\Spec S)$ is non-empty, $\Mod_{\Spec S_i}(\Cscr)$ has a
        compact generator for all $i$. Write
        \begin{equation*}
            V_k=\bigcup_{i=1}^k\Spec S_i,
        \end{equation*}
        and assume that $\Mod_{V_k}(\Cscr)$ has a compact generator for some $k$ in $[1,n)$. Then,
        $\Spec S_{k+1}\rightarrow V_{k+1}$ and the open $V_k\subseteq V_{k+1}$ satisfy the
        hypotheses of Lemma~\ref{lem:nisgluing} (take $X=V_{k+1}$, $X'=\Spec S_{k+1}$, and $U=V_k$).
        Therefore, $\Mod_{V_{k+1}}(\Cscr)$ has a compact generator. By induction, we see that
        $\Mod_{U\times_X\Spec S}(\Cscr)$ has a compact generator, as desired.
    \end{proof}
\end{theorem}

\subsection{Lifting theorems}

Now we put together the local-global principles of the previous sections into one of the
main theorems of the paper. In the case of schemes built from simplicial commutative rings, this was proved
in~\cite{toen-derived}*{Theorem 4.7}. Our proof is rather different, as the \'etale
local-global principle requires different methods for connective $\EE_\infty$-rings.

\begin{theorem}\label{thm:gluing}
    Let $X$ be a quasi-compact, quasi-separated derived scheme. Then, every Morita
    class $\alpha:X\rightarrow\ShMr$ on $X$ lifts to an algebra $X\rightarrow\ShAlg$.
    \begin{proof}
        By definition of sheafification, the Morita class $\alpha:X\rightarrow\ShMr$ lifts
        \'etale locally through $\ShAlg\rightarrow\ShP$. It follows that there is an
        \'etale cover $\pi:\coprod_i\Spec T_i\rightarrow X$ such that
        $\pi^*\alpha:\coprod_i\Spec T_i\rightarrow\ShMr$ factors
        through $\ShAlg\rightarrow\ShP$; in other words, $\Mod_{T_i}^{\alpha}$ has a compact
        generator for all $i$. Since \'etale maps are open, we can assume,
        possibly by refining the cover, that the image of $\Spec T_i$ in $X$ is an affine
        subscheme $\Spec S_i$. By Theorem~\ref{thm:etalegluing}, it follows that
        $\Mod_{S_i}^{\alpha}$ has a compact generator. By Theorem~\ref{thm:zariskigluing},
        it follows that $\Mod_X^{\alpha}$ has a perfect generator. This completes the proof
        by Proposition~\ref{prop:factorization}.
    \end{proof}
\end{theorem}

We now consider several applications, which show the power of this theorem in establishing the
compact generation of various stable presentable $\infty$-categories. These are motivated by
the results of~\cite{schwede-shipley} in the affine case.

\begin{example}
    If $X$ is a quasi-compact and quasi-separated derived scheme, and if $E$ is an
    $R$-module such that localization with respect to $E$ is smashing, then $L_E\Mod_X$, the
    full subcategory of $E$-local objects in $\Mod_X$,
    is compactly generated by a single compact object.
\end{example}

\begin{example}
    If $X$ is a quasi-compact and quasi-separated derived scheme over the $p$-local sphere,
    consider the localization $L_{K(n)}\Mod_X$, where $K(n)$ is the $n$th Morava
    $K$-theory at the prime $p$. In this case, the $K(n)$-localization of $\Oscr_X$ need not
    be compact in $L_{K(n)}\Mod_X$. However, if $F$ is a finite type $n$ complex, then over any affine $\Spec S\rightarrow X$, the
    $K(n)$-localization of $S\otimes F$ \emph{is} a compact generator of
    $L_{K(n)}\Mod_S$. It follows from Theorem~\ref{thm:gluing} that there is a compact
    generator of $L_{K(n)}\Mod_X$.
\end{example}

Our main application of the theorem is the following statement.

\begin{corollary}\label{cor:azumayagluing}
    Let $X$ be a quasi-compact, quasi-separated derived scheme. Then, every Brauer
    class $\alpha:X\rightarrow\ShBr$ on $X$ lifts to an Azumaya algebra $X\rightarrow\ShAz$.
\end{corollary}

Note that this theorem is false in non-derived algebraic geometry. There is a non-separated,
but quasi-compact and quasi-separated, surface $X$ and a non-zero cohomological Brauer class
$\alpha\in\Hoh^2_{\et}(X,\Gm)_{\tors}$ that is not represented by an ordinary Azumaya
algebra~\cite{ehkv}*{Corollary 3.11}. In this case,
the Brauer class vanishes on a Zariski cover of $X$. However, there is no global
$\alpha$-twisted vector bundle, so there cannot be a non-derived Azumaya algebra. The
corollary shows that, even in this case, there \emph{is} a derived Azumaya algebra with
class $\alpha$.

\section{Brauer groups}\label{sec:brauer}

We prove our main theorems on the Brauer group, which will, in particular, allow
us to show that the Brauer group of the sphere spectrum vanishes.

\subsection{The Brauer space}

Classically, there are two Brauer groups of a commutative ring or a scheme $X$.
The first is the algebraic Brauer group, which is the group of Morita equivalence
classes of Azumaya algebras over $X$. This notion goes back to Azumaya~\cite{azumaya} for
algebras free over commutative rings, Auslander
and Goldman~\cite{auslander-goldman} for the general affine case, and
Grothendieck~\cite{grothendieck-brauer-1} for schemes. The second is the cohomological
Brauer group $\Hoh^2_{\et}(X,\Gm)_{\tors}$ introduced by Grothendieck~\cite{grothendieck-brauer-1}.
There is an inclusion from the algebraic Brauer group into the cohomological Brauer
group (under the reasonable assumption that $X$ has only finitely many connected components),
but they are not always identical, as noted above. As a result of Theorem~\ref{cor:azumayagluing}, the natural
generalizations of these definitions to quasi-compact, quasi-separated schemes \emph{do} agree.
Moreover, these generalizations yield not just groups but in fact grouplike
$\Einfinity$-spaces; the Brauer groups are the groups of connected components of these
spaces. We work again over some fixed connective $\Einfinity$-ring $R$.

\begin{definition}
    Let $X$ be an \'etale sheaf. Then, the Brauer space of $X$ is $\ShBr(X)$, the space of
    maps from $X$ to $\ShBr$ in $\Shv_R^{\et}$. The Brauer group of $X$ is $\pi_0\ShBr(X)$.
\end{definition}

When $X$ is an arbitrary \'etale sheaf, we cannot say too much about the algebraic nature of
the classes in $\pi_0\ShBr(X)$. However, write $\ShBr_{\alg}(X)$ for the full subspace of
$\ShBr(X)$ of classes $\alpha:X\rightarrow\ShBr$ that factor through
$\ShAz\rightarrow\ShBr$. In other words, $\pi_0\ShBr_{\alg}(X)$ is the subgroup of the
Brauer group consisting of those classes representable by an Azumaya algebra over $X$.
When $X=\Spec S$, we will write $\ShBr(S)$ for $\ShBr(\Spec S)$.

We now can answer the analogue of the $\Br=\Br'$ question of Grothendieck.

\begin{theorem}
    For any quasi-compact and quasi-separated derived scheme $X$, 
    $\ShBr_{\alg}(X)\we\ShBr(X)$.
    \begin{proof}
        This is the content of Corollary~\ref{cor:azumayagluing}.
    \end{proof}
\end{theorem}

An important fact about the Brauer space of a connective commutative ring spectrum is that
it has a purely categorical formulation.
Recall that $\Cat_{S,\omega}$ is the symmetric
monoidal $\infty$-category of compactly generated $S$-linear categories together with
colimit preserving functors that preserve compact objects. We saw in
Theorem~\ref{thm:azumayaiffinvertible} that if $A$ is an
$S$-algebra, then $\Mod_A$ is invertible in $\Cat_{S,\omega}$ if and only if $A$ is Azumaya.
Write $\Cat_{S,\omega}^{\times}$ for the grouplike $\Einfinity$-space of invertible objects
in $\Cat_{S,\omega}$.

\begin{proposition}
    If $S$ is a connective commutative $R$-algebra, then the natural morphism
    $\Cat_{S,\omega}^{\times}\rightarrow\ShBr(S)$ is an equivalence.
    \begin{proof}
        Consider the diagram
        \begin{equation*}
            \Cat_{S,\omega}^\times\xrightarrow{i}\ShBr(S)\xrightarrow{j}\ShP(S).
        \end{equation*}
        The composition $j\circ i$ is fully faithful, by definition. The map $j$ is fully
        faithful by construction of $\ShBr$. Thus, $i$ is fully faithful. On the other hand, by
        Corollary~\ref{cor:azumayagluing}, the map $i$ is essentially surjective. Thus, $i$
        is an equivalence.
    \end{proof}
\end{proposition}

This proposition has the following two interesting corollaries, which will not be used in
the sequel.

\begin{corollary}
    The presheaf of spaces that sends a connective commutative $R$-algebra $S$ to $\Cat_{S,\omega}^{\times}$ is an \'etale sheaf.
\end{corollary}

\begin{corollary}
    The space $\ShBr(X)$ is a grouplike $\Einfinity$-space.
    \begin{proof}
        The space $\ShBr(S)$ is a grouplike $\Einfinity$-space for every connective
        commutative $R$-algebra $S$, and the grouplike $\Einfinity$-structure is natural
        in $S$. Thus, $\ShBr$ is a grouplike $\Einfinity$-object in $\Shv_R^{\et}$. The
        mapping space
        \begin{equation*}
            \ShBr(X)=\Map_{\Shv_R^{\et}}(X,\ShBr)
        \end{equation*}
        thus inherits a grouplike $\Einfinity$-structure from that on $\ShBr$.
    \end{proof}
\end{corollary}

As a result of the corollary, when $X$ is an \'etale sheaf, we may construct via delooping a spectrum $\Shbr(X)$, with
$\Omega^{\infty}\Shbr(X)\we\ShBr(X)$. A similar idea has been pursued recently by
Szymik~\cite{szymik}, but with rather different methods.

We will need the following proposition, as well as the computations in the following
section, to tell us the homotopy sheaves of $\ShBr$. This will be used to give a complete computation
of $\ShBr(X)$ using a descent spectral sequence when $X$ is affine.

\begin{proposition}
    There is a natural equivalence of \'etale sheaves $\Omega\ShBr\we\ShPic$, where $\ShPic$
    is the sheaf of line bundles.
    \begin{proof}
        By the \'etale local triviality of Azumaya algebras proven in Theorem~\ref{thm:etalelocaltriviality},
        it follows that $\ShBr$ is a connected sheaf and that
        it is equivalent to the classifying space of the trivial Brauer class.
        But, the sheaf of auto-equivalences of $\StMod$ is precisely the sheaf of line bundles in
        $\StMod$.
    \end{proof}
\end{proposition}

\subsection{Picard groups of connective ring spectra}

In the previous section, we showed that $\Omega\ShBr\we\ShPic$, and by the \'etale local
triviality of Azumaya algebras, we know that the sheaf $\pi_0\ShBr$ vanishes. Thus, to compute the
homotopy sheaves of $\ShBr$, it is enough to compute them for $\ShPic$, which is what we now
do.

If $R$ is a discrete commutative ring, let $\Pic(R)$ be the Picard group of invertible
$R$-modules. This should be distinguished from $\ShPic(\Hrm R)$, the grouplike
$\Einfinity$-space of invertible $\Hrm R$-modules, and from $\Pic(\Hrm R)$.

\begin{proposition}[Fausk \cite{fausk}]\label{prop:hr}
    Let $R$ be a discrete commutative ring. Then, there is an exact sequence
    \begin{equation*}
        0\rightarrow\Pic(R)\rightarrow\pi_0\ShPic(\Hrm R)\xrightarrow{c}\Hoh^0(\Spec R,\ZZ)\rightarrow 0,
    \end{equation*}
    where the inclusion comes from the monoidal functor $\Mod_{R}\rightarrow\Mod_{HR}$,
    and the map $c$ sends an invertible element $L$ to its degree of connectivity on each
    connected component of $\Spec R$. Thus, $c(L)=n$ if and only
    if $\pi_m(L)=0$ for $m<n$, and $\pi_n(L)\neq 0$.
\end{proposition}

The purpose of this section is to extend Proposition~\ref{prop:hr} to all connective commutative rings.
The following lemma is essentially found in~\cite{hopkins-mahowald-sadofsky}*{Page 90}.
We remark that if $L$ is an invertible $R$-module, then $L$ is perfect and $L^{-1}$ is the dual of $L$,
$\Map_R(L,R)$. It follows that there is a canonical evaluation map $\ev:L\otimes_R
L^{-1}\rightarrow R$, which is an equivalence.

\begin{lemma}\label{lem:equivalences}
    Let $R$ be an $\Einfinity$-ring spectrum, and let $L$ be an invertible $R$-module.
    Suppose that there are $R$-module maps $\phi:\Sigma^nR\rightarrow L$ and $\omega:\Sigma^{-n}R\rightarrow L^{-1}$
    such that
    \begin{equation*}
        \ev\circ\phi\otimes_R\omega:R\we\Sigma^nR\otimes_R\Sigma^{-n}R\rightarrow L\otimes_R
        L^{-1}\rightarrow R
    \end{equation*}
    is homotopic to the identity. Then, $\phi$ and $\omega$ are weak equivalences.
    \begin{proof}
        The $n$th suspension of $\ev\circ\phi\otimes_R\omega$ is homotopic to the composition
        \begin{equation*}
            \Sigma^nR\xrightarrow{\phi}L\we
            L\otimes_{R}R\xrightarrow{1\otimes\Sigma^{n}\omega}L\otimes_{R}\Sigma^nL^{-1}\rightarrow\Sigma^n R.
        \end{equation*}
        Therefore, $\Sigma^n R$ is a retract of $L$; specifically, there exists a perfect
        $R$-module $M$ and an equivalence $L\we\Sigma^nR\oplus M$.
        Similarly, $L^{-1}\we\Sigma^{-n}\oplus N$ for some perfect $R$-module $N$. But,
        \begin{equation*}
            R\we L\otimes_R L^{-1}\we(\Sigma^nR\oplus M)\otimes_R(\Sigma^{-n}R\oplus N)\we
            R\oplus \Sigma^{-n}M\oplus \Sigma^n N\oplus \left(M\otimes_R N\right),
        \end{equation*}
        which shows that $M$ and $N$ are zero, and hence that $\phi$ and $\omega$ are
        equivalences.
    \end{proof}
\end{lemma}

\begin{theorem}\label{thm:rlocal}
    Let $R$ be a connective local $\Einfinity$-ring spectrum (that is, $\pi_0 R$ is a local ring).
    Then, $R\rightarrow \tau_{\leq 0}R\we\Hrm\pi_0 R$ induces an isomorphism
    $\pi_0\ShPic(R)\rightarrow\pi_0\ShPic(\tau_{\leq 0}R)\iso\ZZ$.
    \begin{proof}
        Since $\pi_0 R$ is local, $\pi_0\ShPic(\tau_{\leq 0}R)=\ZZ$ by
        Proposition~\ref{prop:hr}. Thus, it suffices to show that if $L$ is an invertible $R$-module, then $L\we\Sigma^n R$ for some $n$. Fixing $L$,
        we first identify the appropriate integer $n$.

        The invertibility of $L$ implies that $L$ is a perfect $R$-module.
        By Proposition~\ref{prop:compactness}, it follows that $L$ has a bottom homotopy group, say $\pi_n L$.
        This means that for $m<n$, $\pi_m L=0$, while $\pi_n L\neq 0$. Similarly, let $\pi_m L^{-1}$
        be the bottom homotopy group of $L^{-1}$.   We  will  show  that  $n=-m$,  and  that
        $L\we\Sigma^n R$.  Consider the $\Tor$ spectral  sequence  for  $L\otimes_R  L^{-1}$:
        \begin{equation*}
            \Eoh^2_{p,q}=\Tor_{p}^{\pi_* R}(\pi_* L,\pi_* L^{-1})_q\Rightarrow \pi_{p+q}R.
        \end{equation*}
        The differential $d^r$ is of degree $(-r,r-1)$. Thus, for degree reasons, $\Eoh^2_{0,n+m}=\Eoh^{\infty}_{0,n+m}$. In this case, we have
        \begin{equation*}
            (\pi_* L\otimes_{\pi_* R} \pi_* L^{-1})_{n+m}\iso \pi_n L\otimes_{\pi_0 R}
            \pi_m L^{-1}.
        \end{equation*}
        Since $\pi_n L$ and $\pi_m L^{-1}$ are non-zero and $\pi_0 R$ is local, the term $\Eoh^2_{0,n+m}$ is non-zero.
        It is the term of smallest total degree that is non-zero. Thus, since it is
        permanent in the spectral sequence,
        \begin{equation*}
            \pi_n L\otimes_{\pi_0 R} \pi_m L^{-1}\iso \pi_0 R,
        \end{equation*}
        and $n=-m$. Again, since $\pi_0 R$ is local, $\pi_n L$ and $\pi_m L^{-1}$ are both
        in fact isomorphic to $\pi_0 R$.

        Choose $\phi\in\pi_n L$ and $\omega\in\pi_m L^{-1}$ such that the isomorphism above
        gives $\phi\otimes_R\omega=1_R\in \pi_0 R$.
        The homotopy classes $\phi$ and $\omega$ are represented by $R$-module maps
        \begin{align*}
            \phi:\Sigma^n R &\rightarrow L,\\
            \omega:\Sigma^m R   &\rightarrow L^{-1}.
        \end{align*}
        Then,
        \begin{equation*}
            R\xrightarrow{\phi\otimes_R\omega}L\otimes_R L^{-1}\rightarrow R
        \end{equation*}
        is homotopic to $\phi\otimes_R\sigma\we 1_R$. Thus, applying
        Lemma~\ref{lem:equivalences}, the $R$-module maps $\phi$ and $\omega$ are in fact equivalences.
        This completes the proof.
    \end{proof}
\end{theorem}

Consider the \'etale sheaf $\ShGL_1$, which sends a connective commutative $R$-algebra $S$ to
the space of units in $S$. That is, $\ShGL_1(S)$ is defined as the pullback in the diagram
of spaces
\begin{equation*}
    \xymatrix{
        \ShGL_1(S)    \ar[r]\ar[d] &    \Omega^{\infty}S\ar[d]\\
        \pi_0S^\times   \ar[r] &    \pi_0 S.
    }
\end{equation*}
The classifying space $\Brm\ShGL_1(S)$ of this grouplike $\Einfinity$-space naturally includes
as the identity component into $\ShPic(S)$. Thus, there is a natural map
$\ShB\ShGL_1\rightarrow\ShPic$ from the classifying sheaf of $\ShGL_1$ into $\ShPic$. When $S$
is a local connective commutative $R$-algebra, then $\ShPic(S)$ decomposes as the product
$\ShB\ShGL_1(S)\times\ZZ$, where the map
\[
\ZZ\longrightarrow\ShPic(S)
\]
sends $n$ to $\Sigma^n S$.
Thus, we have the following corollary.

\begin{corollary}
    The sequence $\ShB\ShGL_1\rightarrow\ShPic\rightarrow\ZZ$ is a split fiber sequence of hypersheaves.
    \begin{proof}
        Since $\ShGL_1$ is a hypersheaf,
        so is $\ShB\ShGL_1$. We also know that $\ShPic$ is a hypersheaf by
        Proposition~\ref{cor:unitsaresheaves}. Finally, $\ZZ$ is by definition the hypersheaf
        associated to the constant presheaf with values $\ZZ$. Evidently, the sequence is in
        fact a sequence of sheaves of grouplike $\EE_1$-spaces. Since $\ZZ$ is freely generated as a
        sheaf of grouplike $\EE_1$-spaces by a single object, the splitting is obtained by
        taking the canonical basepoint of $\ShPic$.
    \end{proof}
\end{corollary}

With this corollary, we can give the computation of the homotopy sheaves of $\ShBr$, which
we need in the next section in order to actually compute the Brauer groups of some
connective $\Einfinity$-rings.

\begin{corollary}
    The homotopy sheaves of $\ShBr$ are
    \begin{equation}\label{eq:pistarbr}
        \pi_i\ShBr\iso\begin{cases}
            0   &   \text{if $i=0$,}\\
            \ZZ &   \text{if $i=1$,}\\
            \pi_0\Oscr^* &   \text{if $i=2$,}\\
            \pi_{i-2}\Oscr  &   \text{if $i\geq 3$},
        \end{cases}
    \end{equation}
    where $\Oscr$ is the structure sheaf on $\Shv_R^{\et}$.
\end{corollary}

\subsection{The exact sequence of Picard and Brauer groups}

Suppose that $X=U\cup V$ is a derived scheme, written as the union of two Zariski open
subschemes. Then, because $\ShBr$ is an
\'etale sheaf, there is a fiber sequence of spaces
\begin{equation*}
    \ShBr(X)\rightarrow\ShBr(U)\times\ShBr(V)\rightarrow\ShBr(U\cap V).
\end{equation*}
Taking long exact sequences, we obtain the following exact sequence:
\begin{gather*}
    \pi_2\ShBr(U\cap
    V)\rightarrow\pi_1\ShBr(X)\rightarrow\pi_1\ShBr(U)\oplus\pi_1\ShBr(V)\rightarrow
    \pi_1\ShBr(U\cap V)\\
    \rightarrow\pi_0\ShBr(X)\rightarrow\pi_0\ShBr(U)\oplus\pi_0\ShBr(V)\rightarrow\pi_0\ShBr(U\cap V),
\end{gather*}
which generalizes the classical Picard-Brauer exact sequence
\begin{equation*}
    \Pic(X)\rightarrow\Pic(U)\oplus\Pic(V)\rightarrow\Pic(U\cap
    V)\rightarrow\Br(X)\rightarrow\Br(U)\oplus\Br(V)\rightarrow\Br(U\cap V),
\end{equation*}
when $U$, $V$, and $X$ are ordinary schemes.
The computations in the next section can be used to show that the sequence is not, in general, exact on the right.

The important connecting morphism $\delta:\pi_1\ShBr(U\cap V)\rightarrow\pi_0\ShBr(X)$ can be described in the
following Morita-theoretic way. The $\infty$-category $\Mod_X$ of quasi-coherent sheaves on
$X$ can be glued from $\Mod_U$ and $\Mod_V$ by taking the natural equivalence
$\Mod_U|_{U\cap V}\we\Mod_V|_{U\cap V}$. On the other hand, given a line bundle $L$ over
$U\cap V$, we can twist the gluing data by tensoring with $L$. The resulting category is
$\Mod_X^{\delta(L)}$, the $\infty$-category of quasi-coherent $\delta(L)$-twisted sheaves.

\subsection{The Brauer space spectral sequence}\label{sec:brss}

In this section, we obtain a spectral sequence
converging conditionally to the homotopy groups of $\ShBr(X)$. In most cases of interest, for instance
when $X$ is affine or has finite \'etale cohomological dimension, we show that the spectral sequence
converges completely (see~\cite{bousfield-kan}*{Section IX.5}). In particular, the graded
pieces of the filtration on the abuttment of the spectral sequence are in fact computed by the spectral sequence.
As an application, in the next section, we give various example computations of Brauer
groups. For now, we fix a connective $\Einfinity$-ring spectrum $R$.

If $A$ is a grouplike $\Einfinity$-object of $\Shv_R^{\et}$, and if $X$ is any object of
$\Shv_R^{\et}$, then for every $p\geq 0$, there is a cohomology group
\begin{equation*}
    \Hoh^p_{\et}(X,A)=\pi_0\Map_{\Shv_R^{\et}}(X,\ShB^p A),
\end{equation*}
where $\ShB^p A$ denotes a $p$-fold delooping of $A$. In particular, if $A$ is a sheaf of abelian
groups in $\Shv_X^{\et}$, then we can view $A$ canonically as a grouplike
$\Einfinity$-space. An $\infty$-topos $\Xscr$ has cohomological dimension $\leq n$ if
$\Hoh^m(\Xscr,A)=0$ for all abelian sheaves $A$ in $\Xscr$ and all
$m>n$~\cite{htt}*{Definition~7.2.2.18}.

Recall that by~\cite{ha}*{Theorem~8.5.0.6},
the small \'etale site on $\Spec S$ is equivalent to the nerve of the small
\'etale site on $\Spec\pi_0 S$. Therefore, by~\cite{htt}*{Remark~7.2.2.17}, for any sheaf of
abelian groups $A$ over $S$, there is a natural isomorphism
\begin{equation*}
    \Hoh^p_{\et}(\Spec S,A)\iso\Hoh^p_{\et}(\Spec\pi_0S,A),
\end{equation*}
where the right-hand side denotes the classical \'etale cohomology groups over $\Spec\pi_0 S$.

\begin{theorem}
    Let $X$ be an object of $\Shv_R^{\et}$. Then, there is a conditionally convergent spectral sequence
    \begin{equation}\label{eq:brauerss}
        \Eoh_2^{p,q}=\begin{cases}
            \Hoh^p_{\et}(X,\pi_q\ShBr) & p\leq q\\
            0                           & p>q
        \end{cases}
        \,\,\,\Rightarrow\,\,\,\pi_{q-p}\ShBr(X)
    \end{equation}
    with differentials $d_r$ of degree $(r,r-1)$. If $X$ is affine, discrete, or if $\left(\Shv_R^{\et}\right)_{/X}$ has finite
    cohomological dimension, then  the  spectral  sequence  converges  completely. 
    \begin{proof}
        Because $\ShBr$ is hypercomplete, the map from $\ShBr$ to the limit
        of its Postnikov tower $\ShBr\rightarrow\lim_n\tau_{\leq n}\ShBr$ is an equivalence. See~\cite{htt}*{Section 6.5}.
        Taking sections preserves limits, so that
        $$\ShBr(X)\rightarrow\lim_n\left(\left(\tau_{\leq n}\ShBr\right)(X)\right)$$ is also an equivalence. Thus,
        $\ShBr(X)$ is the limit of a tower, and to any such tower
        there is an associated spectral sequence \cite{bousfield-kan}*{Chapter IX} which converges conditionally to the homotopy groups of the limit. Using the methods
        of Brown and Gersten~\cite{brown-gersten}, one identifies the $\Eoh_2$-page
        as~\eqref{eq:brauerss}.

        If $X$ is affine, discrete, or if $\left(\Shv_R^{\et}\right)_{/X}$ has finite cohomological
        dimension, then the spectral sequence degenerates at some finite page. This is
        clear in the latter case, and if $X$ is discrete the spectral sequence collapses
        entirely at the $\Eoh_2$-page. So, suppose that $X=\Spec S$. Then, $\ShBr(X)$
        can be computed on the small \'etale site on $\Spec S$. But, as mentioned above, this site is the nerve
        of a discrete category, the small \'etale site on $\Spec\pi_0 S$. Therefore,
        \begin{equation*}
            \Hoh^p_{\et}(\Spec S,\pi_q\ShBr)\iso\Hoh^p_{\et}(\Spec\pi_0 S,\pi_q\ShBr).
        \end{equation*}
        Since $\pi_q\ShBr\we\pi_{q-2}\Oscr$ for $q\geq 3$, and since these are all
        quasi-coherent $\pi_0\Oscr$-modules, it follows that
        \begin{equation*}
            \Hoh^p_{\et}(\Spec S,\pi_q\ShBr)\iso\Hoh^p_{\et}(\Spec\pi_0 S,\pi_{q-2}\Oscr)=0
        \end{equation*}
        for $q\geq 3$ and $p\geq 1$ by Grothendieck's vanishing theorem.
        Thus, the only possible differentials are
        \begin{equation*}
            d_2:\Hoh^p(\Spec S,\ZZ)\rightarrow\Hoh^{p+2}(\Spec S,\pi_0\Oscr^{\times}).
        \end{equation*}
        However, these differentials vanish because $\Brm\ZZ$ is in fact a split retract of
        $\ShBr$. Therefore, if $X$ is affine, the spectral sequences degenerates at the $\Eoh_2$ page.
        It follows from the degeneration and the complete convergence
        lemma~\cite{bousfield-kan}*{IX.5.4} that the
        spectral sequence converges completely to $\pi_*\ShBr(X)$. This completes the proof.
    \end{proof}
\end{theorem}

Using the theorem and the remarks preceeding it, we deduce the following corollary, which
completely computes the homotopy groups of the Brauer space of a connective commutative ring
$R$.

\begin{corollary}\label{cor:brcomputation}
    If $R$ is a connective $\Einfinity$-ring spectrum, then the homotopy groups of
    $\ShBr(R)$ are described by
    \begin{equation*}
        \pi_k\ShBr(R)\iso\begin{cases}
            \Hoh^1_{\et}(\Spec\pi_0 R,\ZZ)\times\Hoh^2_{\et}(\Spec\pi_0 R,\Gm) & k=0\\
            \Hoh^0_{\et}(\Spec\pi_0 R,\ZZ)\times\Hoh^1_{\et}(\Spec\pi_0 R,\Gm) & k=1\\
            \pi_0 R^{\times}    & k=2\\
            \pi_{k-2}R & k\geq 3.
        \end{cases}
    \end{equation*}
    \begin{proof}
        This follows immediately from the degeneration of the Brauer spectral sequence for $\Spec R$
        together with the fact that $\ShB\ZZ$ splits off of $\ShBr$.
    \end{proof}
\end{corollary}

Note that in the special case where $R$ is a discrete commutative ring, Szymik obtained similar
computations for the purely algebraic Brauer spectrum of $\Hrm R$ defined in~\cite{szymik}.
The computations also follow from the next corollary.

\begin{corollary}\label{cor:brclassical}
    If $X$ is a quasi-compact and quasi-separated ordinary scheme, then
    \begin{equation*}
         \pi_k\ShBr(X)\iso\begin{cases}
            \Hoh^1_{\et}(X,\ZZ)\times\Hoh^2_{\et}(X,\Gm) & k=0\\
            \Hoh^0_{\et}(X,\ZZ)\times\Hoh^1_{\et}(X,\Gm) & k=1\\
            \Hoh^0_{\et}(X,\Gm)    & k=2\\
            0 & k\geq 3.
        \end{cases}
    \end{equation*}
\end{corollary}

\subsection{Computations of Brauer groups of ring spectra}

In this section, we give several examples of Brauer groups of ring spectra and of derived
schemes. Our convention throughout this section is to write $\Br(R)$ for the Brauer group of
Azumaya algebras over a discrete commutative ring $R$. This injects but is not, in general, the same as
$\pi_0\ShBr(\Hrm R)$, as we will see below. Note that $\Br(R)\iso\Hoh^2_{\et}(\Spec
R,\Gm)_{\tors}$, by Gabber~\cite{gabber}. If $R$ is a regular domain, then
by~\cite{grothendieck-brauer-2}*{Corollaire 1.8}, we have
$\Hoh^2_{\et}(\Spec R,\Gm)_{\tors}=\Hoh^2_{\et}(\Spec R,\Gm)$.

\begin{lemma}
    If $X$ is a normal ordinary scheme, then $\Hoh^1_{\et}(X,\ZZ)=0$.
    \begin{proof}
        Using the exact sequence $0\rightarrow\ZZ\rightarrow\QQ\rightarrow\QQ/\ZZ\rightarrow
        0$, it is enough to show that $\Hoh^1_{\et}(X,\QQ)=0$. This is can be shown as in~\cite{deninger}*{2.1}.
    \end{proof}
\end{lemma}

However, the $\Hoh^1_{\et}(X,\ZZ)$ term does not always vanish, even when $X$ is ordinary
and affine, so there are some truly exotic elements in the derived Brauer group, even over discrete rings.
Here is an example: let $k$ be an algebraically closed field, and let $R=k[x,y]/(y^2-x^3+x^2)$. Then, $\Spec R$ is
a non-normal affine curve with singular point at $(0,0)$. The normalization of $\Spec R$
is $\AA^1_k$. It follows from~\cite{demeyer}*{Page~19} that $\Br(R)=0$. It is also
known that $\Hoh^1_{\et}(R,\ZZ)\iso\ZZ$. Therefore, we have computed that
$\pi_0\ShBr(\Hrm R)\iso\ZZ$.\footnote{We thank Angelo Vistoli for pointing out this example to us
at~\texttt{mathoverflow.net/questions/84414}.}

We can show that the Brauer group vanishes in many cases.

\begin{theorem}
    Let $R$ be a connective commutative ring spectrum such that $\pi_0 R$ is either $\ZZ$ or
    the ring of Witt vectors $\WW_q$ of $\FF_q$. Then,
    \begin{equation*}
        \pi_0\ShBr(R)=0.
    \end{equation*}
    \begin{proof}
        Both $\ZZ$ and $\WW_p$ are normal, so that $\Hoh^1_{\et}(\pi_0R,\ZZ)=0$. The ring of
        Witt vectors $\WW_q$ is a hensel local ring with residue field $\FF_q$. Thus, by a
        theorem of Azumaya (see~\cite{grothendieck-brauer-1}*{Th\'eor\`eme 1}), there is an
        isomorphism $\Br(\WW_q)\iso\Br(\FF_q)$. But, $\Br(\FF_q)=0$ by a theorem of
        Wedderburn. The Albert-Brauer-Hasse-Noether theorem from class field theory implies that
        $\Hoh^2_{\et}(\Spec\ZZ,\Gm)=0$~\cite{grothendieck-brauer-3}*{Proposition~2.4}.
        Thus, in both cases, we have established the required vanishing.
    \end{proof}
\end{theorem}

\begin{corollary}
    The Brauer group of the sphere spectrum is zero.
\end{corollary}

Of course, it would be nice to have some more examples where the Brauer group does not vanish.
We can give several. First, we recall some standard results, all of which can be found
in~\cite{grothendieck-brauer-3}*{Section 2}. There is a
residue isomorphism
\begin{equation*}
    h_p:\Br(\QQ_p)\rightarrow\Hoh^1_{\et}(\Spec\FF_p,\QQ/\ZZ)\iso\QQ/\ZZ,
\end{equation*}
and, for any open subscheme $U$ of $\Spec\ZZ$, there is an exact sequence
\begin{equation*}
    0\rightarrow\Br(U)\rightarrow\Br(\QQ)\rightarrow\bigoplus_{p\in
    U}\Br(\QQ_p),
\end{equation*}
where the sum is over all prime integers $p$ in $U$. We may also identify
$h_{\RR}:\Br(\RR)\iso\ZZ/2\subseteq\QQ/\ZZ$; the unique non-zero class is represented by the
real quaternions. Finally, there is an exact sequence
\begin{equation*}
    0\rightarrow\Br(\QQ)\rightarrow\Br(\RR)\oplus\bigoplus_p\Br(\QQ_p)\rightarrow\QQ/\ZZ\rightarrow 0,
\end{equation*}
where the right-hand map is induced by mapping $\Br(\RR)$ or $\Br(\QQ_p)$ to $\QQ/\ZZ$ and
summing. These two exact sequences are compatible in the obvious way.

If $\alpha\in\Br(\QQ)$ write $\alpha_p$ for the image of $\alpha$ in $\Br(\QQ_p)$, and write
$\alpha_\RR$ for the image of $\alpha$ in $\Br(\RR)$.
By examining the two exact sequences above, it follows that
\begin{equation*}
    \Br\left(\ZZ\left[\frac{1}{p}\right]\right)\iso\ZZ/2.
\end{equation*}
Indeed, if $\alpha$ is a class of $\Br(\QQ)$ that lifts to $\Br(\ZZ[1/p])$, then it follows
that $h_q(\alpha_q)=0$ for all primes $q\neq p$. Therefore,
$h_p(\alpha_p)+h_\RR(\alpha_\RR)=0$. Since there is a unique non-zero class in $\Br(\RR)$,
the result follows.

Similarly, if $\alpha\in\Br(\QQ)$ lifts to $\Br(\ZZ_{(p)})$, then $h_p(\alpha_p)=0$. Thus, there is an exact sequence
\begin{equation*}
    0\rightarrow\Br(\ZZ_{(p)})\rightarrow \ZZ/2\oplus\bigoplus_{q\neq p}\QQ/\ZZ\rightarrow\QQ/\ZZ\rightarrow 0.
\end{equation*}
We have therefore proven the following corollary to Corollary~\ref{cor:brcomputation}.

\begin{corollary}
    \begin{enumerate}
        \item The Brauer group of the sphere with $p$ inverted is $\pi_0\ShBr(\SS[1/p])\iso\ZZ/2$.
        \item the Brauer group of the $p$-local sphere fits into the exact sequence
            \begin{equation*}
                0\rightarrow\pi_0\ShBr(\SS_{(p)})\rightarrow\ZZ/2\oplus\bigoplus_p\QQ/\ZZ\rightarrow\QQ/\ZZ\rightarrow 0.
            \end{equation*}
        \item There is an isomorphism $\pi_0\ShBr(L_{\QQ_p}\SS)\iso\QQ/\ZZ$, where $L_{\QQ_p}\SS$ is the rational $p$-adic sphere.
    \end{enumerate}
\end{corollary}

Note the important fact that the first two cases in the corollary give examples of
non-Eilenberg-MacLane commutative ring spectra with non-zero Brauer groups.

Finally, we mention two examples of ordinary schemes, where the derived Brauer group exhibits
different behavior than the classical Brauer group. The first is the scheme $X$ used
in~\cite{ehkv}*{Corollary 3.11}, which is the gluing of two affine quadric cones along the
non-singular locus, viewed as a derived scheme over the complex numbers.
This is a normal, qauasi-compact, non-separated, quasi-separated scheme, so it satisfies
the hypotheses of the theorems. One can check that $\pi_0\ShBr(X)=\ZZ/2$ by
Corollary~\ref{cor:brclassical}. This example was
studied originally because the classical Brauer group of the scheme $X$ viewed as an
ordinary geometric object over $\CC$ is $\Br(X)=0$, while the cohomological Brauer group is
$\Br'(X)=\Hoh^2_{\et}(X,\Gm)=\ZZ/2$. In other words, the non-zero class $\alpha\in\Br'(X)$
is represented by an Azumaya algebra, but not by an ordinary Azumaya algebra (an algebra
concentrated in degree $0$).

The second example is the surface of Mumford~\cite{grothendieck-brauer-2}*{Remarques 1.11(b)}.
He constructs a normal surface $Y$ such that
$\Hoh^2_{\et}(Y,\Gm)$ has non-torsion elements. Of course, these can never be the classes of
ordinary Azumaya algebras over $Y$. On the other hand, by Corollary~\ref{cor:azumayagluing},
they are represented by (derived) Azumaya algebras over $Y$.

\vspace{20pt}
\noindent
Benjamin Antieau
[\texttt{\href{mailto:antieau@math.ucla.edu}{antieau@math.ucla.edu}}]\\
UCLA\\
Math Department\\
520 Portola Plaza\\
Los Angeles, CA 90095-1555\\
USA\\

\noindent
David Gepner 
[\texttt{\href{mailto:djgepner@gmail.com}{djgepner@gmail.com}}]\\
Fakult\"at f\"ur Mathematik\\
Universit\"at Regensburg\\
93040 Regensburg\\
Germany\\

\end{document}